\documentclass[12pt,a4paper,draft]{article}%
\usepackage{amsfonts}
\usepackage{makeidx}
\usepackage{latexsym}
\usepackage{amsmath}
\usepackage{color}
\usepackage{extarrows}
\usepackage{amssymb}
\usepackage{graphicx, amsmath, amsthm, amssymb, bm, boxedminipage}
\usepackage{graphicx}%
\providecommand{\U}[1]{\protect\rule{.1in}{.1in}}
\allowdisplaybreaks[4]
\allowdisplaybreaks[4]

\numberwithin{equation}{section}
\theoremstyle{plain}
\newtheorem{thm}{Theorem}[section]
\newtheorem{lemma}[thm]{Lemma}
\newtheorem{remark}[thm]{Remark}

\newtheorem{corol}[thm]{Corollary}

\begin{document}

\title{A Central Limit Theorem for Sets of Probability Measures\thanks{Chen is at
School of Mathematics, Shandong University, zjchen@sdu.edu.cn, and Epstein is
at Department of Economics, Boston University, lepstein@bu.edu. We thank Juan
Li and Shige Peng for helpful comments. Chen gratefully acknowledges the
support of the National Key R\&D Program of China (grant No. 2018YFA0703900).
}}
\author{Zengjing Chen
\and Larry G. Epstein}
\maketitle
\date{}

\begin{abstract}
We prove a central limit theorem for a sequence of random variables whose
means are ambiguous and vary in an unstructured way. Their joint distribution
is described by a set of measures. The limit is (not the normal distribution
and is) defined by a backward stochastic differential equation that can be
interpreted as modeling an ambiguous continuous-time random walk.

\end{abstract}


\pagestyle{plain}


\section{Introduction}

We present a Central Limit Theorem (CLT) for situations where random events
(or experiments) are describable by nonsingleton sets of probability measures.
Such sets arise in economics and finance as the subjective prior beliefs of an
agent within a model who does not have sufficient information to justify
reliance on a single probability measure\ (e.g. \cite{gs89, gilboabook,
annrev}), in mathematical statistics and econometrics where, for example, they
represent the predictions of the theory being tested or estimated empirically
and where predictions are multivalued because the theory is incomplete (e.g.
\cite{huber, walley, tamer, eks}). We refer to such situations as featuring
ambiguity. Our focus in this paper is on a sequential or temporal context,
where experiments are ordered. The set of probability measures can be taken to
be objective (the set of logically possible probability laws) or subjective
(representing an individual's beliefs about future experiments).

Our first main result can be outlined roughly as follows. Let $\left(
\Omega,\mathcal{G}\right)  $ be a measurable space and let $\left(
X_{i}\right)  $ be a sequence of (real-valued) random variables, where $X_{i}$
describes the outcome of experiment $i$. Let $\mathcal{P}$ be a set of
probability measures on $\left(  \Omega,\mathcal{G}\right)  $. Information is
represented by the filtration $\{\mathcal{G}_{i}\}${, (}${\mathcal{G}%
_{0}=\{\emptyset,\Omega\}}$), such that $\left(  X_{i}\right)  $ is adapted to
$\{\mathcal{G}_{i}\}$ and $\mathcal{G}=\sigma(\cup_{1}^{\infty}\mathcal{G}%
_{i})$. Assume that the upper and lower conditional means of the $X_{i}$s
satisfy:
\begin{equation}
ess\sup\limits_{Q\in\mathcal{P}}E_{Q}[X_{i}|\mathcal{G}_{i-1}]=\overline{\mu
}\text{ and }ess\inf_{Q\in\mathcal{P}}E_{Q}[X_{i}|\mathcal{G}_{i-1}%
]=\underline{\mu}\text{, for all }i\geq1\text{.} \label{mubar}%
\end{equation}
Ambiguity about means is indicated if $\overline{\mu}>$\underline{$\mu$}.
Conditional variances are taken to be unambiguous and common to all $X_{i}$s:%
\begin{equation}
E_{Q}\left[  (X_{i}-E_{Q}[X_{i}|\mathcal{G}_{i-1}])^{2}|\mathcal{G}%
_{i-1}\right]  =\sigma^{2}>0\text{ for all }Q\in\mathcal{P}\text{ and all
}i\text{.} \label{condvar}%
\end{equation}
Then, under suitable additional assumptions, we show that for every
$\varphi\in C\left(  \left[  -\infty,\infty\right]  \right)  $, the class of
all bounded continuous functions with finite limits at $\pm\infty$,{\small
\begin{equation}
\lim\limits_{n\rightarrow\infty}\sup\limits_{Q\in\mathcal{P}}E_{Q}\left[
\varphi\left(  \frac{1}{n}{\sum_{i=1}^{n}X_{i}}+\frac{1}{\sqrt{n}}{\sum
_{i=1}^{n}}\frac{1}{\sigma}{(X_{i}-E_{Q}[X_{i}|\mathcal{G}_{i-1}])}\right)
\right]  =\mathbb{E}_{\left[  \underline{\mu},\overline{\mu}\right]  }\left[
\varphi\left(  B_{1}\right)  \right]  \text{,} \label{CLT0}%
\end{equation}
}where the right side of this equation is defined to be $Y_{0}$, given that
$(Y_{t},Z_{t})$ is the solution of the backward stochastic differential
equation (BSDE)
 {\normalsize
\begin{equation}
Y_{t}=\varphi\left(  B_{1}\right)  +\int_{t}^{1}{\max_{\underline{\mu}\leq
\mu\leq\overline{\mu}}(\mu Z_{s})}ds-\int_{t}^{1}Z_{s}dB_{s},\;0\leq t\leq1,
\label{BSDE}%
\end{equation}
}and $(B_{t})$ is a standard Brownian motion on a probability space
$(\Omega^{\ast},\mathcal{F}^{\ast},P^{\ast})$.

{\normalsize The result highlights the connection between CLTs and BSDEs. If
$\overline{\mu}=$\underline{$\mu$}$=\mu$, and given (\ref{mubar}) and any
fixed measure in $\mathcal{P}$, then $\left(  X_{i}-\mu\right)  $ is a
martingale difference and the limit result reduces to a form of the classical
martingale CLT (applying the strong {Law of Large Numbers (LLN)} for
martingales which gives a.s. convergence of $\frac{1}{n}\sum_{i=1}^{n}X_{i}$
to $\mu$). In addition, the right side reduces to a linear BSDE that, through
its solution, yields the expectation of $\varphi\left(  B_{1}\right)  $ under
the normal distribution $\boldsymbol{N}\left(  \mu,1\right)  $. More
generally, in our CLT accommodating ambiguity about means, the associated BSDE
is nonlinear. Rather it corresponds to a model in which a Brownian motion is
augmented by a drift that can vary stochastically thru time subject only to
remaining in the interval $\left[  \underline{\mu},\overline{\mu}\right]  $.
For example, when $\varphi$ is the indicator function $I_{\left[  a,b\right]
}$, \cite{chenBSDE} shows that $sgn\left(  Z_{s}\right)  =-sgn\left(
B_{s}-\frac{a+b-\left(  \overline{\mu}+\underline{\mu}\right)  \left(
1-s\right)  }{2}\right)  $, and hence%
\[
{\arg\max_{\underline{\mu}\leq\mu\leq\overline{\mu}}(\mu Z_{s})}=\left\{
\begin{array}
[c]{cc}%
\underline{\mu} & \text{ if }B_{s}\geq\frac{a+b-\left(  \overline{\mu
}+\underline{\mu}\right)  \left(  1-s\right)  }{2},\\
\overline{\mu} & \text{ \ if }B_{s}<\frac{a+b-\left(  \overline{\mu
}+\underline{\mu}\right)  \left(  1-s\right)  }{2}.
\end{array}
\right.
\]
This stochastic variability of the maximizing mean $\mu$ leads to a non-normal
limiting distribution. }

{\normalsize Two important points regarding tractability should be noted.
First, from \cite{chenBSDE} and also Lemma \ref{explicit-solution} below, the
indicated BSDE can be solved in closed-form for some specifications of
$\varphi$. For example, when $\varphi$ is the indicator for the interval
$\left[  a,b\right]  $, then the right side of (\ref{CLT0}) is given by
\begin{equation}
{\mathbb{E}_{\left[  \underline{\mu},\overline{\mu}\right]  }[I_{[a,b]}\left(
B_{1}\right)  ]=\left\{
\begin{array}
[c]{lc}%
\Phi_{-\overline{\mu}}\left(  -a\right)  -e^{-\frac{(\overline{\mu}%
-\underline{\mu})(b-a)}{2}}\Phi_{-\overline{\mu}}\left(  -b\right)  & \text{if
}a+b\geq d,\\
\Phi_{\underline{\mu}}\left(  b\right)  -e^{-\frac{(\overline{\mu}%
-\underline{\mu})(b-a)}{2}}\Phi_{\underline{\mu}}\left(  a\right)  & \text{if
}a+b<d,
\end{array}
\right.  \label{upperBSDE}}%
\end{equation}
where $d\equiv\overline{\mu}+\underline{\mu}$ and {$\Phi_{\mu}$ is the normal
cdf with mean $\mu$ and unit variance. The second point concerns the left side
of (\ref{CLT0}) which is nonstandard in that the argument of }$\varphi$, whose
distribution is at issue, includes the measures $Q$ in $\mathcal{P}$ and hence
is not a function only of past realizations of the $X_{i}$s. However, our
second principal result (Theorem \ref{thm-special}) is that for a class of
functions $\varphi$, including indicators and quadratics, both of which are
prominent in statistical theory and methods, (\ref{CLT0}) is valid also when
each conditional expectation ${E_{Q}[X_{i}|\mathcal{G}_{i-1}]}$ is replaced by
a suitable (and explicit) function of $\left(  X_{1},...,X_{i-1}\right)  $
alone. Potential usefulness of this result is illustrated by an application to
hypothesis testing. }

{\normalsize The key additional assumption underlying both theorems is that
the set $\mathcal{P}$ is "rectangular", or closed with respect to the pasting
of alien marginals and conditionals. (Rectangularity was introduced in
\cite{es2003} in the context of recursive utility theory, where an axiomatic
analysis demonstrated its central role in modeling dynamic behavior. It has
been studied and applied also in robust stochastic dynamic optimization
\cite{shapiro}, in the literature on dynamic risk measures
\cite{reidel,cheridito,penner}, and in continuous-time modeling in finance
\cite{CE}.) It can be understood as endowing $\mathcal{P}$ with a recursive
structure that yields a form of the law of iterated expectations. If
$\mathcal{P}=\{P\}$, which implies (and, for our purposes, is essentially
equivalent to) $\overline{\mu}=$\underline{$\mu$}, then the law of iterated
expectations is a consequence of updating by Bayes rule and rectangularity is
vacuously satisfied. (Sections \ref{section-prelim} and \ref{section-IID}
provide a precise definition of rectangularity and some motivating informal
interpretation.) }

{\normalsize Some connections to the literature conclude this introduction. In
the classical probability framework, there are numerous CLTs with non-normal
limiting distributions (with stable laws, for example) \cite{KG,dasgupta}, all
of which have much different motivation and limits than our result. There
exist alternative generalizations of the classical theorem that are motivated
by robustness to ambiguity.\textbf{ }In \cite{eks} (see also the
generalization in \cite{shi}), experiments are not ordered and the analysis is
intended for a cross-sectional context. In addition, $\mathcal{P}$ is assumed
to be the core of a convex (that is, supermodular) capacity, which renders it
inconsistent with a recursive structure \cite{chen-davison}. Finally, the
limiting distribution in their result is the normal, in contrast to our novel
BSDE-based limit. }

{\normalsize Closer to this paper is the CLT due to Peng \cite{peng2008,
peng2009} who also assumes that experiments are ordered. Peng's focus is on
ambiguity about variance (or at least about the second moment), while{ our
focus is on ambiguity about the mean.}
A more recent paper \cite{peng2019} provides a CLT (Theorem 3.2) with
ambiguity about both mean and variance. (Their theorem also considers rates of
convergence, which are ignored here.) To compare it with this paper, consider
the special case where there is ambiguity about means only. Then their CLT is
related primarily to our Theorem \ref{thm-CLT2}, rather than to our central
results Theorems \ref{thm-CLT} and \ref{thm-special}. In particular, only in
the latter are limits defined by a BSDE rather than by a normal distribution
(as in \cite{peng2019}). See section \ref{section-discuss} for elaboration.
Another difference is that our approach is more probability-theoretic: Peng
and coauthors take a nonlinear expectation as the core primitive and adopt the
PDE approach, while our primitive is a set of probability measures and
conditionals are central only in our analysis.}

{\normalsize The next section describes the model's primitives and key
assumptions formally. These are illustrated in section \ref{section-IID} via a
canonical example that can be understood as generalizing the classical random
walk to accommodate ambiguity. The two main CLT results (Theorems
\ref{thm-CLT} and \ref{thm-special}) are presented in section
\ref{section-CLT}. Section \ref{section-discuss} provides perspective on our
main results by relating them to an alternative CLT (Theorem \ref{thm-CLT2})
and a weak LLN for our setting (Corollary \ref{cor-LLN}). Proofs of Theorems
\ref{thm-CLT} and \ref{thm-special} are presented in section
\ref{section-proofs}. An appendix contains other proofs and supplementary
material. }

\section{Primitives and assumptions\label{section-prelim}}

{\normalsize Let $\left(  \Pi_{1}^{\infty}\Omega_{i},\{\mathcal{G}_{n}%
\}_{n=1}^{\infty}\right)  $ be a filtered space modeling a sequence of
experiments. The set of possible outcomes for the $i^{th}$ experiment is
$\Omega_{i}$. For each $n$, $\mathcal{G}_{n}$ is a $\sigma$-algebra on
$\Pi_{1}^{n}\Omega_{i}$ representing the observable events regarding
experiments $1,...,n$. (Accordingly, we assume that $\mathcal{G}_{n}$ is
increasing with $n$ and we take $\mathcal{G}_{0}$ to be the trivial $\sigma
$-algebra.) The observable events for the collection of all experiments are
given by $\mathcal{G}$,
\[
\mathcal{G}=\sigma(\cup_{1}^{\infty}\mathcal{G}_{n})\text{,}%
\]
a $\sigma$-algebra on $\Omega$, where
\[
\Omega=\Pi_{1}^{\infty}\Omega_{i}\text{.}%
\]
(Here and in the sequel, we identify each $\mathcal{G}_{n}$ in the obvious way
with a }$\sigma$-algebra on $\Omega$.) {\normalsize The ex ante probabilities
of observable events are not known precisely and are represented by a set
$\mathcal{P}$ of probability measures,\footnote{For any measurable space
$\left(  Y,\mathcal{F}\right)  $, the corresponding set of probability
measures is denoted $\Delta\left(  Y,\mathcal{F}\right)  $.}
\[
\mathcal{P}\subset\Delta\left(  \Omega,\mathcal{G}\right)  \text{.}%
\]
We limit ambiguity about which events are possible and assume that all
measures in $\mathcal{P}$ are equivalent on each $\mathcal{G}_{n}$. }

{\normalsize Below we assume that for each measure $P$ in $\mathcal{P}$ and
each $n$, there exists a regular conditional measure $P\left(  \cdot
\mid\mathcal{G}_{n}\right)  $. For example, a well-known \cite[Theorem
7.1]{partha} sufficient condition for regular $\mathcal{G}_{n}$-conditionals
to exist for every $P$ in $\Delta\left(  \Omega,\mathcal{G}\right)  $ is that
$\left(  \Omega,\mathcal{G}\right)  $ is a separable standard Borel space (a
special case is where $\Omega$ is a complete separable metric space and
$\mathcal{G}$ is its Borel $\sigma$-algebra). }

{\normalsize Finally, we consider a sequence $\left(  X_{i}\right)  $ of
real-valued random variables (r.v.), $X_{i}:$ $\Pi_{1}^{\infty}\Omega
_{j}\longrightarrow\mathbb{R}$, such that $X_{i}$ is $\mathcal{G}_{i}%
$-measurable (using the Borel $\sigma$-algebra on $\mathbb{R}$). Think of
$X_{i}$ as a scalar measure of the outcome of experiment $i$ or of the value
(or utility) of that outcome. In general, $X_{i}$ can depend also on the
outcomes of earlier experiments. }

\begin{remark}
{\normalsize We presume a particular ordering of experiments, which may be
arbitrary in cross-sectional contexts. Thus we view the analysis and the
resulting CLT as more relevant to sequential or time-series contexts where an
ordering is given. }
\end{remark}

{\normalsize In the rest of this section, we describe our assumptions on the
above primitives. We use the following notation. $\mathcal{H}$ denotes the set
of all r.v. $X$ on $\left(  \Omega,\mathcal{G}\right)  $ satisfying
$\sup_{Q\in\mathcal{P}}E_{Q}[|X|]<\infty$; $E_{Q}[\cdot]$ is the expectation
under the probability measure $Q.$ For any $X$ in $\mathcal{H}$, its
\emph{upper and lower expectations} are defined respectively by
\[
\mathbb{E}[X]\equiv\sup_{Q\in\mathcal{P}}E_{Q}[X],\;\quad\mathcal{E}%
[X]\equiv\inf_{Q\in\mathcal{P}}E_{Q}[X]=-\mathbb{E}[-X]\text{,}%
\]
and its \emph{conditional upper and lower expectations} are defined
respectively by%
\[
\mathbb{E}\left[  X\mid\mathcal{G}_{n}\right]  \equiv ess\sup_{Q\in
\mathcal{P}}E_{Q}\left[  X\mid\mathcal{G}_{n}\right]  \text{, \ \ }%
\mathcal{E}\left[  X\mid\mathcal{G}_{n}\right]  \equiv ess\inf_{Q\in
\mathcal{P}}E_{Q}\left[  X\mid\mathcal{G}_{n}\right]  \text{.}%
\]
Obviously, the conditional expectations are well-defined due to equivalence of
all measures in $\mathcal{P}$ on each $\mathcal{G}_{n}.$ (See section
\ref{section-proofs}{ for key properties of these expectations.) Rewritten
with this notation, (\ref{mubar}) takes the form}%
\[
\mathbb{E}\left[  X_{i}\mid\mathcal{G}_{i-1}\right]  =\overline{\mu}\text{ and
}\mathcal{E}\left[  X_{i}\mid\mathcal{G}_{i-1}\right]  =\underline{\mu}\text{
for all }i\text{.}%
\]
{ } }

{\normalsize Say that $(X_{i})$ has an \emph{unambiguous conditional variance
}$\sigma^{2}\,$\ if (\ref{condvar}) is satisfied. Say that $\left(
X_{i}\right)  $ satisfies the \emph{Lindeberg condition} if
\begin{equation}
\lim_{n\rightarrow\infty}\frac{1}{n}\sum\limits_{i=1}^{n}\mathbb{E}\left[
\left\vert X_{i}\right\vert ^{2}I_{\{\left\vert X_{i}\right\vert >\sqrt
{n}\epsilon\}}\right]  =0,\quad\forall\epsilon>0. \label{linder}%
\end{equation}
}

{\normalsize To formulate the remaining assumption requires additional
notation and terminology. Write
\begin{align*}
\omega_{\left(  n\right)  }  &  =\left(  \omega_{n},...\right)  \text{,
}\omega^{\left(  n\right)  }=\left(  \omega_{1},...,\omega_{n}\right)  ,\\
\mathcal{P}_{0,n}  &  =\{P_{\mid\mathcal{G}_{n}}:P\in\mathcal{P}\}\text{
\ and}\\
\mathcal{G}_{(n+1)}  &  =\left\{  A\subset\Pi_{n+1}^{\infty}\Omega_{i}:\Pi
_{1}^{n}\Omega_{i}\times A\in\mathcal{G}\right\}  \text{.}%
\end{align*}
}

{\normalsize A probability kernel from $\left(  \Pi_{1}^{n}\Omega
_{i},\mathcal{G}_{n}\right)  $ to $\left(  \Pi_{n+1}^{\infty}\Omega
_{i},\mathcal{G}_{(n+1)}\right)  $ is a function $\lambda:\Pi_{1}^{n}%
\Omega_{i}\times\mathcal{G}_{(n+1)}\longrightarrow\left[  0,1\right]  $
satisfying: }

\begin{enumerate}
\item[Kernel1] {\normalsize $\forall\omega^{(n)}\in\Pi_{1}^{n}\Omega_{i}$,
$\lambda\left(  \omega^{\left(  n\right)  },\cdot\right)  $ is a probability
measure on $\left(  \Pi_{n+1}^{\infty}\Omega_{i},\mathcal{G}_{(n+1)}\right)
$, }

\item[Kernel2] {\normalsize $\forall A\in\mathcal{G}_{(n+1)}$, $\lambda\left(
\cdot,A\right)  $ is a $\mathcal{G}_{n}$-measurable function on $\Pi_{1}%
^{n}\Omega_{i}$. }
\end{enumerate}

{\normalsize \noindent Any pair $\left(  p_{n},\lambda\right)  $ consisting of
a probability measure $p_{n}$ on $\left(  \Pi_{1}^{n}\Omega_{i},\mathcal{G}%
_{n}\right)  $ and a probability kernel $\lambda$ as above, induces a unique
probability measure $P$ on $\left(  \Pi_{1}^{\infty}\Omega_{i},\mathcal{G}%
\right)  $ that coincides with $p_{n}$ on $\mathcal{G}_{n}$. It is given by,
$\forall A\in\mathcal{G}$,%
\begin{equation}
P\left(  A\right)  =\int\limits_{\Pi_{1}^{n}\Omega_{i}}\int\limits_{\Pi
_{n+1}^{\infty}\Omega_{i}}I_{A}\left(  \omega^{\left(  n\right)  }%
,\omega_{(n+1)}\right)  \lambda\left(  \omega^{\left(  n\right)  }%
,d\omega_{(n+1)}\right)  p_{n}\left(  d\omega^{\left(  n\right)  }\right)
\text{.} \label{Pn}%
\end{equation}
}

{\normalsize For $Q\in\mathcal{P}$, let $Q\left(  \cdot\mid\mathcal{G}%
_{n}\right)  $ denote its regular conditional. Then it defines a probability
kernel $\lambda$ by: $\forall\omega^{\left(  n\right)  }\in\Pi_{1}^{n}%
\Omega_{i}$,
\begin{equation}
\lambda\left(  \omega^{\left(  n\right)  },A\right)  =Q\left(  \Pi_{1}%
^{n}\Omega_{i}\times A\mid\mathcal{G}_{n}\right)  \left(  \omega^{\left(
n\right)  }\right)  \text{, }\forall A\in\mathcal{G}_{(n+1)}\text{. ~}
\label{kernel}%
\end{equation}
A feature of such a kernel is that the single measure $Q$ is used to define
the conditional at every $\omega^{\left(  n\right)  }$. We are interested in
kernels for which the measure to be conditioned can vary with $\omega^{\left(
n\right)  }$. Say that the probability kernel $\lambda$ is a $\mathcal{P}%
$-\emph{kernel} if: $\forall\omega^{\left(  n\right)  }\in\Pi_{1}^{n}%
\Omega_{i}$, $\exists Q\in\mathcal{P}$ satisfying (\ref{kernel}). }

{\normalsize Finally, say that $\mathcal{P}$ is \emph{rectangular }(with
respect to the filtration $\{\mathcal{G}_{n}\}$) if: $\forall n,~\forall
p_{n}\in$ $\mathcal{P}_{0,n}$ and for every $\mathcal{P}$-\emph{kernel
}$\lambda$, if $P$ is defined as in (\ref{Pn}), then $P\in\mathcal{P}$. (Note
that a measure $P\in\Delta\left(  \Omega,\mathcal{G}\right)  $ is well-defined
by (\ref{Pn}), for any $p_{n}\in$ $\mathcal{P}_{0,n}$ and $\mathcal{P}%
$-\emph{kernel }$\lambda$, because of the assumption that all measures in
$\mathcal{P}$ are equivalent on $\mathcal{G}_{n}$). When $\mathcal{P}$ is the
singleton $\{P\}$, rectangularity is trivially implied by Bayesian updating,
specifically by the fact that after decomposing $P$ into a marginal and
conditional, these can be pasted together to recover $P$. More generally,
rectangularity requires that the set $\mathcal{P}$ is closed also with respect
to pasting together conditionals and marginals that are \emph{alien}, that is,
induced by possibly different measures in $\mathcal{P}$. In this sense,
$\mathcal{P}$ does not restrict the pattern of heterogeneity across
experiments (see the next section for elaboration). }

{\normalsize The significance of rectangularity is illuminated by the
following lemma. (See Appendix \ref{app-rect} for a partial proof. The
complement of any $A\subset\Omega$ is denoted }$A^{c}$.{\normalsize ) }

\begin{lemma}
{\normalsize \label{lemma-rect} $\mathcal{P}$ rectangular implies the
following (for any $0\leq m\leq n\in\mathbb{N}$).}

{\normalsize \noindent(i) \textbf{Stability by composition:} For any
$Q,R\in\mathcal{P}$, $\exists P\in\mathcal{P}$ such that, for any
$X\in\mathcal{H}$,
\[
E_{P}[X|\mathcal{G}_{m}]=E_{Q}[E_{R}[X|\mathcal{G}_{n}]|\mathcal{G}%
_{m}]\text{.}%
\]
(ii) \textbf{Stability by bifurcation:} For any $Q,R\in\mathcal{P}$, and any
$A_{n}\in\mathcal{G}_{n}$, $\exists P\in\mathcal{P}$ such that, for any
$X\in\mathcal{H}$,
\[
E_{P}[X|\mathcal{G}_{n}]=I_{A_{n}}E_{Q}[X|\mathcal{G}_{n}]+I_{A_{n}^{c}}%
E_{R}[X|\mathcal{G}_{n}].
\]
}

{\normalsize \noindent(iii) \textbf{Law of iterated upper expectations}: For
any $X\in\mathcal{H}$,
\begin{equation}
\mathbb{E}[\mathbb{E}[X|\mathcal{G}_{n}]|\mathcal{G}_{m}]=\mathbb{E}%
[X|\mathcal{G}_{m}]. \label{LIE}%
\end{equation}
(iv) Let $\{X_{i}\}$ be a sequence in $\mathcal{H}$. {Set $S_{n-1}=\sum
_{i=1}^{n-1}X_{i}$ and, for any $Q\in\mathcal{P}$},{ $S_{n-1}^{Q}=S_{n-1}%
-\sum_{i=1}^{n-1}E_{Q}[X_{i}|\mathcal{G}_{i-1}]$. Then, for any continuous
bounded functions $f,h$:
\begin{align}
&  \sup\limits_{Q\in\mathcal{P}}E_{Q}\left[  f\left(  \frac{S_{n-1}}{n}%
+\frac{S_{n-1}^{Q}}{\sqrt{n}}\right)  +h\left(  \frac{S_{n-1}}{n}%
+\frac{S_{n-1}^{Q}}{\sqrt{n}}\right)  X_{n}\right] \nonumber\\
=  &  \sup\limits_{Q\in\mathcal{P}}E_{Q}\left[  ess\sup\limits_{R\in
\mathcal{P}}E_{R}\left[  f\left(  \frac{S_{n-1}}{n}+\frac{S_{n-1}^{Q}}%
{\sqrt{n}}\right)  +h\left(  \frac{S_{n-1}}{n}+\frac{S_{n-1}^{Q}}{\sqrt{n}%
}\right)  X_{n}|\mathcal{G}_{n-1}\right]  \right]  .\nonumber
\end{align}
}(v) If $\{X_{i}\}$ is a sequence in $\mathcal{H}$ satisfying (\ref{mubar}),
then
\begin{equation}
\mathbb{E}\left[  X_{n}\mid\mathcal{G}_{n-1}\right]  =\mathbb{E}\left[
X_{n}\right]  \text{ and }\mathcal{E}\left[  X_{n}\mid\mathcal{G}%
_{n-1}\right]  =\mathcal{E}\left[  X_{n}\right]  \text{. } \label{condindep}%
\end{equation}
}
\end{lemma}

{\normalsize (i) and (ii) make explicit two senses in which rectangularity of
$\mathcal{P}$ implies that combinations of distinct measures ($Q\not =R$) from
$\mathcal{P}$ leave one within $\mathcal{P}$. Together they lead to (iii). The
latter is built into the classical model but must be adopted explicitly, via
rectangularity, for upper (or lower) expectations. For a general set
$\mathcal{P}$, one would expect the supremum on the left in (\ref{LIE}) to be
(weakly) larger because it permits the choices of measures conditional on each
history {$\omega^{\left(  n-1\right)  }$ } and the ex ante measure on
$\mathcal{G}_{n-1}$ to be alien. However, rectangularity implies that any such
combination of measures yields a measure in $\mathcal{P}$, and thus the
single-stage supremum on the right is no smaller. The proof of our CLT employs
a similar recursive relation also in instances when the r.v. itself depends on
$Q\in\mathcal{P}$ as in (iv), the intuition for which is similar to that for
(iii). (v) states that conditional upper and lower expectations do not vary
with the outcomes of previous experiments. It is an immediate consequence of
(\ref{mubar}) and (iii); for example,
\[
\mathbb{E}\left[  X_{n}\right]  =\mathbb{E}\left[  \mathbb{E}\left[  X_{n}%
\mid\mathcal{G}_{n-1}\right]  \right]  =\mathbb{E}\left[  \overline
{\mathbb{\mu}}\right]  =\overline{\mathbb{\mu}}=\mathbb{E}\left[  X_{n}%
\mid\mathcal{G}_{n-1}\right]  \text{.}%
\]
}

\section{Example: IID\label{section-IID}}

{\normalsize Our canonical example (adapted from \cite{esIID}) is as follows.
Specialize the above framework by assuming that there exists a measurable
space $(\overline{\Omega},\overline{\mathcal{F}})$ such that, for all {$1\leq
i\leq n$,}%
\[
(\Omega_{i},\mathcal{F}_{i})=(\overline{\Omega},\overline{\mathcal{F}})\text{
and }\mathcal{G}_{n}=\Pi_{1}^{n}\mathcal{F}_{i}\text{.}%
\]
That is, experiments have a common set of possible outcomes $\overline{\Omega
}$ and an associated common $\sigma$-algebra $\overline{\mathcal{F}}$. In
addition, suppose that, for all $i$,
\[
X_{i}=\overline{X}:(\Omega_{i},\mathcal{F}_{i})\rightarrow\mathbb{R}\text{.}%
\]
}

{\normalsize One-step-ahead conditionals are central. Thus, for each $P$ in
$\mathcal{P}$, and each $n$, let $P_{n,n+1}\left(  \omega^{\left(  n\right)
}\right)  $ denote the restriction to $\mathcal{G}_{n+1}$ of $P\left(
\cdot\mid\mathcal{G}_{n}\right)  \left(  \omega^{\left(  n\right)  }\right)
$. }

{\normalsize Fix a subset $\mathcal{L}$ of $\Delta(\overline{\Omega}%
,\overline{\mathcal{F}})$, all of whose measures are equivalent. Then the IID
model is defined via the set $\mathcal{P}^{IID}$,%
\begin{equation}
\mathcal{P}^{IID}=\left\{  P\in\Delta\left(  \Omega,\mathcal{G}\right)
:P_{n,n+1}\left(  \omega^{\left(  n\right)  }\right)  \in\mathcal{L},\text{
}\forall n,\omega^{\left(  n\right)  }\in\mathcal{G}_{n}\right\}  \text{.}
\label{Pca}%
\end{equation}
The set consists of all measures whose one-step-ahead conditionals, at every
history, lie in $\mathcal{L}$. Thus, $\mathcal{L}$ is the set of plausible
probability laws for each experiment, independent of history, modeling partial
ignorance about each experiment separately. There remains the question of the
perception of, or information about, the sequence of experiments, that is, how
experiments are related to one another. In spite of $\mathcal{L}$ being common
to all $i$, in this model experiments are not necessarily identical.
(Accordingly, we refer to experiments as being \textit{indistinguishable
rather than identical} and take IID to mean "\textit{indistinguishably} and
independently distributed".) Indeed, any measure in $\mathcal{L}$ is plausible
as the law describing the {$i^{th}$ experiment} in conjunction with any
possibly different measure in $\mathcal{L}$ being the law describing the
{$j^{th}$ experiment}. Indeed, $\mathcal{P}^{IID}$ imposes no restrictions on
joint distributions thus capturing \emph{agnosticism about the pattern of
heterogeneity across experiments}. As demonstrated below, this feature is
closely related to rectangularity. }

{\normalsize In the special case where $\mathcal{L}=\{P\}$, $\mathcal{P}%
^{IID}$ consists of the single i.i.d. product of $P$, as in a random walk. One
might think of $\mathcal{P}^{IID}$ as modeling an ``ambiguous random walk". }

{\normalsize The following lemma gives some readily verified properties of
$\mathcal{P}^{IID}$ (see Appendix \ref{app-IID} for some proof details). }

\begin{lemma}
{\normalsize \label{lemma-IID}The set $\mathcal{P}^{IID}$ satisfies (for every
$n\in\mathbb{N}$): }\newline{\normalsize (i) $\mathcal{P}^{IID}$ is
rectangular.\newline(ii) Measures in $\mathcal{P}^{IID}$ are mutually
equivalent on each $\mathcal{G}_{n}$.\newline(iii) For any $\varphi\in
C\left(  \mathbb{R}\right)  $, with $\varphi\left(  X_{n}\right)
\in\mathcal{H}$,
\[
\mathbb{E}\left[  \varphi\left(  X_{n}\right)  \mid\mathcal{G}_{n-1}\right]
=\sup_{q\in\mathcal{L}}E_{q}\left[  \varphi\left(  \overline{X}\right)
\right]  =\mathbb{E}\left[  \varphi(X_{n})\right]  =\mathbb{E}\left[
\varphi(X_{1})\right]  .
\]
(iv) Conditional variances satisfy:
\begin{align*}
\sup_{Q\in\mathcal{P}^{IID}}E_{Q}\left[  (X_{n}-E_{Q}[X_{n}|\mathcal{G}%
_{n-1}])^{2}|\mathcal{G}_{n-1}\right]   &  =\sup_{q\in\mathcal{L}}E_{q}\left[
(\overline{X}-E_{q}[\overline{X}])^{2}\right]  ,\\
\inf_{Q\in\mathcal{P}^{IID}}E_{Q}\left[  (X_{n}-E_{Q}[X_{n}|\mathcal{G}%
_{n-1}])^{2}|\mathcal{G}_{n-1}\right]   &  =\inf_{q\in\mathcal{L}}E_{q}\left[
(\overline{X}-E_{q}[\overline{X}])^{2}\right]  .
\end{align*}
}
\end{lemma}

{\normalsize The key property of $\mathcal{P}^{IID}$ is rectangularity.
Because of its centrality, we verify rectangularity here: Let $p_{n}$,
$\lambda$ and $P$ be as in (\ref{kernel}). Then, for the given $\omega
^{\left(  n\right)  }$,
\begin{align*}
P\left(  \Pi_{1}^{n}\Omega_{i}\times\cdot\mid\mathcal{G}_{n}\right)  \left(
\omega^{\left(  n\right)  }\right)   &  =\lambda\left(  \omega^{\left(
n\right)  },\cdot\right) \\
&  =Q\left(  \Pi_{1}^{n}\Omega_{i}\times\cdot\mid\mathcal{G}_{n}\right)
\left(  \omega^{\left(  n\right)  }\right)  ,
\end{align*}
for some $Q\in\mathcal{P}^{IID}$. Therefore, the one-step-ahead conditional of
$P$ at history $\omega^{\left(  n\right)  }$ equals that of $Q$ and hence lies
in $\mathcal{L}$. Therefore, $P\in\mathcal{P}^{IID}$.\hfill\hfill\ \hfill}

{\normalsize The lemma implies that $\mathcal{P}^{IID}$ readily accommodates
also the other assumptions in the CLT below. For example, (\ref{mubar}) is
implied by (iii) and conditional variances are common and unambiguous if and
only if\footnote{In decision theory \noindent(in \cite{gs89}, for example), it
is often innocuous and a convenient normalization to take sets of measures to
be convex. But because variances are not linear in the measure $q$, convexity
of $\mathcal{L}$ precludes (\ref{varL}) except in the degenerate case where
means are also unambiguous. Thus we do not assume that $\mathcal{L}$ is
convex.}%
\begin{equation}
var_{q}\left(  \overline{X}\right)  \equiv E_{q}\left[  (\overline{X}%
-E_{q}[\overline{X}])^{2}\right]  =\sigma^{2}\text{, \ for all }%
q\in\mathcal{L}\text{. } \label{varL}%
\end{equation}
}

{\normalsize For perspective, consider also the set $\mathcal{P}^{prod}$,
consisting of all (nonidentical) product measures that can be constructed from
$\mathcal{L}$ - refer to this as the \emph{product model}. The set
$\mathcal{P}^{prod}$ also implies a degree of agnosticism about heterogeneity
-- after all, it consists of product measures $\Pi_{i=1}^{\infty}\ell_{i}$,
where $\ell_{i}\not =\ell_{j}$ in general, and these measures are restricted
only by the requirement that they lie in $\mathcal{L}$. However, the two
models differ in a significant way in that $\mathcal{P}^{prod}$ violates
rectangularity, and hence also (\ref{LIE}), for example. This is because
$\mathcal{P}^{prod}$ is "too small" in the sense of not being closed with
respect to the pasting of alien marginals and conditionals (note that
$\mathcal{P}^{prod}$ is a strict subset of $\mathcal{P}^{IID}$). Our
interpretation of $\mathcal{P}^{prod}$ is that it models \emph{certainty} that
the probability law for experiment $i$ does not vary with the outcomes of
preceding experiments (note that invariance to these outcomes is exhibited by
each individual measure in $\mathcal{P}^{prod}$). In contrast, in
$\mathcal{P}^{IID}$ one-step-ahead conditionals can vary arbitrarily across
different histories subject only to lying in $\mathcal{L}$. Thus
$\mathcal{P}^{IID}$ permits heterogeneity across experiments to vary
stochastically and thereby models greater agnosticism regarding heterogeneity.
}

{\normalsize A simple concrete example illustrates both models and the
difference between them. Each experiment can produce one of three outcomes:
success ($s$), failure ($f$) and the neutral outcome ($n$). Thus
$\overline{\Omega}=\{s,f,n\}$ and $\overline{\mathcal{F}}$ is the power set.
Outcomes are valued by $\overline{X}$ according to%
\[
\overline{X}\left(  s\right)  =1\text{, }\overline{X}\left(  f\right)
=-1\text{, }\overline{X}\left(  n\right)  =0\text{.}%
\]
Outcomes are uncertain but their probabilities are not known precisely. Let%
\[
0<q<p\text{, \ }p+q\leq1\text{.}%
\]
}

{\normalsize It is known that, for each experiment, and regardless of the
outcomes in preceding experiments, the outcomes $s$, $f$ and $n$ (in that
order) are given \emph{either} by the favorable distribution $\left(
p,q,1-p-q\right)  $ \emph{or} by the unfavorable distribution $\left(
q,p,1-p-q\right)  $, that is,
\[
\mathcal{L}=\left\{  \left(  p,q,1-p-q\right)  ,\left(  q,p,1-p-q\right)
\right\}  \text{.}%
\]
}

{\normalsize There is no additional information provided that would justify,
for example, assigning weights (or probabilities) to these two distributions
and then using the average as the Bayesian model would require - there is
complete ignorance about which distribution applies for any given experiment.
Consequently, conditional on any history, the implied upper and lower means of
each $X_{i}$ equal $\overline{\mu}=p-q$ and \underline{$\mu$}$=-(p-q)$
respectively, and the implied conditional variance $\sigma^{2}$ of each
$X_{i}$ is unambiguous and equals $p+q-\left(  p-q\right)  ^{2}$%
\thinspace.\ Thus $p+q$ and $p-q$ parametrize risk (measured by $\sigma^{2}$)
and ambiguity (measured by $\frac{\overline{\mu}-\underline{\mu}}{2}$)
respectively in the sense that a change in $p+q$ alone changes only risk and a
change in $p-q$ alone changes only ambiguity. }

{\normalsize The final issue is the relation between experiments. Arguably,
ignorance about which probability law applies to any given experiment,
logically implies (or at least suggests) ignorance about how experiments are
related. Accordingly, $\mathcal{P}^{IID}$ does not restrict measures on the
entire sequence of experiments beyond requiring that each one-step-ahead
conditional lie in $\mathcal{L}$. In contrast, $\mathcal{P}^{prod}$ admits
only measures for which the conditional law for the $i^{th}$ experiment,
though it can be either favorable or unfavorable, is necessarily the same for
all histories of outcomes. Thus, for example, $\mathcal{P}^{prod}$ excludes
measures that specify both (1) the favorable law for experiment $i$ after a
successful outcome in $i-1$, and (2) the unfavorable law for experiment $i$
after a failure in $i-1$. }

\section{The main results\label{section-CLT}}

\subsection{Two theorems}

{\normalsize We extend (a version of) the classical martingale CLT to admit
ambiguity about means while maintaining the assumption of unambiguous
variances. Though the theorems deal with real-valued random variables,
multidimensional versions can be proven in a similar fashion and will be
reported elsewhere. }

\newpage

\begin{thm}
{\normalsize \label{thm-CLT} Let the sequence $\left(  X_{i}\right)  $ be such
that $X_{i}\in\mathcal{H}$ for each $i$, and where $\left(  X_{i}\right)  $
satisfies (\ref{mubar}) and (\ref{condvar}), with conditional upper and lower
means $\overline{\mu}$ and $\underline{\mu}$, and{ unambiguous conditional
variance} $\sigma^{2}>0$. Assume also the Lindeberg condition (\ref{linder})
and that $\mathcal{P}$ is rectangular. Then, for any $\varphi\in C\left(
\left[  -\infty,\infty\right]  \right)  $,{\small
\begin{equation}
\lim\limits_{n\rightarrow\infty}\sup_{Q\in\mathcal{P}}E_{Q}\left[
\varphi\left(  \frac{1}{n}{\sum_{i=1}^{n}X_{i}}+\frac{1}{\sqrt{n}}%
\sum\limits_{i=1}^{n}\frac{1}{\sigma}{(X_{i}-E_{Q}[X_{i}|\mathcal{G}_{i-1}%
])}\right)  \right]  =\mathbb{E}_{\left[  \underline{\mu},\overline{\mu
}\right]  }[\varphi\left(  B_{1}\right)  ], \label{CLT}%
\end{equation}
or equivalently,%
\begin{equation}
\lim\limits_{n\rightarrow\infty}\inf_{Q\in\mathcal{P}}E_{Q}\left[
\varphi\left(  \frac{1}{n}{\sum_{i=1}^{n}X_{i}}+\frac{1}{\sqrt{n}}{\sum
_{i=1}^{n}\frac{1}{\sigma}(X_{i}-E_{Q}[X_{i}|\mathcal{G}_{i-1}])}\right)
\right]  =\mathcal{E}_{\left[  \underline{\mu},\overline{\mu}\right]
}[\varphi\left(  B_{1}\right)  ], \label{CLTinf}%
\end{equation}
where $\mathbb{E}_{\left[  \underline{\mu},\overline{\mu}\right]  }%
[\varphi\left(  B_{1}\right)  ]\equiv Y_{0}$ is called $g$-expectation
by Peng in \cite{peng1997}, given that $(Y_{t},Z_{t})$ is the solution of the
BSDE%
\[
Y_{t}=\varphi(B_{1})+\int_{t}^{1}\max\limits_{\underline{\mu}\leq\mu
\leq\overline{\mu}}(\mu Z_{s})ds-\int_{t}^{1}Z_{s}dB_{s},\ 0\leq t\leq1,
\]
and $\mathcal{E}_{\left[  \underline{\mu},\overline{\mu}\right]  }%
[\varphi\left(  B_{1}\right)  ]\equiv y_{0}$, given that $(y_{t},z_{t})$ is
the solution of the BSDE
\begin{equation}
y_{t}=\varphi(B_{1})+\int_{t}^{1}\min\limits_{\underline{\mu}\leq\mu
\leq\overline{\mu}}(\mu z_{s})ds-\int_{t}^{1}z_{s}dB_{s},\ 0\leq
t\leq1\text{.} \label{BSDE-1-1}%
\end{equation}
Here $(B_{t})$ is a standard Brownian motion on a probability space
$(\Omega^{\ast},\mathcal{F}^{\ast},P^{\ast})$. }}
\end{thm}

\begin{remark}
{\normalsize By standard limiting arguments, (\ref{CLT}) can be extended to
indicator functions for intervals. Such indicators are sufficient in the
classical CLT, because of the additivity of a single probability measure. But
when dealing with sets of measures, (\ref{CLT}) is strictly stronger. Another
remark is that while in (\ref{CLT}) the second term inside $\varphi\left(
\cdot\right)  $ is normalized by the standard deviation $\sigma$, a change of
variables delivers a CLT without that normalization. (Set $\alpha=\sigma$ and
$\beta=1$ in the statement of Theorem \ref{thm-CLTgen} in the appendix.) }
\end{remark}

{\normalsize Three differences from classical results stand out. First, the
limiting distribution is not normal but rather is given by the BSDE
(\ref{BSDE}). Another notable difference is that the r.v. on the left in
(\ref{CLT}) combines the sample average, typical of LLNs, with a term that is
more typical of CLTs. Both of these features will be discussed in section
\ref{section-discuss} below. }

{\normalsize Here we consider the fact that the argument of $\varphi$ above,
whose distribution is the focus, includes measures $Q$ from $\mathcal{P}$,
which might raise concerns about tractability. To partially alleviate such
concerns, we show that (\ref{CLT}) takes on a more tractable form when
restricted to "symmetric" functions $\varphi$. Say that $\varphi
:\mathbb{R}\rightarrow\mathbb{R}$ is \emph{symmetric with center}
$c\in\mathbb{R}$ if $\varphi\left(  c-x\right)  =\varphi\left(  c+x\right)  $
for all $x\in\mathbb{R}$. Examples include indicator(s) {$\varphi(t)=\pm
I_{\left[  a,b\right]  }(t)$} with $c=\frac{a+b}{2}$, and quadratic functions
$\varphi\left(  t\right)  =\pm\left(  t-c\right)  ^{2}$, both of which are
prominent in statistical methods. It is important to emphasize also that for
both of these classes of functions \cite{chenBSDE} provides closed-form
expressions for the BSDE-based limits appearing on the right sides of
(\ref{CLT}) and (\ref{CLTinf}) above, and (\ref{sym-clt1}) and
(\ref{sym-clt2-0}) below; recall (\ref{upperBSDE}), for example. Section
\ref{section-hyp} exploits these closed-forms in an application to hypothesis
testing. }

{\normalsize The next theorem is the second major result of the paper.
(Throughout sums of the form $\Sigma_{n}^{0}x_{i}$, }$n\geq1$,{\normalsize
are taken to equal $0$, and increasing/decreasing are intended in the weak
sense.) }

{\normalsize
}

\begin{thm}
{\normalsize \label{thm-special} Adopt the assumptions in Theorem
\ref{thm-CLT} and let the function $\varphi\in C([-\infty,\infty])$ be
symmetric with center $c\in\mathbb{R}$. For $n\geq1$ and $0\leq m\leq n$,
define
\begin{align}
M_{m,n}  &  =\frac{1}{n}\sum\limits_{i=1}^{m}X_{i}+\frac{1}{\sqrt{n}}%
\sum\limits_{i=1}^{m}\frac{1}{\sigma}\left(  X_{i}-\mu_{i}^{n}\right)
\text{,}\quad M_{0,n}\equiv0,\label{Mmn}\\
\widetilde{M}_{m,n}  &  =\frac{1}{n}\sum\limits_{i=1}^{m}X_{i}+\frac{1}%
{\sqrt{n}}\sum\limits_{i=1}^{m}\frac{1}{\sigma}\left(  X_{i}-\widetilde{\mu
}_{i}^{n}\right)  , \quad\widetilde{M}_{0,n}\equiv0 \text{,}\; \label{tMmn}%
\end{align}
where }%
\begin{align}
\mu_{m}^{n}  &  =\overline{\mu}I_{A_{m-1,n}}+\underline{\mu}I_{A_{m-1,n}^{c}%
},\label{mumn}\\
A_{m-1,n}  &  =\left\{  M_{m-1,n}\leq-\tfrac{\overline{\mu}+\underline{\mu}%
}{2}\left(  {1-\tfrac{m-1}{n}}\right)  +c\right\} ,\nonumber
\end{align}
and {\normalsize
\begin{align}
\widetilde{\mu}_{m}^{n}=  &  \overline{\mu}I_{\widetilde{A}_{m-1,n}%
}+\underline{\mu}I_{\widetilde{A}_{m-1,n}^{c}},\label{tmumn}\\
\widetilde{A}_{m-1,n}  &  =\left\{  \widetilde{M}_{m-1,n}\geq-\tfrac
{\overline{\mu}+\underline{\mu}}{2}\left(  {1-\tfrac{m-1}{n}}\right)
+c\right\}  \text{.}\nonumber
\end{align}
}

\begin{description}
\item {\normalsize {(1)} Assume that $\varphi$ is decreasing on $(c,\infty)$.
Then
\begin{equation}
\lim\limits_{n\rightarrow\infty}\sup_{Q\in\mathcal{P}}E_{Q}\left[
\varphi\left(  M_{n,n}\right)  \right]  =\mathbb{E}_{\left[  \underline{\mu
},\overline{\mu}\right]  }[\varphi\left(  B_{1}\right)  ]. \label{sym-clt1}%
\end{equation}
}

\item {\normalsize {(2)} Assume that $\varphi$ is increasing on $(c,\infty).$
Then
\begin{equation}
\liminf\limits_{n\rightarrow\infty}\sup_{Q\in\mathcal{P}}E_{Q}\left[
\varphi\left(  \widetilde{M}_{n,n}\right)  \right]  \geq\mathbb{E}_{\left[
\underline{\mu},\overline{\mu}\right]  }[\varphi\left(  B_{1}\right)  ].
\label{sym-clt2-0}%
\end{equation}
}\noindent{\normalsize
}\noindent{\normalsize Furthermore, assume also that
\begin{equation}
\lim_{\delta\rightarrow0}\limsup_{n\rightarrow\infty}\frac{1}{n}\sum_{m=1}%
^{n}\sup_{Q\in\mathcal{P}}E_{Q}\left[  \left\vert E_{Q}[X_{m}|\mathcal{G}%
_{m-1}]-\widetilde{\mu}_{m}^{n}\right\vert I_{\widetilde{A}_{m-1,n}^{\delta}%
}\right]  =0, \label{condition-1}%
\end{equation}
}\noindent{\normalsize where
\[
\widetilde{A}_{m-1,n}^{\delta}=\left\{  \left\vert \widetilde{M}%
_{m-1,n}+\tfrac{\overline{\mu}+\underline{\mu}}{2}\left(  {1-\tfrac{m-1}{n}%
}\right)  -c\right\vert \leq\delta\right\}  ,\quad\delta>0.
\]
}\noindent{\normalsize Then
\begin{equation}
\lim\limits_{n\rightarrow\infty}\sup_{Q\in\mathcal{P}}E_{Q}\left[
\varphi\left(  \widetilde{M}_{n,n}\right)  \right]  =\mathbb{E}_{\left[
\underline{\mu},\overline{\mu}\right]  }[\varphi\left(  B_{1}\right)  ].
\label{sym-clt2}%
\end{equation}
}

{\normalsize
{ }}
\end{description}
\end{thm}

{\normalsize Consider (1). Given $n\geq1$, $\{\mu_{m}^{n}:m\leq n\}$ are
defined recursively with $\mu_{m}^{n}$ being a function of $(X_{1}%
,...,X_{m-1})$. The definition is clearer in the special case where%
\begin{equation}
c=0\text{ and }\underline{\mu}+\overline{\mu}=0\text{.} \label{symm-complete}%
\end{equation}
Then
\begin{equation}
\mu_{m}^{n}=\left\{
\begin{array}
[c]{ccc}%
\overline{\mu} & \text{if} & \frac{1}{n}\sum\limits_{i=1}^{m-1}X_{i}+\frac
{1}{\sqrt{n}}\sum\limits_{i=1}^{m-1}\frac{1}{\sigma}\left(  X_{i}-\mu_{i}%
^{n}\right)  \leq0,\\
\underline{\mu} & \text{if} & \frac{1}{n}\sum\limits_{i=1}^{m-1}X_{i}+\frac
{1}{\sqrt{n}}\sum\limits_{i=1}^{m-1}\frac{1}{\sigma}\left(  X_{i}-\mu_{i}%
^{n}\right)  >0.
\end{array}
\right.  \label{mu-symm-complete}%
\end{equation}
That is, $\mu_{m}^{n}$ is set as large (small) as possible when }%
$M_{m-1,n}\leq(>)0,$ {\normalsize hence lying in the region where $\varphi$ is
increasing (decreasing).
}

Conclude that\ the theorem delivers the statistic $M_{n}=M_{n,n}$ defined in
(\ref{Mmn}),
and, through the upper expectation of $\varphi\left(  M_{n}\right)  $ for the
indicated set of functions $\varphi$, (\ref{sym-clt1}) gives information about
its asymptotic distribution. Moreover, in combination with (\ref{upperBSDE}),
this information can be expressed in closed-form when $\varphi$ is the
indicator for an interval (for a simpler proof than in \cite{chenBSDE} see
Lemma \ref{explicit-solution} below). In particular, we have: For any $a<b\in
${\normalsize $\mathbb{R}$,{
\begin{align}
&  \lim\limits_{n\rightarrow\infty}\sup_{Q\in\mathcal{P}}Q\left(  a\leq
M_{n}\leq b\right) \nonumber\\
=  &  \left\{
\begin{array}
[c]{lcc}%
\Phi_{-\overline{\mu}}\left(  -a\right)  -e^{-\frac{(\overline{\mu}%
-\underline{\mu})(b-a)}{2}}\Phi_{-\overline{\mu}}\left(  -b\right)  &  &
\text{if }a+b\geq\overline{\mu}+\underline{\mu},\\
\Phi_{\underline{\mu}}\left(  b\right)  -e^{-\frac{(\overline{\mu}%
-\underline{\mu})(b-a)}{2}}\Phi_{\underline{\mu}}\left(  a\right)  &  &
\text{if }a+b<\overline{\mu}+\underline{\mu}.
\end{array}
\right.  \label{sym-indi-clt}%
\end{align}
}
}

Similarly, part (2) produces the statistic $\widetilde{M}_{n}=\widetilde
{M}_{n,n}$ defined in (\ref{tMmn}),
that plays a corresponding role. A difference is that only the inequality
(\ref{sym-clt2-0}) is proven in general, though equality obtains under the
condition (\ref{condition-1}). In that case one obtains (as above) that: For
any $a<b\in\mathbb{R}$,{\normalsize {
\begin{align*}
&  \lim\limits_{n\rightarrow\infty}\inf_{Q\in\mathcal{P}}Q\left(
a\leq\widetilde{M}_{n}\leq b\right) \\
=  &  \left\{
\begin{array}
[c]{lcc}%
\Phi_{-\underline{\mu}}\left(  -a\right)  -e^{\frac{(\overline{\mu}%
-\underline{\mu})(b-a)}{2}}\;\Phi_{-\underline{\mu}}\left(  -b\right)  &  &
\text{if }a+b\geq\overline{\mu}+\underline{\mu},\\
\Phi_{\overline{\mu}}\left(  b\right)  -e^{\frac{(\overline{\mu}%
-\underline{\mu})(b-a)}{2}}\;\Phi_{\overline{\mu}}\left(  a\right)  &  &
\text{if }a+b<\overline{\mu}+\underline{\mu}.
\end{array}
\right.
\end{align*}
}}

Finally, we note that {\normalsize (\ref{condition-1}) is easily verified when
$\overline{\mu}=\underline{\mu}=\mu$, because then $E_{Q}[X_{m}|\mathcal{G}%
_{m-1}]=\mu=$}${\widetilde{\mu}_{m}^{n}}$, {\normalsize for any $Q\in
\mathcal{P}$} {\normalsize and $1\leq m\leq n$. }\noindent More generally,
({\normalsize \ref{condition-1}) is satisfied if
\[
\lim_{\delta\rightarrow0}\limsup_{n\rightarrow\infty}\frac{1}{n}\sum_{m=1}%
^{n}\sup_{Q\in\mathcal{P}}Q(\widetilde{A}_{m-1,n}^{\delta})=0.
\]
}

{\normalsize When $c=\pm\infty$, the assumptions in the theorem imply global
monotonicity conditions for }$\varphi${\normalsize , and lead to the fixed
means $\underline{\mu}$ and $\overline{\mu}$ replacing the stochastic means
appearing in { (\ref{CLT}), }(\ref{sym-clt1}) and (\ref{sym-clt2-0})
respectively, and to the normal as the limiting distribution. These features
apply, in particular, to one-sided indicators $I_{(-\infty,b]}$ and
$I_{[a,\infty)}$, and stand in contrast to the implications described above
for two-sided indicators $I_{[a,b]}$. }

\begin{corol}
{\normalsize \label{cor-monotonic} Adopt the assumptions in Theorem
\ref{thm-CLT} and assume that $\varphi\in C([-\infty,\infty])$. }

\begin{description}
\item {\normalsize {(1)} If $\varphi$ is decreasing on $\mathbb{R}$, then
\begin{equation}
\lim\limits_{n\rightarrow\infty}\sup_{Q\in\mathcal{P}}E_{Q}\left[
\varphi\left(  \frac{1}{n}\sum\limits_{i=1}^{n}X_{i}+\frac{1}{\sqrt{n}}%
\sum\limits_{i=1}^{n}\frac{1}{\sigma}(X_{i}-\underline{\mu})\right)  \right]
=\int\varphi\left(  t\right)  d\Phi_{\underline{\mu}}\left(  t\right)
\text{.} \label{decreasing-clt}%
\end{equation}
}

\item {\normalsize {(2)} If $\varphi$ is increasing on $\mathbb{R}$, then
\begin{equation}
\lim\limits_{n\rightarrow\infty}\sup_{Q\in\mathcal{P}}E_{Q}\left[
\varphi\left(  \frac{1}{n}\sum\limits_{i=1}^{n}X_{i}+\frac{1}{\sqrt{n}}%
\sum\limits_{i=1}^{n}\frac{1}{\sigma}(X_{i}-\overline{\mu})\right)  \right]
=\int\varphi\left(  t\right)  d\Phi_{\overline{\mu}}\left(  t\right)  .
\label{increasing-clt}%
\end{equation}
{\Large \noindent} }
\end{description}
\end{corol}

\subsection{{\protect\normalsize \noindent An application to hypothesis
testing\label{section-hyp}}}

{\normalsize We give an illustrative application of Theorem \ref{thm-special}
to hypothesis testing that demonstrates tractability; a more comprehensive
study of statistical applications is beyond the scope of this paper. Here we
exploit also explicit solutions to BSDEs established in \cite{chenBSDE}, an
example of which is provided in (\ref{upperBSDE}). }

{\normalsize Consider the model%
\[
X_{i}=\theta+Y_{i}\text{, }i=1,2,...\text{,}%
\]
where $\theta\in\mathbb{R}$ is a parameter of interest, $(X_{i})$ describes
observable data, and $\left(  Y_{i}\right)  $ is an unobservable error
process. The usual assumption on errors is that they are i.i.d. with zero
mean. Since errors are unobservable, a weaker a priori specification is
natural. Thus, for example, assume the IID model $\mathcal{P}^{IID}$, and for
simplicity, that errors have means that lie in the interval $\left[
-\kappa,\kappa\right]  $. Both the variance $\sigma$ and $\kappa$, which
$\,$measures ambiguity, are assumed known. In the special case $\kappa=0$,
$\theta$ is the unknown mean of each $X_{i}$ and one can test hypotheses about
its value by exploiting the classical CLT. Here we generalize that test
procedure to cover $\kappa>0$. }

{\normalsize Let $\varphi=I_{\left[  a,b\right]  }$, which is symmetric with
center $c=\frac{a+b}{2}$, and define the statistic{ $M_{n}=M_{n,n}$ by
(\ref{Mmn})}. It follows from Theorem \ref{thm-special}(1) and
(\ref{upperBSDE}) that, for any $\theta$, (see Appendix \ref{app-hyp}), {
\begin{align}
&  \lim\limits_{n\rightarrow\infty}\sup_{Q\in\mathcal{P}^{IID}}Q\left(
\{M_{n}-b\leq\theta\leq M_{n}-a\}\right) \nonumber\\
&  =\lim\limits_{n\rightarrow\infty}\sup_{Q\in\mathcal{P}^{IID}}Q\left(
\{a\leq M_{n}-\theta\leq b\}\right)  =\mathbb{E}_{\left[  -\kappa
,\kappa\right]  }[I_{[a,b]}\left(  B_{1}\right)  ]\label{Mtheta}\\
&  =\left\{
\begin{array}
[c]{lcc}%
\Phi_{-\kappa}\left(  -a\right)  -e^{-\kappa\left(  b-a\right)  }\Phi
_{-\kappa}\left(  -b\right)  &  & \text{if }a+b\geq0,\\
\Phi_{-\kappa}\left(  b\right)  -e^{-\kappa\left(  b-a\right)  }\Phi_{-\kappa
}\left(  a\right)  &  & \text{if }a+b<0.
\end{array}
\right. \nonumber
\end{align}
} }

{\normalsize The null hypothesis is $H_{0}:\theta\in\Theta$ and the
alternative is $H_{1}$: $\theta\not \in \Theta$, for some $\Theta
\subset\mathbb{R}$. A nonstandard feature is that there are several
probability laws that conceivably describe the data even given a specific
$\theta$. One test procedure is to accept $H_{0}$ if and only if the realized
statistic $M_{n}$ is "sufficiently consistent" with some $\theta\in\Theta$ and
some probability law in $\mathcal{P}^{IID}$. Precisely, choose $\left[
a,b\right]  $ so that $\mathbb{E}_{\left[  -\kappa,\kappa\right]  }[
I_{[a,b]}\left(  B_{1}\right)  ]=1-\alpha$, for a suitable $\alpha$, and
accept $H_{0}$ if and only if $\mathcal{C}_{n}\cap\Theta\not =\varnothing$,
where the random interval $\mathcal{C}_{n}$ is given by%
\[
\mathcal{C}_{n}=\left[  M_{n}-b,M_{n}-a\right]  \text{.}%
\]
Then, if $H_{0}$ is true, in the limit for large samples the (upper)
probability of acceptance is approximately $1-\alpha$. The upper probability
of wrongly rejecting $H_{0}$ is typically greater than $\alpha$ because of the
multiplicity of measures in $\mathcal{P}^{IID}$:
\begin{align*}
\sup_{Q\in\mathcal{P}^{IID}}Q\left(  \{\mathcal{C}_{n}\cap\Theta
=\varnothing\}\right)   &  =1-\inf_{Q\in\mathcal{P}^{IID}}Q\left(
\{\mathcal{C}_{n}\cap\Theta\not =\varnothing\}\right) \\
&  \geq1-\sup_{Q\in\mathcal{P}^{IID}}Q\left(  \{\mathcal{C}_{n}\cap
\Theta\not =\varnothing\}\right)  \text{.}%
\end{align*}
}

{\normalsize Let $\Theta=\{\theta_{0}\}$ and suppose that the truth is
$\theta=\theta_{1}\equiv\theta_{0}+\xi$, $\xi\not =0$. Then the limiting upper
probability of wrongly accepting $\theta_{0}$ is given by
\begin{align*}
&  \lim\limits_{n\rightarrow\infty}\sup_{Q\in\mathcal{P}^{IID}}Q\left(
\{a\leq M_{n}-\theta_{0}\leq b\}\right) \\
&  =\lim\limits_{n\rightarrow\infty}\sup_{Q\in\mathcal{P}^{IID}}Q\left(
\{a\leq M_{n}-\theta_{1}+\xi\leq b\}\right) \\
&  =\mathbb{E}_{\left[  -\kappa+\xi,\kappa+\xi\right]  }[I_{[a,b]}\left(
B_{1}\right)  ]=\mathbb{E}_{\left[  -\kappa,\kappa\right]  }[I_{[a-\xi,b-\xi
]}\left(  B_{1}\right)  ]
\end{align*}
We emphasize that, given $a$ and $b$, $\mathbb{E}_{\left[  -\kappa
,\kappa\right]  }[I_{[a-\xi,b-\xi]}\left(  B_{1}\right)  ]$ can be expressed
in closed-form (using (\ref{upperBSDE})); and $a$ and $b$ might be chosen by
solving
\begin{equation}
\min_{a\leq b}\mathbb{E}_{\left[  -\kappa,\kappa\right]  }[I_{[a-\xi,b-\xi
]}\left(  B_{1}\right)  ]\text{ ~s.t. }\mathbb{E}_{\left[  -\kappa
,\kappa\right]  }[I_{[a,b]}\left(  B_{1}\right)  ]\geq1-\alpha\text{.}
\label{min}%
\end{equation}
}

\section{Further discussion\label{section-discuss}}

{\normalsize We turn attention to two nonstandard features of the CLT Theorem
\ref{thm-CLT} mentioned only briefly above. One novel feature is that the
limit is defined by the BSDE (\ref{BSDE}). It is shown in \cite[Theorem
2.2]{CE}, using the Girsanov Theorem, that $\mathbb{E}_{\left[  \underline
{\mu},\overline{\mu}\right]  }\left[  \cdot\right]  $ is also an upper
expectation for a set of probability measures, where these are defined on
$C\left(  \left[  0,1\right]  \right)  $, the space of continuous
trajectories. Moreover, measures in this set define differing models of the
underlying stochastically varying (instantaneous) drift. Stochastic
variability of the drift is suggested by (\ref{BSDE}), according to which it
varies between $\underline{\mu}$ and $\overline{\mu}$ depending on the sign of
$Z_{s}$. When the mean is unambiguous ($\underline{\mu}=\overline{\mu}=\mu$),
then the drift is constant and $\mathbb{E}_{\left[  \underline{\mu}%
,\overline{\mu}\right]  }\left[  \cdot\right]  $ reduces to a linear
expectation with normal distribution. However, in general, \emph{the limit is
given by a two parameter (}$\underline{\mu}$\emph{ and }$\overline{\mu}%
$\emph{) family of upper expectations} that model stochastically varying drift
in a continuous-time context. This limiting family is common to a large class
of models (for example, to all IID models in section \ref{section-IID}), thus
endowing the BSDE with special significance for asymptotic approximations in a
sequential context with considerable unstructured heterogeneity in means. }

{\normalsize The other notable feature is that the r.v. on the left in
(\ref{CLT}) combines the sample average, typical of LLNs, with a term that is
more typical of CLTs. In the classical i.i.d. or martingale model, including
the empirical average $\frac{1}{n}{\sum_{i=1}^{n}X_{i}}$ is of little
consequence for the CLT because the LLN permits replacing it by the common
mean of the $X_{i}$s, thereby merely shifting the mean of the limiting normal
distribution. This supports the common view that, in large samples, sample
average reveals location of the population distribution while the ($\sqrt{n}%
$-scaled) average deviation from the mean reflects the distribution about that
location. But this separation of roles is not true in our framework because
empirical averages need not converge given ambiguity (see related LLNs in
\cite{esIID, mm, peng2008,chen}, for example). Next we show that both a LLN
and a "more standard-looking" CLT can be obtained from Theorem \ref{thm-CLT} -
the former as a corollary and the latter by adapting the proof of our CLT.
However, our CLT is more than the "sum of these parts"; for example, a
BSDE-based limit as in (\ref{CLT}) is not present or at all evident from
inspection of the two derivative results. }

\begin{thm}
{\normalsize \label{thm-CLT2}Let the sequence $\left(  X_{i}\right)  $ be such
that $X_{i}\in\mathcal{H}$ for each $i$, and where $\left(  X_{i}\right)  $
satisfies (\ref{mubar}) and (\ref{condvar}), with conditional upper and lower
means $\overline{\mu}$ and $\underline{\mu}$, and {unambiguous conditional
variance} $\sigma^{2}>0$. Suppose also that $(X_{i})$ satisfies the Lindeberg
condition (\ref{linder}). Then, for any $\varphi\in C\left(  \left[
-\infty,\infty\right]  \right)  $,%
\begin{equation}
\lim\limits_{n\rightarrow\infty}\sup_{Q\in\mathcal{P}}E_{Q}\left[
\varphi\left(  \frac{1}{\sqrt{n}}\sum\limits_{i=1}^{n}\frac{1}{\sigma}%
{(X_{i}-E_{Q}[X_{i}|\mathcal{G}_{i-1}]}\right)  \right]  =\int\varphi\left(
t\right)  d\Phi_{0}\left(  t\right)  \text{.} \label{CLT2}%
\end{equation}
}
\end{thm}

{\normalsize \noindent A proof can be constructed along the lines of that of
Theorem \ref{thm-CLT} as indicated in Remarks \ref{remark-CLT2} and
\ref{remark-CLT2b} and in Appendix \ref{app-CLT2}. }

{\normalsize In comparison with Theorem \ref{thm-CLT}, the above theorem drops
rectangularity and yields a limit given by the normal distribution as in the
classical martingale CLT. This is intuitive since, as argued earlier, the
non-normal limit in Theorem \ref{thm-CLT} reflects agnosticism about the
\emph{stochastic} variation in means, which is implicit in rectangularity. The
difference between the two theorems can be seen clearly through their
canonical examples, the IID model $\mathcal{P}^{IID}$ for Theorem
\ref{thm-CLT} and, we would argue, the product model $\mathcal{P}^{prod}$ for
the second theorem. The noted agnosticism motivates $\mathcal{P}^{IID}$ but is
excluded by $\mathcal{P}^{prod}~$(section \ref{section-IID}).
}

{\normalsize Another point of comparison is that while Theorem \ref{thm-CLT2}
adopts weaker assumptions, there is a sense in which it also produces a weaker
result. For example, it does not discriminate between the IID and product
models - the limit is the same for both. In contrast, it can be shown that
Theorem \ref{thm-CLT}, where the sample average term is included, is not valid
for the product model. }

{\normalsize Theorem \ref{thm-CLT2} also clarifies the relation (outlined in
the introduction) between this paper and CLTs by Peng and coauthors. In
particular, in common with (\ref{CLT2}) and unlike (\ref{CLT}), \cite[Theorem
3.2]{peng2019} excludes the sample average term and delivers a normal
distribution in the limit. }

{\normalsize Finally, we show that if Theorem \ref{thm-CLT} is modified so as
to include only the sample average term, then one obtains the following LLN.
(The idea in the proof, found in Appendix \ref{app-LLN}, is first to note the
appropriate form of (\ref{CLT}) when the deviation term is weighted by
$\alpha>0$, and then to let $\alpha\to0$.) }

\begin{corol}
{\normalsize \label{cor-LLN} Adopt the assumptions in Theorem \ref{thm-CLT}.
Then, for any $\varphi\in C\left(  \left[  -\infty,\infty\right]  \right)  $,
\begin{equation}
\lim\limits_{n\rightarrow\infty}\sup_{Q\in\mathcal{P}}E_{Q}\left[
\varphi\left(  \frac{1}{n}{\sum_{i=1}^{n}X_{i}}\right)  \right]
=\sup\limits_{\underline{\mu}\leq\mu\leq\overline{\mu}}\varphi\left(
\mu\right)  \text{.} \label{LLN}%
\end{equation}
}
\end{corol}

{\normalsize \noindent For example, if $\varphi=I_{\left[  a,b\right]  }$,
then (\ref{LLN}) takes the form{
\[
\lim\limits_{n\rightarrow\infty}\sup_{Q\in\mathcal{P}}Q\left(  a\leq\frac
{1}{n}{\sum_{i=1}^{n}X_{i}}\leq b\right)  =\left\{
\begin{array}
[c]{cc}%
1 & \text{ if }\left[  a,b\right]  \cap\left[  \underline{\mu},\overline{\mu
}\right]  \not =\varnothing\\
0 & \text{otherwise.}%
\end{array}
\right.
\]
} }

\section{Main proofs\label{section-proofs}}

{\normalsize This section proves Theorems \ref{thm-CLT} and \ref{thm-special}.
Throughout we use the following well-known properties of (conditional) upper
expectations, understood to hold for all $X$ and $Y$ in $\mathcal{H}$, and all
$n\geq0$. }

\begin{enumerate}
\item {\normalsize Monotonicity: $X\geq Y$ implies $\mathbb{E}\left[
X\mid\mathcal{G}_{n}\right]  \geq\mathbb{E}\left[  Y\mid\mathcal{G}%
_{n}\right]  . $ }

\item {\normalsize Sub-additivity: $\mathbb{E}\left[  X+Y\mid\mathcal{G}%
_{n}\right]  \leq\mathbb{E}\left[  X\mid\mathcal{G}_{n}\right]  +\mathbb{E}%
\left[  Y\mid\mathcal{G}_{n}\right]  . $ }

\item {\normalsize Homogeneity: If $Z$ is $\mathcal{G}_{n}$ measurable,
\[
\mathbb{E}\left[  ZX\mid\mathcal{G}_{n}\right]  =Z^{+}\mathbb{E}\left[
X\mid\mathcal{G}_{n}\right]  -Z^{-}\mathcal{E}\left[  X\mid\mathcal{G}%
_{n}\right]  .
\]
}

\item {\normalsize Translation homogeneity: If $Z$ is $\mathcal{G}_{n}$
measurable,
\[
\mathbb{E}[Z+X\mid\mathcal{G}_{n}]=Z+\mathbb{E}[X\mid\mathcal{G}_{n}].
\]
}
\end{enumerate}

{\normalsize The assumptions in Theorem \ref{thm-CLT} are adopted throughout.
As indicated following (\ref{BSDE-1-1}), $(B_{t})$ is a standard Brownian
motion on a filtered probability space {$(\Omega^{\ast},\mathcal{F}^{\ast
},\{\mathcal{F}_{t}\},P^{\ast})$}; $\{\mathcal{F}_{t}\}$ is the natural
filtration generated by $(B_{t})$. }

{\normalsize For both theorems, we prove them first for the special case
where
\begin{equation}
-\underline{\mu}=\overline{\mu}=\kappa\geq0, \label{kappa}%
\end{equation}
that is,
\[
\mathbb{E}[X_{i}\mid\mathcal{G}_{i-1}]=\kappa,\quad\mathcal{E}[X_{i}%
\mid\mathcal{G}_{i-1}]=-\kappa.
\]
Then the results asserted for general $\underline{\mu}$ and $\overline{\mu}$
are established by applying the preceding special case to $\left(
Y_{i}\right)  $, where $Y_{i}=X_{i}-\tfrac{\overline{\mu}+\underline{\mu}}{2}%
$, and thus
\[
\mathbb{E}[Y_{i}\mid\mathcal{G}_{i-1}]=\frac{\overline{\mu}-\underline{\mu}%
}{2},\quad\mathcal{E}[Y_{i}\mid\mathcal{G}_{i-1}]=-\frac{\overline{\mu
}-\underline{\mu}}{2}\text{.}%
\]
}

\subsection{Lemmas}

{\normalsize The following lemmas prepare the groundwork for proofs of both
Theorems \ref{thm-CLT} and \ref{thm-special}. The special case (\ref{kappa})
is assumed throughout unless specified otherwise. }

{\normalsize For any fixed $\epsilon>0$, define $g_{\epsilon}:\mathbb{R}%
\rightarrow\mathbb{R}$ by
\begin{equation}
g_{\epsilon}(z)=\kappa\left(  \sqrt{z^{2}+\epsilon^{2}}-\epsilon\right)
\text{.} \label{gep}%
\end{equation}
Obviously, $g_{\epsilon}(0)=0$ and $g_{\epsilon}$ is symmetric with center
$c=0$. For any suitably integrable random variable $\xi\in\mathcal{F}_{1}$,
define $g$-expectation by $\mathbb{E}_{g_{\epsilon}}[\xi]=Y_{0}^{\epsilon}$,
where $\left(  Y_{t}^{\epsilon},Z_{t}^{\epsilon}\right)  $ is the unique
solution to the BSDE
\begin{equation}
Y_{t}^{\epsilon}=\xi+\int_{t}^{1}g_{\epsilon}\left(  Z_{s}^{\epsilon}\right)
ds-\int_{t}^{1}Z_{s}^{\epsilon}dB_{s},\;0\leq t\leq1\text{.} \label{BSDEep}%
\end{equation}
(Existence of a unique solution follows from \cite{pp}.) Moreover, by
\cite[Proposition 2.1]{KPQ}, { for any suitably integrable $\xi\in
\mathcal{F}_{1},$}
\[
\mathbb{E}_{g_{\epsilon}}[\xi]\rightarrow\mathbb{E}_{g_{0}}[\xi]\text{, as
}\epsilon\rightarrow0,
\]
where $\mathbb{E}_{g_{0}}[\xi]=Y_{0}^{0}$, and $\left(  Y_{t}^{0},Z_{t}%
^{0}\right)  $ is the unique solution to the BSDE (\ref{BSDEep}) for the
extreme case corresponding to $\epsilon=0$, where
\begin{equation}
g_{0}\left(  z\right)  =\kappa|z|\text{,} \label{g0}%
\end{equation}
and $\mathbb{E}_{g_{0}}$ is alternative notation for $\mathbb{E}_{\left[
-\kappa,\kappa\right]  }$. We consider $g_{\epsilon}$ for $\epsilon>0$ in
order to overcome the nondifferentiability of $g_{0}$ at $z=0$. (The relevant
smoothness is exploited in Lemma \ref{lemma-ddp}.)}

{\normalsize We introduce a sequence of functions generated by $g$-expectation
$\mathbb{E}_{g_{\epsilon}}$. Some properties of $g$-expectations can be found
in \cite{peng1997}, we need to prove the following properties. }

{\normalsize Given $\varphi\in C_{b}^{3}(\mathbb{R})$, let $\xi=\varphi\left(
x+B_{1}-B_{\frac{m}{n}}\right)  $ in (\ref{BSDEep}) and define the functions
$\{H_{m,n}\}_{m=0}^{n}$ by
\begin{equation}
H_{m,n}\left(  x\right)  \equiv\mathbb{E}_{g_{\epsilon}}\left[  \varphi\left(
x+B_{1}-B_{\frac{m}{n}}\right)  \right]  ,\;m=0,\cdots,n. \label{BSDEep-1}%
\end{equation}
($\epsilon>0$ is fixed and dependence on $\epsilon$ is suppressed
notationally.) Obviously,
\[
H_{n,n}(x)=\varphi(x),\;H_{0,n}(x)=\mathbb{E}_{g_{\epsilon}}\left[
\varphi\left(  x+B_{1}\right)  \right]  .
\]
The following lemma shows that the functions $\{H_{m,n}\}_{m=0}^{n}$ are
suitably differentiable given that $\varphi\in C_{b}^{3}(\mathbb{R})$ and
$\epsilon>0$. }

\begin{lemma}
{\normalsize \label{lemma-ddp}The functions $\{H_{m,n}\}_{m=0}^{n}$ satisfy: }

\begin{description}
\item {\normalsize {(1)} $H_{m,n}\in C_{b}^{2}(\mathbb{R}),$ for
$n\geq1,\;m=0,1,\cdots n.$ }

\item {\normalsize {(2)} The second derivatives of $H_{m,n}$ are uniformly
bounded and Lipschitz continuous with uniform Lipschitz constant for
$\ \{\left(  m,n\right)  :0\leq m\leq n\}$.
}

\item {\normalsize {(3)} Dynamic programming principle: for $n\geq1$,
$m=1,...,n,$ }%
\[
{\normalsize H_{m-1,n}(x)=\mathbb{E}_{g_{\epsilon}}\left[  H_{m,n}\left(
x+B_{\frac{m}{n}}-B_{\frac{m-1}{n}}\right)  \right]  ,\;x\in\mathbb{R}.}%
\]

\item {\normalsize {(4)} Identically distributed: for $n\geq1$, $m=1,...,n,$
\[
\mathbb{E}_{g_{\epsilon}}\left[  H_{m,n}\left(  x+B_{\frac{m}{n}}%
-B_{\frac{m-1}{n}}\right)  \right]  =\mathbb{E}_{g_{\epsilon}}\left[
H_{m,n}\left(  x+B_{\frac{1}{n}}\right)  \right]  ,\;x\in\mathbb{R}.
\]
}
\end{description}
\end{lemma}

{\normalsize \noindent\textbf{Proof: }{(1) and (2):} From the nonlinear
Feynman-Kac Formula \cite[Proposition 4.3]{KPQ}, we have $H_{m,n}%
(x)=u(\frac{m}{n},x)$ and $u$ is the solution of the PDE%
\begin{equation}
\left\{
\begin{array}
[c]{ll}%
\partial_{t} u+\frac{1}{2}\partial^{2}_{xx}u+\kappa\left(  \sqrt{|\partial
_{x}u|^{2}+\epsilon^{2}}-\epsilon\right)  =0, & \\
u(1,x)=\varphi(x). &
\end{array}
\right.  \label{pde1}%
\end{equation}
}

{\normalsize Next we prove that for any $t\in\lbrack0,1]$, $u(t,\cdot)\in
C_{b}^{2}(\mathbb{R})$; $u(t,\cdot),$ $\partial_{x}u(t,\cdot),$ $\partial
_{xx}^{2}u(t,\cdot)$ are bounded uniformly in $t\in\lbrack0,1]$; and for any
$x,x^{\prime}\in\mathbb{R}$, $\exists C>0$ such that $|\partial_{xx}%
^{2}u(t,x)-\partial_{xx}^{2}u(t,x^{\prime})|\leq C|x-x^{\prime}|\text{,
\ }\forall t\in\lbrack0,1]$. }

{\normalsize By the definition of $g_{\epsilon}$,
\[
g_{\epsilon}^{\prime}(z)=\kappa\frac{z}{\sqrt{\epsilon^{2}+z^{2}}}%
\Rightarrow|g_{\epsilon}^{\prime}(z)|\leq\kappa,
\]%
\[
g_{\epsilon}^{\prime\prime}(z)=\kappa\frac{\epsilon^{2}}{(\epsilon^{2}%
+z^{2})^{3/2}}\Rightarrow|g_{\epsilon}^{\prime\prime}(z)|\leq\frac{\kappa
}{\epsilon},
\]%
\[
g_{\epsilon}^{\prime\prime\prime}(z)=-\kappa\frac{3z\epsilon^{2}}%
{(\epsilon^{2}+z^{2})^{5/2}}\Rightarrow|g_{\epsilon}^{\prime\prime\prime
}(z)|\leq\frac{3\kappa}{\epsilon^{2}}.
\]
Consider the following BSDE,
\begin{equation}
Y_{s}^{t,x}=\varphi(x+B_{1}-B_{t})+\int_{s}^{1}g_{\epsilon}(Z_{r}%
^{t,x})dr-\int_{s}^{1}Z_{r}^{t,x}dB_{r},~s\in\lbrack t,1].\label{bsde}%
\end{equation}
Then $u(t,x)=Y_{t}^{t,x}$ is the classical unique solution of PDE
(\ref{pde1}), and
\begin{align}
\partial_{x}Y_{s}^{t,x}= &  \varphi^{\prime}(x+B_{1}-B_{t})+\int_{s}%
^{1}g_{\epsilon}^{\prime}(Z_{r}^{t,x})\partial_{x}Z_{r}^{t,x}dr-\int_{s}%
^{1}\partial_{x}Z_{r}^{t,x}dB_{r},~s\in\lbrack t,1].\label{1stbsde}\\
\partial_{x}^{2}Y_{s}^{t,x}= &  \varphi^{\prime\prime}(x+B_{1}-B_{t})+\int
_{s}^{1}g_{\epsilon}^{\prime\prime}(Z_{r}^{t,x})|\partial_{x}Z_{r}^{t,x}%
|^{2}dr+\int_{s}^{1}g_{\epsilon}^{\prime}(Z_{r}^{t,x})\partial_{x}^{2}%
Z_{r}^{t,x}dr\nonumber\\
&  -\int_{s}^{1}\partial_{x}^{2}Z_{r}^{t,x}dB_{r},~s\in\lbrack
t,1].\label{2ndbsde}%
\end{align}
From standard estimates of BSDEs (\cite{KPQ}), we have, $\forall p\geq2$,
$\forall x\in\mathbb{R}$,
\begin{align*}
&  E_{P^{\ast}}\left[  \sup_{s\in\lbrack t,1]}|Y_{s}^{t,x}|^{p}|\mathcal{F}%
_{t}\right]  +E_{P^{\ast}}\left[  (\int_{t}^{1}|Z_{s}^{t,x}|^{2}ds)^{\frac
{p}{2}}|\mathcal{F}_{t}\right]  \\
\leq &  C_{p}^{0}E_{P^{\ast}}\left[  |\varphi(x+B_{1}-B_{t})|^{p}+(\int
_{t}^{1}|g_{\epsilon}(0)|dr)^{p}|\mathcal{F}_{t}\right]  \leq{C_{p}^{0}%
\Vert\varphi\Vert^{p}};\\
&  E_{P^{\ast}}\left[  \sup_{s\in\lbrack t,1]}|\partial_{x}Y_{s}^{t,x}%
|^{p}|\mathcal{F}_{t}\right]  +E_{P^{\ast}}\left[  (\int_{t}^{1}|\partial
_{x}Z_{s}^{t,x}|^{2}ds)^{\frac{p}{2}}|\mathcal{F}_{t}\right]  \\
\leq &  C_{p}^{1}E_{P^{\ast}}\left[  |\varphi^{\prime}(x+B_{1}-B_{t}%
)|^{p}|\mathcal{F}_{t}\right]  \leq C_{p}^{1}\Vert\varphi^{\prime}\Vert^{p};\\
&  E_{P^{\ast}}\left[  \sup_{s\in\lbrack t,1]}|\partial_{x}^{2}Y_{s}%
^{t,x}|^{p}|\mathcal{F}_{t}\right]  +E_{P^{\ast}}\left[  (\int_{t}%
^{1}|\partial_{x}^{2}Z_{s}^{t,x}|^{2}ds)^{\frac{p}{2}}|\mathcal{F}_{t}\right]
\\
\leq &  C_{p}^{2}E_{P^{\ast}}\left[  |\varphi^{\prime\prime}(x+B_{1}%
-B_{t})|^{p}+(\int_{t}^{1}|g_{\epsilon}^{\prime\prime}({Z_{s}^{t,x}%
})||\partial_{x}Z_{s}^{t,x}|^{2}ds)^{p}|\mathcal{F}_{t}\right]  \\
\leq &  C_{p}^{2}(\Vert\varphi^{\prime\prime}\Vert^{p}+(\tfrac{\kappa
}{\epsilon})^{p}C_{2p}^{1}\Vert\varphi^{\prime}\Vert^{2p}),
\end{align*}
where $C_{p}^{0},C_{p}^{1},C_{2p}^{1},C_{p}^{2}$ are constants independent of
$t$ and $\Vert f\Vert=\sup_{x\in\mathbb{R}}f(x)$ denote the sup norm of
function $f$. Then, for any $t\in\lbrack0,1]$, $u(t,\cdot)\in C_{b}%
^{2}(\mathbb{R})$ and $u(t,\cdot),\partial_{x}u(t,\cdot),\partial_{xx}%
^{2}u(t,\cdot)$ are bounded uniformly in $t\in\lbrack0,1]$. }

{\normalsize From (\ref{bsde}), the Malliavin derivative satisfies, for
$u\in\lbrack t,s)$,
\[
D_{u}Y_{s}^{t,x}=\varphi^{\prime}(x+B_{1}-B_{t})+\int_{s}^{1}g_{\epsilon
}^{\prime}(Z_{r}^{t,x})D_{u}Z_{r}^{t,x}dr-\int_{s}^{1}D_{u}Z_{r}^{t,x}%
dB_{r},~s\in\lbrack t,1].
\]
From standard estimates for BSDEs, for $s\in\lbrack t,1]$, we have
\[
E_{P^{\ast}}\left[  \left(  \int_{s}^{1}|D_{u}Z_{r}^{t,x}|^{2}dr\right)
^{\frac{p}{2}}|\mathcal{F}_{s}\right]  \leq C_{p}^{1}E_{P^{\ast}}\left[
|\varphi^{\prime}(s+B_{1}-B_{t})|^{p}|\mathcal{F}_{s}\right]  \leq C_{p}%
^{1}\Vert\varphi^{\prime}\Vert^{p}\text{,}%
\]
and from (\ref{1stbsde}), we have
\begin{align*}
D_{u}\left[  \partial_{x}Y_{s}^{t,x}\right]  = &  \varphi^{\prime\prime
}(x+B_{1}-B_{t})+\int_{s}^{1}g_{\epsilon}^{\prime\prime}(Z_{r}^{t,x}%
)D_{u}[Z_{r}^{t,x}]\partial_{x}Z_{r}^{t,x}dr\\
&  +\int_{s}^{1}g_{\epsilon}^{\prime}(Z_{r}^{t,x})D_{u}\left[  \partial
_{x}Z_{r}^{t,x}\right]  dr-\int_{s}^{1}D_{u}\left[  \partial_{x}Z_{r}%
^{t,x}\right]  dB_{r},~s\in\lbrack t,1].
\end{align*}
Let $d\tilde{B}_{s}=dB_{s}-g_{\epsilon}^{\prime}(Z_{s}^{t,x})ds$, $\rho
_{s}=\exp\left\{  \int_{0}^{s}g_{\epsilon}^{\prime}(Z_{r}^{t,x})dB_{r}%
-\frac{1}{2}\int_{0}^{s}|g_{\epsilon}^{\prime}(Z_{r}^{t,x})|^{2}dr\right\}  $,
and $E_{P^{\ast}}[\frac{d\tilde{P}}{dP^{\ast}}|_{\mathcal{F}_{s}}]=\rho_{s}$.
Then,} {\footnotesize
\begin{align*}
&  \left\vert D_{u}\left[  \partial_{x}Y_{s}^{t,x}\right]  \right\vert \\
= &  \left\vert E_{\tilde{P}}\left[  \varphi^{\prime\prime}(x+B_{1}%
-B_{t})+\int_{s}^{1}g_{\epsilon}^{\prime\prime}(Z_{r}^{t,x})D_{u}[Z_{r}%
^{t,x}]\partial_{x}Z_{r}^{t,x}dr|\mathcal{F}_{s}\right]  \right\vert \\
\leq &  \Vert\varphi^{\prime\prime}\Vert+\frac{\kappa}{\epsilon}E_{P^{\ast}%
}\left[  \rho_{1}(\rho_{s})^{-1}\cdot\int_{s}^{1}|D_{u}[Z_{r}^{t,x}%
]|\cdot|\partial_{x}Z_{r}^{t,x}|dr|\mathcal{F}_{s}\right]  \\
\leq &  \Vert\varphi^{\prime\prime}\Vert+\frac{\kappa}{\epsilon}%
M_{s}\;\bigg(E_{P^{\ast}}\Big[\big(\int_{s}^{1}|D_{u}[Z_{r}^{t,x}%
]|^{2}dr\big)^{2}|\mathcal{F}_{s}\Big]\bigg)^{\frac{1}{4}}\bigg(E_{P^{\ast}%
}\Big[\big(\int_{s}^{1}|\partial_{x}Z_{r}^{t,x}|^{2}dr\big)^{2}|\mathcal{F}%
_{s}\Big]\bigg)^{\frac{1}{4}}\\
\leq &  K\text{,}%
\end{align*}
}{\normalsize where $M_{s}\equiv\left(  E_{P^{\ast}}\left[  \rho_{1}^{2}%
(\rho_{s})^{-2}|\mathcal{F}_{s}\right]  \right)  ^{\frac{1}{2}}$, {and
satisfies
\begin{align*}
M_{s}^{2}= &  E_{P^{\ast}}\left[  e^{\int_{s}^{1}2g_{\epsilon}^{\prime}%
(Z_{r}^{t,x})dB_{r}-\int_{s}^{1}|g_{\epsilon}^{\prime}(Z_{r}^{t,x})|^{2}%
dr}|\mathcal{F}_{s}\right]  \\
= &  E_{P^{\ast}}\left[  e^{\int_{s}^{1}2g_{\epsilon}^{\prime}(Z_{r}%
^{t,x})dB_{r}-\frac{1}{2}\int_{s}^{1}|2g_{\epsilon}^{\prime}(Z_{r}^{t,x}%
)|^{2}dr}e^{\int_{s}^{1}|g_{\epsilon}^{\prime}(Z_{r}^{t,x})|^{2}%
dr}|\mathcal{F}_{s}\right]  \leq e^{\kappa^{2}}\text{.}%
\end{align*}
} Here $K$ is a constant that depends on $\kappa,\epsilon,p$,$\Vert
\varphi^{\prime}\Vert$,$\Vert\varphi^{\prime\prime}\Vert.$ With $u\in\lbrack
t,s)$, from (\ref{1stbsde}), the Malliavin derivative satisfies, $\forall
s\in\lbrack t,1]$,
\begin{align*}
D_{u}\left[  \partial_{x}Y_{s}^{t,x}\right]  = &  -\int_{u}^{s}g_{\epsilon
}^{\prime\prime}(Z_{r}^{t,x})D_{u}\left[  Z_{r}^{t,x}\right]  \partial
_{x}Z_{r}^{t,x}dr-\int_{u}^{s}g_{\epsilon}^{\prime}(Z_{r}^{t,x})D_{u}\left[
\partial_{x}Z_{r}^{t,x}\right]  dr\\
&  +\int_{u}^{s}D_{u}\left[  \partial_{x}Z_{r}^{t,x}\right]  dB_{r}%
+\partial_{x}Z_{u}^{t,x}\text{, and}%
\end{align*}%
\[
\lim_{s\downarrow u}D_{u}\left[  \partial_{x}Y_{s}^{t,x}\right]  =\partial
_{x}Z_{u}^{t,x}~~P^{\ast}\text{-a.s.}%
\]
We have,%
\[
|\partial_{x}Z_{u}^{t,x}|\leq K~~du\times dP^{\ast}\text{-}a.s.
\]
Thus, from (\ref{2ndbsde}), by standard estimates for BSDEs again, $\forall
p\geq2$, $\forall x,x^{\prime}\in\mathbb{R}$,}{\small
\begin{align*}
&  E_{P^{\ast}}\left[  \sup_{s\in\lbrack t,1]}|\partial_{x}^{2}Y_{s}%
^{t,x}-\partial_{x}^{2}Y_{s}^{t,x^{\prime}}|^{p}|\mathcal{F}_{t}\right]
+E_{P^{\ast}}\left[  \left(  \int_{t}^{1}|\partial_{x}^{2}Z_{r}^{t,x}%
-\partial_{x}^{2}Z_{r}^{t,x^{\prime}}|^{2}dr\right)  ^{\frac{p}{2}%
}|\mathcal{F}_{t}\right]  \\
\leq &  C_{p}E_{P^{\ast}}\left[  |\varphi^{\prime\prime}(x+B_{1}%
-B_{t})-\varphi^{\prime\prime}(x^{\prime}+B_{1}-B_{t})|^{p}|\mathcal{F}%
_{t}\right]  \\
&  +C_{p}E_{P^{\ast}}\left[  \left(  \int_{t}^{1}|g_{\epsilon}^{\prime\prime
}(Z_{r}^{t,x})(\partial_{x}Z_{r}^{t,x})^{2}-g_{\epsilon}^{\prime\prime}%
(Z_{r}^{t,x^{\prime}})(\partial_{x}Z_{r}^{t,x^{\prime}})^{2}|dr\right)
^{p}|\mathcal{F}_{t}\right]  \\
&  +C_{p}E_{P^{\ast}}\left[  \left(  \int_{t}^{1}|g_{\epsilon}^{\prime}%
(Z_{r}^{t,x})-g_{\epsilon}^{\prime}(Z_{r}^{t,x^{\prime}})||\partial_{x}%
^{2}Z_{r}^{t,x}|dr\right)  ^{p}|\mathcal{F}_{t}\right]  \\
\equiv &  I_{1}+I_{2}+I_{3},
\end{align*}
}{\normalsize where $C_{p}$ is a constant independent of $t$, and $I_{1}%
,I_{2},I_{3}$ satisfied}{\small
\begin{align*}
I_{1}= &  C_{p}E_{P^{\ast}}\left[  |\varphi^{\prime\prime}(x+B_{1}%
-B_{t})-\varphi^{\prime\prime}(x^{\prime}+B_{1}-B_{t})|^{p}|\mathcal{F}%
_{t}\right]  \\
\leq &  C_{p}\Vert\varphi^{\prime\prime\prime}\Vert^{p}|x-x^{\prime}%
|^{p}=C_{1,p}|x-x^{\prime}|^{p},\quad(\text{where }C_{1,p}=C_{p}\Vert
\varphi^{\prime\prime\prime}\Vert^{p})\\
I_{2}= &  C_{p}E_{P^{\ast}}\left[  \left(  \int_{t}^{1}\Big|g_{\epsilon
}^{\prime\prime}(Z_{r}^{t,x})(\partial_{x}Z_{r}^{t,x})(\partial_{x}Z_{r}%
^{t,x}-\partial_{x}Z_{r}^{t,x^{\prime}})\right.  \right.  \\
&  \qquad\qquad+(g_{\epsilon}^{\prime\prime}(Z_{r}^{t,x})-g_{\epsilon}%
^{\prime\prime}(Z_{r}^{t,x^{\prime}}))(\partial_{x}Z_{r}^{t,x})(\partial
_{x}Z_{r}^{t,x^{\prime}})\\
&  \qquad\qquad+\left.  \left.  g_{\epsilon}^{\prime\prime}(Z_{r}%
^{t,x^{\prime}})(\partial_{x}Z_{r}^{t,x^{\prime}})(\partial_{x}Z_{r}%
^{t,x}-\partial_{x}Z_{r}^{t,x^{\prime}})\Big|dr\right)  ^{p}|\mathcal{F}%
_{t}\right]  \\
\leq &  2\cdot3^{p-1}C_{p}(\tfrac{\kappa}{\epsilon})^{p}K^{p}E_{P^{\ast}%
}\left[  \left(  \int_{t}^{1}|\partial_{x}Z_{r}^{t,x}-\partial_{x}%
Z_{r}^{t,x^{\prime}}|dr\right)  ^{p}|\mathcal{F}_{t}\right]  \\
&  +3^{p-1}C_{p}K^{2p}(\tfrac{3\kappa}{\epsilon^{2}})^{p}E_{P^{\ast}}\left[
\left(  \int_{t}^{1}|Z_{r}^{t,x}-Z_{r}^{t,x^{\prime}}|dr\right)
^{p}|\mathcal{F}_{t}\right]  \\
\leq &  C_{2,p}|x-x^{\prime}|^{p},
\end{align*}
}{\normalsize where $C_{2,p}$ is a constant depend on $\kappa,\epsilon
,p,\Vert\varphi^{\prime}\Vert$ and $\Vert\varphi^{\prime\prime}\Vert$%
,}{\small
\begin{align*}
I_{3} &  =C_{p}E_{P^{\ast}}\left[  \left(  \int_{t}^{1}|g_{\epsilon}^{\prime
}(Z_{r}^{t,x})-g_{\epsilon}^{\prime}(Z_{r}^{t,x^{\prime}})||\partial_{x}%
^{2}Z_{r}^{t,x}|dr\right)  ^{p}|\mathcal{F}_{t}\right]  \\
&  \leq C_{p}(\tfrac{\kappa}{\epsilon})^{p}E_{P^{\ast}}\left[  \left(
\int_{t}^{1}|Z_{r}^{t,x}-Z_{r}^{t,x^{\prime}}||\partial_{x}^{2}Z_{r}%
^{t,x}|dr\right)  ^{p}|\mathcal{F}_{t}\right]  \\
&  \leq C_{p}(\tfrac{\kappa}{\epsilon})^{p}\Big(E_{P^{\ast}}\Big[\big(\int
_{t}^{1}|Z_{r}^{t,x}-Z_{r}^{t,x^{\prime}}|^{2}|dr\big)^{p}|\mathcal{F}%
_{t}\Big]\Big)^{\frac{1}{2}}\Big(E_{P^{\ast}}\Big[\big(\int_{t}^{1}%
|\partial_{x}^{2}Z_{r}^{t,x}|^{2}dr\big)^{p}|\mathcal{F}_{t}\Big]\Big)^{\frac
{1}{2}}\\
&  \leq C_{3,p}|x-x^{\prime}|^{p},
\end{align*}
}{\normalsize where $C_{3,p}$ is a constant depend on $\kappa,\epsilon
,p,\Vert\varphi^{\prime}\Vert$ and $\Vert\varphi^{\prime\prime}\Vert$.
Therefore,}{\small
\begin{align*}
&  E_{P^{\ast}}\left[  \sup_{s\in\lbrack t,1]}|\partial_{x}^{2}Y_{s}%
^{t,x}-\partial_{x}^{2}Y_{s}^{t,x^{\prime}}|^{p}|\mathcal{F}_{t}\right]
+E_{P^{\ast}}\left[  \left(  \int_{t}^{1}|\partial_{x}^{2}Z_{r}^{t,x}%
-\partial_{x}^{2}Z_{r}^{t,x^{\prime}}|^{2}dr\right)  ^{\frac{p}{2}%
}|\mathcal{F}_{t}\right]  \\
\leq &  (C_{1,p}+C_{2,p}+C_{3,p})|x-x^{\prime}|^{p}%
\end{align*}
}{\normalsize  Thus we obtain Claims (1) and (2). }

{\normalsize \noindent(3): Follows from Peng's dynamic programming principle
\cite[Theorem 3.2]{peng1992}. }

{\normalsize \noindent(4): It is a direct consequence of \cite[Theorem
3.1]{CE}. \hfill\hfill$\blacksquare$ }

{\normalsize \medskip}

{\normalsize The following lemma is adapted from \cite[Proposition 2.3]%
{BCHMP}. }

\begin{lemma}
{\normalsize \label{lemma6} Suppose that $(b_{t})$ and $(\sigma_{t})$ are two
continuous, bounded $\mathcal{F}_{t}$-adapted processes and that $(X_{t})$ is
of the form
\[
X_{t}=x+\int_{0}^{t}b_{s}ds+\int_{0}^{t}\sigma_{s}dB_{s},\quad x\in
\mathbb{R}.
\]
Then
\[
\lim_{n\rightarrow\infty}n\sup_{x\in\mathbb{R}}\left\vert \mathbb{E}%
_{g_{\epsilon}}\left[  X_{\frac{1}{n}}\right]  -x-\frac{1}{n}g_{\epsilon
}(\sigma_{0})-\frac{1}{n}b_{0}\right\vert =0.
\]
}
\end{lemma}

{\normalsize The next lemma is an immediate consequence. }

\begin{lemma}
{\normalsize \label{lemma-gexpectation} For any $\varphi\in C_{b}%
^{2}(\mathbb{R})$,
\begin{equation}
\lim_{n\rightarrow\infty}n\sup_{x\in\mathbb{R}}\left\vert \mathbb{E}%
_{g_{\epsilon}}\left[  \varphi\left(  x+B_{\frac{1}{n}}\right)  \right]
-\varphi(x)-\frac{1}{n}g_{\epsilon}(\varphi^{\prime}(x))-\frac{1}{2n}%
\varphi^{\prime\prime}(x)\right\vert =0. \label{g-exc}%
\end{equation}
}
\end{lemma}

{\normalsize \noindent\textbf{Proof:} Let $X_{s}\equiv\varphi(x+B_{s})$. By
Ito's formula,
\[
X_{t}=\varphi(x)+\frac{1}{2}\int_{0}^{t}\varphi^{\prime\prime}(x+B_{s}%
)ds+\int_{0}^{t}\varphi^{\prime}(x+B_{s})dB_{s} .
\]
}

{\normalsize \noindent Apply Lemma \ref{lemma6} to complete the proof.
\hfill\hfill$\blacksquare$ }

{\normalsize \smallskip}

\begin{lemma}
{\normalsize \label{lemma-taylor} Let $g_{0}$ be defined by (\ref{g0}). For
any $\varphi\in C_{b}^{3}(\mathbb{R})$, let $\{H_{m,n}\}_{m=0}^{n}$ be the
functions defined in (\ref{BSDEep-1}). Define functions $\{L_{m,n}\}_{m=0}%
^{n}$ by
\begin{equation}
L_{m,n}(x)=H_{m,n}(x)+\frac{1}{n}g_{0}(H_{m,n}^{\prime}(x))+\frac{1}%
{2n}H_{m,n}^{\prime\prime}(x). \label{eq-1}%
\end{equation}
\noindent Let $\{T_{m,n}\}_{m,n\geq0}$ be an array of r.v.s satisfying%
\[
T_{0,n}=0\text{, and }T_{m,n}\in{\mathcal{H}}\text{ is }\mathcal{G}%
_{m}\text{-measurable for all }m\geq1,\ n\geq1\text{,}%
\]
and,{ for any $Q\in\mathcal{P}$,} set $Y_{m}^{Q}=\frac{1}{\sigma}(X_{m}%
-E_{Q}[X_{m}|\mathcal{G}_{m-1}])$. Then{\small
\begin{equation}
\lim_{n\rightarrow\infty}\sum_{m=1}^{n}\left\vert \sup\limits_{Q\in
\mathcal{P}}E_{Q}\left[  H_{m,n}\left(  T_{m-1,n}+\frac{X_{m}}{n}+\frac
{Y_{m}^{Q}}{\sqrt{n}}\right)  \right]  -\sup\limits_{Q\in\mathcal{P}}%
E_{Q}\left[  L_{m,n}(T_{m-1,n})\right]  \right\vert =0. \label{Equa}%
\end{equation}
}}
\end{lemma}

\begin{remark}
{\normalsize \label{remark-CLT2}In Theorem \ref{thm-CLT2}, the sample average
term is absent, and accordingly its proof involves a counterpart of this lemma
where the term $\frac{X_{m}}{n}$ is deleted above. Then the proof of
(\ref{Equa}), so modified, simplifies, in particular, rectangularity is no
longer needed and the generators $g_{\epsilon}$ in (\ref{BSDEep}) and $g_{0}$
in (\ref{eq-1}) can be set equal to $0$. ({Appendix \ref{app-CLT2} provides
some details.) } }
\end{remark}

{\normalsize \noindent\textbf{Proof:} We proceed in two steps. }

{\normalsize \smallskip}

{\normalsize \noindent\textbf{Step 1:} We first give a remainder estimate that
will also be used later in the proof of Lemma \ref{lemma-taylor2}. Let
$\{\theta_{m}\}_{m\geq1}$ be a sequence of $\mathcal{G}_{m-1}$-measurable
random variables satisfying
\[
|\theta_{m}|\leq\kappa,\;\text{for}\ m\geq1.
\]
We prove that%
\begin{align}
\lim_{n\rightarrow\infty}\sum_{m=1}^{n}  &  \left\vert \sup\limits_{Q\in
\mathcal{P}}E_{Q}\left[  H_{m,n}(T_{m-1,n}+\frac{X_{m}}{n}+\frac{X_{m}%
-\theta_{m}}{\sigma\sqrt{n}})\right]  \right. \nonumber\\
&  \qquad\qquad\qquad\qquad\qquad\qquad-\sup\limits_{Q\in\mathcal{P}}%
E_{Q}\left[  F(\theta_{m},m,n)\right]  \bigg\vert=0, \label{remainder-esti}%
\end{align}
where $F(\theta_{m},m,n)\equiv${\footnotesize
\begin{equation}
H_{m,n}(T_{m-1,n})+H_{m,n}^{\prime}(T_{m-1,n})\left(  \frac{X_{m}}{n}%
+\frac{X_{m}-\theta_{m}}{\sigma\sqrt{n}}\right)  +\frac{1}{2}H_{m,n}%
^{\prime\prime}(T_{m-1,n})\left(  \frac{X_{m}-\theta_{m}}{\sigma\sqrt{n}%
}\right)  ^{2}. \label{F}%
\end{equation}
} By Lemma \ref{lemma-ddp}, $\exists C>0$ such that (for all $m$ and $n$),
\[
\sup\limits_{m\leq n}\sup\limits_{x\in\mathbb{R}}|H_{m,n}^{\prime\prime
}(x)|\leq C\quad\text{ and }\quad{\sup\limits_{m\leq n}\sup\limits_{x,y\in
\mathbb{R},x\neq y}\frac{|H_{m,n}^{\prime\prime}(x)-H_{m,n}^{\prime\prime
}(y)|}{|x-y|}\leq C.}%
\]
}

{\normalsize By the Taylor expansion of {$H_{m,n}\in C_{b}^{2}(\mathbb{R})$},
$\forall$$\overline{\epsilon}>0$, $\exists$$\delta>0$ ($\delta$ depends only
on $C$ and $\overline{\epsilon}$), such that $\forall$$x,y\in\mathbb{R}$, and
all $n\geq m\geq1$, {\small
\begin{equation}
\left\vert H_{m,n}(x+y)-H_{m,n}(x)-H_{m,n}^{\prime}(x)y-\frac{1}{2}%
H_{m,n}^{\prime\prime}(x)y^{2}\right\vert \leq\overline{\epsilon}%
|y|^{2}I_{\{|y|<\delta\}}+{C}|y|^{2}I_{\{|y|\geq\delta\}}.\label{le0}%
\end{equation}
}Let $x=T_{m-1,n}$ and $y=\frac{X_{m}}{n}+\frac{X_{m}-\theta_{m}}{\sigma
\sqrt{n}}$ in (\ref{le0}), and obtain{\small
\begin{align*}
&  \sum_{m=1}^{n}\left\vert \sup\limits_{Q\in\mathcal{P}}E_{Q}\left[
H_{m,n}\left(  T_{m-1,n}+\frac{X_{m}}{n}+\frac{X_{m}-\theta_{m}}{\sigma
\sqrt{n}}\right)  \right]  -\sup\limits_{Q\in\mathcal{P}}E_{Q}\left[
F(\theta_{m},m,n)\right]  \right\vert \\
\leq &  \sum_{m=1}^{n}\sup\limits_{Q\in\mathcal{P}}E_{Q}\left[  \left\vert
H_{m,n}\left(  T_{m-1,n}+\frac{X_{m}}{n}+\frac{X_{m}-\theta_{m}}{\sigma
\sqrt{n}}\right)  -F(\theta_{m},m,n)\right\vert \right]  \\
\leq &  R_{1}(\overline{\epsilon},n)+R_{2}(C,n)+R_{3}(C,n)\text{, \ \ where}%
\end{align*}
}%
\begin{align*}
R_{1}(\overline{\epsilon},n) &  :=\overline{\epsilon}\sum\limits_{m=1}^{n}%
\sup\limits_{Q\in\mathcal{P}}E_{Q}\left[  \left\vert \tfrac{X_{m}}{n}%
+\tfrac{X_{m}-\theta_{m}}{\sigma\sqrt{n}}\right\vert ^{2}I_{\left\{
|\frac{X_{m}}{n}+\frac{X_{m}-\theta_{m}}{\sigma\sqrt{n}}|<\delta\right\}
}\right]  ,\\
R_{2}(C,n) &  :={C}\sum\limits_{m=1}^{n}\sup\limits_{Q\in\mathcal{P}}%
E_{Q}\left[  \left\vert \tfrac{X_{m}}{n}+\tfrac{X_{m}-\theta_{m}}{\sigma
\sqrt{n}}\right\vert ^{2}I_{\left\{  |\frac{X_{m}}{n}+\frac{X_{m}-\theta_{m}%
}{\sigma\sqrt{n}}|\geq\delta\right\}  }\right]  ,\\
R_{3}(C,n) &  :=\frac{C}{2}\sum\limits_{m=1}^{n}\sup\limits_{Q\in\mathcal{P}%
}E_{Q}\left[  \left\vert \tfrac{X_{m}}{n}\right\vert ^{2}+2\left\vert
\tfrac{X_{m}}{n}\right\vert \left\vert \tfrac{X_{m}-\theta_{m}}{\sigma\sqrt
{n}}\right\vert \right]  .
\end{align*}
It is readily proven that, for sufficiently large $n$,%
\begin{align*}
R_{1}(\overline{\epsilon},n)\leq &  \frac{2\overline{\epsilon}}{n^{2}}%
\sum\limits_{m=1}^{n}\sup\limits_{Q\in\mathcal{P}}E_{Q}\left[  \left\vert
X_{m}\right\vert ^{2}\right]  +\frac{2\overline{\epsilon}}{n}\sum
\limits_{m=1}^{n}\sup\limits_{Q\in\mathcal{P}}E_{Q}\left[  \left\vert
\tfrac{X_{m}-\theta_{m}}{\sigma}\right\vert ^{2}\right]  \\
\leq &  \frac{4\overline{\epsilon}}{n}\left(  \sigma^{2}+\kappa^{2}\right)
+\frac{4\overline{\epsilon}}{\sigma^{2}}\left(  \sigma^{2}+4\kappa^{2}\right)
,\\
R_{2}(C,n)\leq &  2C\left(  \frac{1}{n}+\frac{1}{\sigma\sqrt{n}}\right)
^{2}\sum\limits_{m=1}^{n}\sup\limits_{Q\in\mathcal{P}}E_{Q}\left[  \left\vert
X_{m}\right\vert ^{2}I_{\left\{  |\frac{X_{m}}{n}+\frac{X_{m}-\theta_{m}%
}{\sigma\sqrt{n}}|\geq\delta\right\}  }\right]  \\
&  +\frac{2C}{\sigma^{2}n}\sum\limits_{m=1}^{n}\sup\limits_{Q\in\mathcal{P}%
}E_{Q}\left[  \left\vert \theta_{m}\right\vert ^{2}I_{\left\{  |\frac{X_{m}%
}{n}+\frac{X_{m}-\theta_{m}}{\sigma\sqrt{n}}|\geq\delta\right\}  }\right]  \\
\leq &  \frac{2C}{\sigma^{2}}\frac{(\sigma+\sqrt{n})^{2}}{n^{2}}%
\sum\limits_{m=1}^{n}\sup\limits_{Q\in\mathcal{P}}E_{Q}\left[  \left\vert
X_{m}\right\vert ^{2}I_{\left\{  |X_{m}|>\frac{\sigma n}{\sigma+\sqrt{n}%
}\delta-\kappa\right\}  }\right]  \\
&  +\frac{2C}{\sigma^{2}n}\frac{\kappa^{2}}{\delta^{2}}\sum\limits_{m=1}%
^{n}\sup\limits_{Q\in\mathcal{P}}E_{Q}\left[  \left\vert \tfrac{X_{m}}%
{n}+\tfrac{X_{m}-\theta_{m}}{\sigma\sqrt{n}}\right\vert ^{2}\right]  ,\\
R_{3}(C,n)\leq &  \frac{C}{2n^{2}}\sum\limits_{m=1}^{n}\sup\limits_{Q\in
\mathcal{P}}E_{Q}\left[  \left\vert X_{m}\right\vert ^{2}\right]  +\frac
{C}{n^{3/2}\sigma}\sum\limits_{m=1}^{n}\sup\limits_{Q\in\mathcal{P}}%
E_{Q}\left[  \left\vert X_{m}\right\vert \left\vert X_{m}-\theta
_{m}\right\vert \right]  \\
\leq &  \left(  \frac{C}{n}+\frac{2C}{\sqrt{n}\sigma}\right)  \left(
\sigma^{2}+\kappa^{2}\right)  +\frac{2C\kappa}{\sqrt{n}\sigma}\sqrt{\sigma
^{2}+\kappa^{2}}\text{.}%
\end{align*}
By the finiteness of $\kappa,\sigma$ and the Lindeberg condition
(\ref{linder}),%
\[
\lim\limits_{\overline{\epsilon}\rightarrow0}\lim\limits_{n\rightarrow\infty
}\left(  R_{1}(\overline{\epsilon},n)+R_{2}(C,n)+R_{3}(C,n)\right)  =0\text{,}%
\]
which proves (\ref{remainder-esti}). }

{\normalsize \medskip}

{\normalsize \noindent\textbf{Step 2}: To prove (\ref{Equa}), it suffices to
prove that if we take $\theta_{m}$$=E_{Q}[X_{m}|\mathcal{G}_{m-1}]$ in
(\ref{F}), then%
\[
\sup\limits_{Q\in\mathcal{P}}E_{Q}\left[  F(\theta_{m},m,n)\right]
=\sup\limits_{Q\in\mathcal{P}}E_{Q}\left[  L_{m,n}(T_{m-1,n})\right]
,\;\forall n\geq m\geq1.
\]
In fact, if $\theta_{m}=$$E_{Q}[X_{m}|\mathcal{G}_{m-1}]$, then by a
generalization of Lemma \ref{lemma-rect}(iii) (in the proof of Theorem
\ref{thm-CLT}, we shall take $T_{m-1,n}=\frac{S_{m-1}}{n}+\frac{S_{m-1}^{Q}%
}{\sqrt{n}}$, in which case part (iv) of Lemma \ref{lemma-rect} suffices),
\begin{align*}
&  \sup\limits_{Q\in\mathcal{P}}E_{Q}\left[  F(\theta_{m},m,n)\right] \\
=  &  \sup\limits_{Q\in\mathcal{P}}E_{Q}\left[  H_{m,n}(T_{m-1,n}%
)+H_{m,n}^{\prime}(T_{m-1,n})\left(  \frac{X_{m}}{n}+\frac{Y_{m}^{Q}}{\sqrt
{n}}\right)  \right. \\
&  \left.  \qquad\qquad\qquad\qquad\qquad\qquad\qquad\qquad\qquad\qquad
+\frac{1}{2}H_{m,n}^{\prime\prime}(T_{m-1,n})\left(  \frac{Y_{m}^{Q}}{\sqrt
{n}}\right)  ^{2}\right] \\
=  &  \sup\limits_{Q\in\mathcal{P}}E_{Q}\left[  H_{m,n}(T_{m-1,n}%
)+H_{m,n}^{\prime}(T_{m-1,n})E_{Q}\left[  \left(  \frac{X_{m}}{n}+\frac
{Y_{m}^{Q}}{\sqrt{n}}\right)  |\mathcal{G}_{m-1}\right]  \right. \\
&  \left.  \qquad\qquad\qquad\qquad\qquad\qquad\qquad\quad+\frac{1}{2n}%
H_{m,n}^{\prime\prime}(T_{m-1,n})E_{Q}\left[  \left(  Y_{m}^{Q}\right)
^{2}|\mathcal{G}_{m-1}\right]  \right] \\
=  &  \sup\limits_{Q\in\mathcal{P}}E_{Q}\left[  H_{m,n}(T_{m-1,n})+\frac{1}%
{n}H_{m,n}^{\prime}(T_{m-1,n})X_{m}+\frac{1}{2n}H_{m,n}^{\prime\prime
}(T_{m-1,n})\right] \\
=  &  \sup\limits_{Q\in\mathcal{P}}E_{Q}\left[  H_{m,n}(T_{m-1,n})+\frac{1}%
{n}\;\mathbb{E}[H_{m,n}^{\prime}(T_{m-1,n})X_{m}|\mathcal{G}_{m-1}]+\frac
{1}{2n}H_{m,n}^{\prime\prime}(T_{m-1,n})\right] \\
=  &  \sup\limits_{Q\in\mathcal{P}}E_{Q}\left[  H_{m,n}(T_{m-1,n})+\frac{1}%
{n}g_{0}\left(  H_{m,n}^{\prime}(T_{m-1,n})\right)  +\frac{1}{2n}%
H_{m,n}^{\prime\prime}(T_{m-1,n})\right] \\
=  &  \sup\limits_{Q\in\mathcal{P}}E_{Q}\left[  L_{m,n}(T_{m-1,n})\right]  .
\end{align*}
Combine with (\ref{remainder-esti}) to complete the proof. \hfill
\hfill$\blacksquare$ }

{\normalsize \medskip}

{\normalsize The next three lemmas consider the special implications of
symmetry and thus relate to the proof of Theorem \ref{thm-special}. }

\begin{lemma}
[\cite{chenBSDE}]{\normalsize \label{lemma-sym0} Let $\varphi\in C_{b}%
^{3}(\mathbb{R})$ be symmetric with {center $c\in\mathbb{R}$}, and $v(t,x)$ be
the unique solution of Cauchy's problem for the parabolic equation
\begin{equation}
\left\{
\begin{array}
[c]{rl}%
\partial_{t}v(t,x) & =~\tfrac{1}{2}\partial_{xx}^{2}v(t,x)+g_{\epsilon
}(\partial_{x}v(t,x))\\
v(0,x) & =~\varphi(x).
\end{array}
\right.  \label{pde}%
\end{equation}
}

\begin{description}
\item {\normalsize (1) For any $t\geq0$, $v(t,\cdot)$ is symmetric with center
$c$. }

\item {\normalsize (2) If $sgn(\varphi^{\prime}(x))=-sgn(x-c)$, then, for any
$t\geq0$,
\[
sgn\left(  \partial_{x}v(t,x)\right)  =-sgn(x-c)\text{.}%
\]
}

\item {\normalsize (3) If $sgn(\varphi^{\prime}(x))=sgn(x-c)$, then, for any
$t\geq0$,
\[
sgn\left(  \partial_{x}v(t,x)\right)  =sgn(x-c)\text{.}%
\]
}
\end{description}
\end{lemma}

{\normalsize
}

\begin{lemma}
{\normalsize \label{lemma-sym} Let $\varphi\in C_{b}^{3}$$(\mathbb{R})$ be
symmetric with {center $c\in\mathbb{R}$}. Then the functions $\{H_{m,n}%
\}_{m=0}^{n}$ defined in (\ref{BSDEep-1}) satisfy, for any $n$ and
$m=0,\cdots,n$: }

\begin{description}
\item {\normalsize (1) $H_{m,n}$ is symmetric with center $c$. }

\item {\normalsize (2) If $sgn(\varphi^{\prime}(x))=-sgn(x-c)$, then
\[
sgn(H_{m,n}^{\prime}(x))=-sgn(x-c),\quad\text{and}\quad H_{m,n}^{\prime\prime
}(c)\leq0.
\]
}

\item {\normalsize (3) If $sgn(\varphi^{\prime}(x))=sgn(x-c)$, then
\[
sgn(H_{m,n}^{\prime}(x))=sgn(x-c),\quad\text{and}\quad H_{m,n}^{\prime\prime
}(c)\geq0.
\]
}
\end{description}
\end{lemma}

{\normalsize \noindent\textbf{Proof:} By the definition of $H_{m,n}(x)$ via
(\ref{BSDEep-1}) and the nonlinear Feynman-Kac formula, we know that
$H_{m,n}(x)=v(1-\tfrac{m}{n},x)$, where $v(t,x)$ is the solution of equation
(\ref{pde}). Then (1)-(3) follows from Lemma \ref{lemma-sym0}. \hfill
\hfill$\blacksquare$}

\begin{lemma}
{\normalsize \label{lemma-taylor2} Adopt the assumptions and notation in
(\ref{mumn}), (\ref{tmumn}) and Lemma \ref{lemma-taylor}, and let $\varphi\in
C_{b}^{3}$$(\mathbb{R})$ be symmetric with {center $c\in\mathbb{R}$}. }

\begin{description}
\item {\normalsize (1) If $sgn(\varphi^{\prime}(x))=-sgn(x-c)$, and if
$Y_{m}^{Q}$ in (\ref{Equa}) is replaced by $Z_{m}^{n}$, where
\[
Z_{m}^{n}=\frac{1}{\sigma}(X_{m}-\mu_{m}^{n}),\;\mu_{m}^{n}=\kappa
I_{A_{m-1,n}}-\kappa I_{A_{m-1,n}^{c}},\;A_{m-1,n}=\{T_{m-1,n}\leq c\},
\]
then
\begin{align}
\lim_{n\rightarrow\infty}\sum_{m=1}^{n}  &  \left\vert \sup\limits_{Q\in
\mathcal{P}}E_{Q}\left[  H_{m,n}\left(  T_{m-1,n}+\frac{X_{m}}{n}+\frac
{Z_{m}^{n}}{\sqrt{n}}\right)  \right]  \right. \nonumber\\
&  \qquad\qquad\qquad\qquad\left.  -\sup\limits_{Q\in\mathcal{P}}E_{Q}\left[
L_{m,n}\left(  T_{m-1,n}\right)  \right]  \right\vert =0. \label{sym-ty-eq}%
\end{align}
}

\item {\normalsize (2) If $sgn(\varphi^{\prime}(x))=sgn(x-c)$, and if
$Y_{m}^{Q}$ in (\ref{Equa}) is replaced by $\widetilde{Z}_{m}^{n}$, where{
\[
\widetilde{Z}_{m}^{n}=\frac{1}{\sigma}(X_{m}-\widetilde{\mu}_{m}%
^{n}),\;\widetilde{\mu}_{m}^{n}=\kappa I_{\widetilde{A}_{m-1,n}}-\kappa
I_{\widetilde{A}_{m-1,n}^{c}},\widetilde{A}_{m-1,n}=\{T_{m-1,n}\geq
c\}\text{,}%
\]
} then%
\begin{align}
\liminf_{n\rightarrow\infty}\sum_{m=1}^{n}  &  \left\{  \sup\limits_{Q\in
\mathcal{P}}E_{Q}\left[  H_{m,n}\left(  T_{m-1,n}+\frac{X_{m}}{n}%
+\frac{\widetilde{Z}_{m}^{n}}{\sqrt{n}}\right)  \right]  \right. \nonumber\\
&  \qquad\qquad\qquad\qquad-\sup\limits_{Q\in\mathcal{P}}E_{Q}\left[
L_{m,n}\left(  T_{m-1,n}\right)  \right]  \Bigg\}\geq0. \label{sym-ty-eq2-1}%
\end{align}
Furthermore, if
\begin{equation}
\lim_{\overline{\delta}\rightarrow0}\limsup_{n\rightarrow\infty}\frac{1}%
{n}\sum_{m=1}^{n}\sup_{Q\in\mathcal{P}}E_{Q}\left[  \left\vert E_{Q}%
[X_{m}|\mathcal{G}_{m-1}]-\widetilde{\mu}_{m}^{n}\right\vert I_{\{|T_{m-1,n}%
-c|\leq\overline{\delta}\}}\right]  =0, \label{lemma-condition-1}%
\end{equation}
then%
\begin{align}
\lim_{n\rightarrow\infty}\sum_{m=1}^{n}  &  \left\vert \sup\limits_{Q\in
\mathcal{P}}E_{Q}\left[  H_{m,n}\left(  T_{m-1,n}+\frac{X_{m}}{n}%
+\frac{\widetilde{Z}_{m}^{n}}{\sqrt{n}}\right)  \right]  \right. \nonumber\\
&  \qquad\qquad\qquad\qquad\quad-\sup\limits_{Q\in\mathcal{P}}E_{Q}\left[
L_{m,n}\left(  T_{m-1,n}\right)  \right]  \Bigg\vert=0. \label{sym-ty-eq2-2}%
\end{align}
}
\end{description}
\end{lemma}

\begin{remark}
{\normalsize The lemma is valid also if $c=\pm\infty$. Taking $c=+\infty$ in
(1) means that $\varphi$ is increasing on $\mathbb{R}$. Then $A_{m-1,n}%
=\Omega$ and $\mu_{m}^{n}=\kappa$ for any $1\leq m\leq n$. If $c=-\infty$ in
(1), then $\varphi$ is decreasing on $\mathbb{R}$, $A_{m-1,n}=\varnothing$ and
$\mu_{m}^{n}=-\kappa$. Similarly for (2). }
\end{remark}

{\normalsize \noindent\textbf{Proof of (1):} We proceed in two steps.  }

{\normalsize \noindent\textbf{Step 1}: Firstly, we prove
\begin{align}
\limsup_{n\rightarrow\infty}\sum_{m=1}^{n} &  \left\{  \sup\limits_{Q\in
\mathcal{P}}E_{Q}\left[  H_{m,n}\left(  T_{m-1,n}+\frac{X_{m}}{n}+\frac
{Z_{m}^{n}}{\sqrt{n}}\right)  \right]  \right.  \nonumber\\
&  \qquad\qquad\qquad\qquad\quad-\sup\limits_{Q\in\mathcal{P}}E_{Q}\left[
L_{m,n}\left(  T_{m-1,n}\right)  \right]  \bigg\}\leq0.\label{symproofstep1}%
\end{align}
By Lemma \ref{lemma-taylor}, we only need to prove the non-positivity of
\begin{align*}
\limsup_{n\rightarrow\infty}\sum_{m=1}^{n} &  \left\{  \sup\limits_{Q\in
\mathcal{P}}E_{Q}\left[  H_{m,n}\left(  T_{m-1,n}+\frac{X_{m}}{n}+\frac
{Z_{m}^{n}}{\sqrt{n}}\right)  \right]  \right.  \\
&  \left.  \quad-\sup\limits_{Q\in\mathcal{P}}E_{Q}\left[  H_{m,n}\left(
T_{m-1,n}+\frac{X_{m}}{n}+\frac{Y_{m}^{Q}}{\sqrt{n}}\right)  \right]
\right\}  .
\end{align*}
{For any $\overline{\delta}>0$,}
we set
\[
D_{m-1,n}^{\overline{\delta},1}=\{T_{m-1,n}>c+\overline{\delta}\},\ D_{m-1,n}%
^{\overline{\delta},2}=\{T_{m-1,n}<c-\overline{\delta}\}\text{,}%
\]%
\[
D_{m-1,n}^{\overline{\delta},3}=\{|T_{m-1,n}-c|\leq\overline{\delta
}\},\ N_{m,n}^{\overline{\delta}}=\left\{  \left\vert X_{m}\right\vert
\leq\tfrac{\sigma n}{\sigma+\sqrt{n}}\overline{\delta}-\kappa\right\}  .
\]
For any $\omega\in D_{m-1,n}^{\overline{\delta},1}\cap N_{m,n}^{\overline
{\delta}}$, we have
\begin{align*}
c &  <\left(  T_{m-1,n}+\frac{X_{m}}{n}+\frac{X_{m}-E_{Q}[X_{m}|\mathcal{G}%
_{m-1}]}{\sqrt{n}\sigma}\right)  (\omega)\\
&  \leq\left(  T_{m-1,n}+\frac{X_{m}}{n}+\frac{X_{m}+\kappa}{\sqrt{n}\sigma
}\right)  (\omega).
\end{align*}
By Lemma \ref{lemma-sym}, $H_{m,n}$ is decreasing on $(c,+\infty)$. Thus%
\begin{align*}
&  I_{D_{m-1,n}^{\overline{\delta},1}\cap N_{m,n}^{\overline{\delta}}}%
H_{m,n}\left(  T_{m-1,n}+\frac{X_{m}}{n}+\frac{Y_{m}^{Q}}{\sqrt{n}}\right)  \\
\geq &  I_{D_{m-1,n}^{\overline{\delta},1}\cap N_{m,n}^{\overline{\delta}}%
}H_{m,n}\left(  T_{m-1,n}+\frac{X_{m}}{n}+\frac{Z_{m}^{n}}{\sqrt{n}}\right)  .
\end{align*}
Also, for any $\omega\in D_{m-1,n}^{\overline{\delta},2}\cap N_{m,n}%
^{\overline{\delta}}$,%
\begin{align*}
c &  >\left(  T_{m-1,n}+\frac{X_{m}}{n}+\frac{X_{m}-E_{Q}[X_{m}|\mathcal{G}%
_{m-1}]}{\sqrt{n}\sigma}\right)  (\omega)\\
&  \geq\left(  T_{m-1,n}+\frac{X_{m}}{n}+\frac{X_{m}-\kappa}{\sqrt{n}\sigma
}\right)  (\omega)\text{.}%
\end{align*}
By Lemma \ref{lemma-sym}, $H_{m,n}$ is increasing on $(-\infty,c)$. Thus
\begin{align*}
&  I_{D_{m-1,n}^{\overline{\delta},2}\cap N_{m,n}^{\overline{\delta}}}%
H_{m,n}\left(  T_{m-1,n}+\frac{X_{m}}{n}+\frac{Y_{m}^{Q}}{\sqrt{n}}\right)  \\
\geq &  I_{D_{m-1,n}^{\overline{\delta},2}\cap N_{m,n}^{\overline{\delta}}%
}H_{m,n}\left(  T_{m-1,n}+\frac{X_{m}}{n}+\frac{Z_{m}^{n}}{\sqrt{n}}\right)
\text{.}%
\end{align*}
For $F(\theta_{m},m,n)$ defined in (\ref{F}), we have
\begin{align*}
&  H_{m,n}\left(  T_{m-1,n}+\frac{X_{m}}{n}+\frac{Z_{m}^{n}}{\sqrt{n}}\right)
-H_{m,n}\left(  T_{m-1,n}+\frac{X_{m}}{n}+\frac{Y_{m}^{Q}}{\sqrt{n}}\right)
\\
= &  F(\mu_{m}^{n},m,n)-F\left(  E_{Q}[X_{m}|\mathcal{G}_{m-1}],m,n\right)  \\
&  +H_{m,n}\left(  T_{m-1,n}+\frac{X_{m}}{n}+\frac{Z_{m}^{n}}{\sqrt{n}%
}\right)  -F(\mu_{m}^{n},m,n)\\
&  -\left(  H_{m,n}\left(  T_{m-1,n}+\frac{X_{m}}{n}+\frac{Y_{m}^{Q}}{\sqrt
{n}}\right)  -F\left(  E_{Q}[X_{m}|\mathcal{G}_{m-1}],m,n\right)  \right)  .
\end{align*}
Therefore,{\small
\begin{align*}
&  \sum_{m=1}^{n}\left\{  \sup\limits_{Q\in\mathcal{P}}E_{Q}\left[
H_{m,n}\left(  T_{m-1,n}+\frac{X_{m}}{n}+\frac{Z_{m}^{n}}{\sqrt{n}}\right)
\right]  \right.  \\
&  \left.  \qquad\qquad\qquad\qquad\qquad\qquad\qquad-\sup\limits_{Q\in
\mathcal{P}}E_{Q}\left[  H_{m,n}\left(  T_{m-1,n}+\frac{X_{m}}{n}+\frac
{Y_{m}^{Q}}{\sqrt{n}}\right)  \right]  \right\}  \\
\leq &  \sum_{m=1}^{n}\sup\limits_{Q\in\mathcal{P}}E_{Q}\left[  H_{m,n}\left(
T_{m-1,n}+\frac{X_{m}}{n}+\frac{Z_{m}^{n}}{\sqrt{n}}\right)  -H_{m,n}\left(
T_{m-1,n}+\frac{X_{m}}{n}+\frac{Y_{m}^{Q}}{\sqrt{n}}\right)  \right]  \\
\leq &  \sum_{m=1}^{n}\sup\limits_{Q\in\mathcal{P}}E_{Q}\left[  I_{D_{m-1,n}%
^{\overline{\delta},3}}\left[  H_{m,n}\left(  T_{m-1,n}+\frac{X_{m}}{n}%
+\frac{Z_{m}^{n}}{\sqrt{n}}\right)  \right.  \right.  \\
&  \left.  \left.  \qquad\qquad\qquad\qquad\qquad\qquad\qquad\qquad
\qquad-H_{m,n}\left(  T_{m-1,n}+\frac{X_{m}}{n}+\frac{Y_{m}^{Q}}{\sqrt{n}%
}\right)  \right]  \right]  \\
&  +\sum_{m=1}^{n}\sup\limits_{Q\in\mathcal{P}}E_{Q}\left[  I_{\{\left\vert
X_{m}\right\vert >\tfrac{\sigma n}{\sigma+\sqrt{n}}\overline{\delta}-\kappa
\}}\left\vert H_{m,n}\left(  T_{m-1,n}+\frac{X_{m}}{n}+\frac{Z_{m}^{n}}%
{\sqrt{n}}\right)  \right.  \right.  \\
&  \left.  \left.  \qquad\qquad\qquad\qquad\qquad\qquad\qquad\qquad
\qquad-H_{m,n}\left(  T_{m-1,n}+\frac{X_{m}}{n}+\frac{Y_{m}^{Q}}{\sqrt{n}%
}\right)  \right\vert \right]  \\
\leq &  I_{n}^{1}+I_{n}^{2}+I_{n}^{3},
\end{align*}
} where $I_{n}^{1},I_{n}^{2},I_{n}^{3}$ are defined by{\small
\begin{align*}
I_{n}^{1}\equiv &  \sum_{m=1}^{n}\sup\limits_{Q\in\mathcal{P}}E_{Q}\left[
I_{D_{m-1,n}^{\overline{\delta},3}}\left(  F(\mu_{m}^{n},m,n)-F\left(
E_{Q}[X_{m}|\mathcal{G}_{m-1}],m,n\right)  \right)  \right]  ,\\
I_{n}^{2}\equiv &  \sum_{m=1}^{n}\sup\limits_{Q\in\mathcal{P}}E_{Q}\left[
2\Vert\varphi\Vert I_{\{\left\vert X_{m}\right\vert >\tfrac{\sigma n}%
{\sigma+\sqrt{n}}\overline{\delta}-\kappa\}}\right]  \qquad(\text{ where
}\Vert\varphi\Vert=\sup_{x\in\mathbb{R}}\varphi(x)),\\
I_{n}^{3}\equiv &  \sum_{m=1}^{n}\sup\limits_{Q\in\mathcal{P}}E_{Q}\left[
\left\vert H_{m,n}\left(  T_{m-1,n}+\frac{X_{m}}{n}+\frac{Z_{m}^{n}}{\sqrt{n}%
}\right)  -F(\mu_{m}^{n},m,n)\right\vert \right]  \\
&  +\sum_{m=1}^{n}\sup\limits_{Q\in\mathcal{P}}E_{Q}\left[  \left\vert
H_{m,n}\left(  T_{m-1,n}+\frac{X_{m}}{n}+\frac{Y_{m}^{Q}}{\sqrt{n}}\right)
-F\left(  E_{Q}[X_{m}|\mathcal{G}_{m-1}],m,n\right)  \right\vert \right]  .
\end{align*}
}{{\small By the Lindeberg condition (\ref{linder}), $\lim
\limits_{n\rightarrow\infty}I_{n}^{2}=0$, and by the remainder estimate in the
proof of Lemma \ref{lemma-taylor}, $\lim\limits_{n\rightarrow\infty}I_{n}%
^{3}=0$. Thus it suffices to show that $\lim\limits_{n\rightarrow\infty}%
I_{n}^{1}=0$, which is proven as follows. }}}
{\normalsize {From Lemma \ref{lemma-sym}(2), $H_{m,n}^{\prime\prime}(c)\leq0$.
Therefore,{\small
\begin{align*}
I_{n}^{1}= &  \sum_{m=1}^{n}\sup\limits_{Q\in\mathcal{P}}E_{Q}\left[
I_{D_{m-1,n}^{\overline{\delta},3}}\left(  F(\mu_{m}^{n},m,n)-F\left(
E_{Q}[X_{m}|\mathcal{G}_{m-1}],m,n\right)  \right)  \right]  \\
\leq &  \sum_{m=1}^{n}\sup\limits_{Q\in\mathcal{P}}E_{Q}\left[  I_{D_{m-1,n}%
^{\overline{\delta},3}}H_{m,n}^{\prime\prime}(T_{m-1,n})\frac{(E_{Q}%
[X_{m}|\mathcal{G}_{m-1}]-\mu_{m}^{n})^{2}}{2n\sigma^{2}}\right]  \\
\leq &  \sum_{m=1}^{n}\sup\limits_{Q\in\mathcal{P}}E_{Q}\left[  I_{D_{m-1,n}%
^{\overline{\delta},3}}\left(  H_{m,n}^{\prime\prime}(T_{m-1,n})-H_{m,n}%
^{\prime\prime}(c))\right)  \frac{(E_{Q}[X_{m}|\mathcal{G}_{m-1}]-\mu_{m}%
^{n})^{2}}{2n\sigma^{2}}\right]  \\
&  +\sum_{m=1}^{n}\sup\limits_{Q\in\mathcal{P}}E_{Q}\left[  I_{D_{m-1,n}%
^{\overline{\delta},3}}H_{m,n}^{\prime\prime}(c)\frac{(E_{Q}[X_{m}%
|\mathcal{G}_{m-1}]-\mu_{m}^{n})^{2}}{2n\sigma^{2}}\right]  \\
\leq &  \sum_{m=1}^{n}\sup\limits_{Q\in\mathcal{P}}E_{Q}\left[  I_{D_{m-1,n}%
^{\overline{\delta},3}}L|T_{m-1,n}-c|\frac{(E_{Q}[X_{m}|\mathcal{G}_{m-1}%
]-\mu_{m}^{n})^{2}}{2n\sigma^{2}}\right]  \leq\frac{L2\kappa^{2}}{\sigma^{2}%
}\overline{\delta},
\end{align*}
}}}where $L$ is the uniform Lipschitz constant for $H_{m,n}^{\prime\prime}$
given in Lemma \ref{lemma-ddp}.

{\normalsize
}

{\normalsize \medskip\smallskip\  }

{\normalsize \noindent\textbf{Step 2}: Next we prove
\begin{align}
\liminf_{n\rightarrow\infty}\sum_{m=1}^{n}  &  \left\{  \sup\limits_{Q\in
\mathcal{P}}E_{Q}\left[  H_{m,n}\left(  T_{m-1,n}+\frac{X_{m}}{n}+\frac
{Z_{m}^{n}}{\sqrt{n}}\right)  \right]  \right. \nonumber\\
&  \qquad\qquad\qquad\qquad-\sup\limits_{Q\in\mathcal{P}}E_{Q}\left[
L_{m,n}\left(  T_{m-1,n}\right)  \right]  \bigg\} \geq0. \label{symproofstep2}%
\end{align}
}

{\normalsize By the remainder estimate in the proof of Lemma
\ref{lemma-taylor}, it suffices to take $\theta_{m}=\mu_{m}^{n}$ in (\ref{F})
and to show that
\[
\sup\limits_{Q\in\mathcal{P}}E_{Q}\left[  F(\theta_{m},m,n)\right]  \geq
\sup\limits_{Q\in\mathcal{P}}E_{Q}\left[  L_{m,n}\left(  T_{m-1,n}\right)
\right]  ,\quad n\geq m\geq1.
\]
}

{\normalsize
}

{\normalsize \noindent By Lemma \ref{Lem_increasing_seq}, there exist
$\{Q_{k}\}_{k\geq1},\{\overline{P}_{j}^{m}\}_{j\geq1},\{\underline{P}_{j}%
^{m}\}_{j\geq1}\subset$$\mathcal{P}$ such that
\begin{align*}
\lim_{j\rightarrow\infty}E_{\overline{P}_{j}^{m}}[X_{m}|\mathcal{G}%
_{m-1}]=\kappa,\ \lim_{j\rightarrow\infty}E_{\underline{P}_{j}^{m}}%
[X_{m}|\mathcal{G}_{m-1}]=-\kappa,
\end{align*}
and
\begin{align*}
&  \sup\limits_{Q\in\mathcal{P}}E_{Q}\left[  H_{m,n}(T_{m-1,n})+\frac{1}%
{n}g_{0}\left(  H_{m,n}^{\prime}(T_{m-1,n})\right)  +\frac{1}{2n}%
H_{m,n}^{\prime\prime}(T_{m-1,n})\right] \\
=  &  \lim_{k\rightarrow\infty}E_{Q_{k}}\left[  H_{m,n}(T_{m-1,n})+\frac{1}%
{n}g_{0}\left(  H_{m,n}^{\prime}(T_{m-1,n})\right)  +\frac{1}{2n}%
H_{m,n}^{\prime\prime}(T_{m-1,n})\right]  .
\end{align*}
By Lemma \ref{lemma-rect}(ii), there exist $\{R_{j}^{m}\}_{j\geq1}$$\subset
$$\mathcal{P}$ satisfying
\[
E_{R_{j}^{m}}[X_{m}|\mathcal{G}_{m-1}]=I_{A_{m-1,n}}E_{\overline{P}_{j}^{m}%
}[X_{m}|\mathcal{G}_{m-1}]+I_{A_{m-1,n}^{c}}E_{\underline{P}_{j}^{m}}%
[X_{m}|\mathcal{G}_{m-1}].
\]
By Lemma \ref{lemma-rect}(iii) and the dominated convergence theorem,%
\begin{align*}
&  \sup\limits_{Q\in\mathcal{P}}E_{Q}\left[  L_{m,n}(T_{m-1,n})\right] \\
=  &  \sup\limits_{Q\in\mathcal{P}}E_{Q}\left[  H_{m,n}(T_{m-1,n})+\frac{1}%
{n}g_{0}\left(  H_{m,n}^{\prime}(T_{m-1,n})\right)  +\frac{1}{2n}%
H_{m,n}^{\prime\prime}(T_{m-1,n})\right] \\
=  &  \lim_{k\rightarrow\infty}E_{Q_{k}}\left[  H_{m,n}(T_{m-1,n})+\frac{1}%
{n}g_{0}\left(  H_{m,n}^{\prime}(T_{m-1,n})\right)  +\frac{1}{2n}%
H_{m,n}^{\prime\prime}(T_{m-1,n})\right] \\
=  &  \lim_{k\rightarrow\infty}E_{Q_{k}}\left[  H_{m,n}(T_{m-1,n}%
)+I_{A_{m-1,n}}H_{m,n}^{\prime}(T_{m-1,n})\frac{\kappa}{n}\right. \\
&  \qquad\qquad\qquad\qquad\qquad\left.  +I_{A_{m-1,n}^{c}}H_{m,n}^{\prime
}(T_{m-1,n})\frac{-\kappa}{n}+\frac{1}{2n}H_{m,n}^{\prime\prime}%
(T_{m-1,n})\right] \\
=  &  \lim_{k\rightarrow\infty}E_{Q_{k}}\Bigg[H_{m,n}(T_{m-1,n})\\
&  \qquad\qquad\quad+\lim_{j\rightarrow\infty}E_{\overline{P}_{j}^{m}}\left[
I_{A_{m-1,n}} H_{m,n}^{\prime}(T_{m-1,n})\left(  \frac{X_{m}}{n}+\frac
{X_{m}-\kappa}{\sigma\sqrt{n}}\right)  |\mathcal{G}_{m-1}\right] \\
&  \qquad\qquad\quad+\lim_{j\rightarrow\infty}E_{\overline{P}_{j}^{m}}\left[
I_{A_{m-1,n}} H_{m,n}^{\prime\prime}(T_{m-1,n})\frac{(X_{m}-\kappa)^{2}%
}{2n\sigma^{2}} |\mathcal{G}_{m-1}\right] \\
&  \qquad\qquad\quad+\lim_{j\rightarrow\infty}E_{\underline{P}_{j}^{m}}\left[
I_{A_{m-1,n}^{c}} H_{m,n}^{\prime}(T_{m-1,n})\left(  \frac{X_{m}}{n}%
+\frac{X_{m}+\kappa}{\sigma\sqrt{n}}\right)  |\mathcal{G}_{m-1}\right] \\
&  \qquad\qquad\quad+\lim_{j\rightarrow\infty}E_{\underline{P}_{j}^{m}}\left[
I_{A_{m-1,n}^{c}} H_{m,n}^{\prime\prime}(T_{m-1,n})\frac{(X_{m}+\kappa)^{2}%
}{2n\sigma^{2}} |\mathcal{G}_{m-1}\right]  \Bigg]\\
=  &  \lim_{k\rightarrow\infty}\lim_{j\rightarrow\infty}E_{Q_{k}}\left[
H_{m,n}(T_{m-1,n})+E_{R_{j}^{m}}\left[  H_{m,n}^{\prime}(T_{m-1,n})\left(
\frac{X_{m}}{n}+\frac{Z_{m}^{n}}{\sqrt{n}}\right)  \right.  \right. \\
&  \qquad\qquad\qquad\qquad\qquad\qquad\qquad\qquad\quad\left.  \left.
+\frac{1}{2n}H_{m,n}^{\prime\prime}(T_{m-1,n})(Z_{m}^{n})^{2}|\mathcal{G}%
_{m-1}\right]  \right] \\
\leq &  \sup\limits_{Q\in\mathcal{P}}E_{Q}\left[  H_{m,n}(T_{m-1,n}%
)+ess\sup_{R\in\mathcal{P}}E_{R}\left[  H_{m,n}^{\prime}(T_{m-1,n})\left(
\frac{X_{m}}{n}+\frac{Z_{m}^{n}}{\sqrt{n}}\right)  \right.  \right. \\
&  \left.  \left.  \qquad\qquad\qquad\qquad\qquad\qquad\qquad\qquad\quad
+\frac{1}{2n}H_{m,n}^{\prime\prime}(T_{m-1,n})(Z_{m}^{n})^{2}|\mathcal{G}%
_{m-1}\right]  \right] \\
=  &  \sup\limits_{Q\in\mathcal{P}}E_{Q}\left[  H_{m,n}(T_{m-1,n}%
)+H_{m,n}^{\prime}(T_{m-1,n})\left(  \frac{X_{m}}{n}+\frac{Z_{m}^{n}}{\sqrt
{n}}\right)  \right. \\
&  \left.  \qquad\qquad\qquad\qquad\qquad\qquad\qquad\qquad\qquad\quad
\quad+\frac{1}{2n}H_{m,n}^{\prime\prime}(T_{m-1,n})(Z_{m}^{n})^{2}\right]
\text{.}%
\end{align*}
Combined with (\ref{remainder-esti}), this implies (\ref{symproofstep2}), thus
completing the proof of (1). \noindent}

{\normalsize \bigskip}

{\normalsize \noindent\textbf{Proof of (2): }\ Proof of inequality
(\ref{sym-ty-eq2-1}) is similar to that of (\ref{symproofstep2}).  }

{\normalsize To prove (\ref{sym-ty-eq2-2}), assuming (\ref{lemma-condition-1}%
), we need only prove
\begin{align}
\limsup_{n\rightarrow\infty}\sum_{m=1}^{n}  &  \left\{  \sup\limits_{Q\in
\mathcal{P}}E_{Q}\left[  H_{m,n}\left(  T_{m-1,n}+\frac{X_{m}}{n}%
+\frac{\widetilde{Z}_{m}^{n}}{\sqrt{n}}\right)  \right]  \right. \nonumber\\
&  \qquad\qquad\qquad\qquad\qquad-\sup\limits_{Q\in\mathcal{P}}E_{Q}\left[
L_{m,n}\left(  T_{m-1,n}\right)  \right]  \Bigg\} \leq0.
\label{symproof-2-step0}%
\end{align}
By assumption (\ref{lemma-condition-1}), $\forall\varepsilon>0$,
$\exists\overline{\delta}_{\varepsilon}>0$ such that
\[
\limsup_{n\rightarrow\infty}\frac{1}{n}\sum_{m=1}^{n}\sup_{Q\in\mathcal{P}%
}E_{Q}\left[  \left\vert E_{Q}[X_{m}|\mathcal{G}_{m-1}]-\widetilde{\mu}%
_{m}^{n}\right\vert I_{\{|T_{m-1,n}-c|\leq\overline{\delta}_{\varepsilon}%
\}}\right]  \leq\varepsilon\text{.}%
\]
By Lemma \ref{lemma-taylor}, we only need to prove the non-positivity of
\begin{align*}
\limsup_{n\rightarrow\infty}\sum_{m=1}^{n}  &  \left\{  \sup\limits_{Q\in
\mathcal{P}}E_{Q}\left[  H_{m,n}\left(  T_{m-1,n}+\frac{X_{m}}{n}%
+\frac{\widetilde{Z}_{m}^{n}}{\sqrt{n}}\right)  \right]  \right. \\
&  \left.  \qquad\qquad-\sup\limits_{Q\in\mathcal{P}}E_{Q}\left[
H_{m,n}\left(  T_{m-1,n}+\frac{X_{m}}{n}+\frac{Y_{m}^{Q}}{\sqrt{n}}\right)
\right]  \right\}  .
\end{align*}
}

{\normalsize Define%
\[
D_{m-1,n}^{\overline{\delta}_{\varepsilon},1}=\{T_{m-1,n}>c+\overline{\delta
}_{\varepsilon}\},\ D_{m-1,n}^{\overline{\delta}_{\varepsilon},2}%
=\{T_{m-1,n}<c-\overline{\delta}_{\varepsilon}\}\text{,}%
\]%
\[
D_{m-1,n}^{\overline{\delta}_{\varepsilon},3}=\{|T_{m-1,n}-c|\leq
\overline{\delta}_{\varepsilon}\}\text{,}\ N_{m,n}^{\overline{\delta
}_{\varepsilon}}=\left\{  \left\vert X_{m}\right\vert \leq\tfrac{\sigma
n}{\sigma+\sqrt{n}}\overline{\delta}_{\varepsilon}-\kappa\right\}  .
\]
For any $\omega\in D_{m-1,n}^{\overline{\delta}_{\varepsilon},1}\cap
N_{m,n}^{\overline{\delta}_{\varepsilon}}$,
\begin{align*}
c &  <\left(  T_{m-1,n}+\frac{X_{m}}{n}+\frac{X_{m}-\kappa}{\sqrt{n}\sigma
}\right)  (\omega)\\
&  \leq\left(  T_{m-1,n}+\frac{X_{m}}{n}+\frac{X_{m}-E_{Q}[X_{m}%
|\mathcal{G}_{m-1}]}{\sqrt{n}\sigma}\right)  (\omega).
\end{align*}
By Lemma \ref{lemma-sym}, $H_{m,n}$ is increasing on $(c,+\infty)$. Thus
\begin{align*}
&  I_{D_{m-1,n}^{\overline{\delta}_{\varepsilon},1}\cap N_{m,n}^{\overline
{\delta}_{\varepsilon}}}H_{m,n}\left(  T_{m-1,n}+\frac{X_{m}}{n}+\frac
{Y_{m}^{Q}}{\sqrt{n}}\right)  \\
\geq &  I_{D_{m-1,n}^{\overline{\delta}_{\varepsilon},1}\cap N_{m,n}%
^{\overline{\delta}_{\varepsilon}}}H_{m,n}\left(  T_{m-1,n}+\frac{X_{m}}%
{n}+\frac{\widetilde{Z}_{m}^{n}}{\sqrt{n}}\right)  .
\end{align*}
Also, for any $\omega\in D_{m-1,n}^{\overline{\delta}_{\varepsilon},2}\cap
N_{m,n}^{\overline{\delta}_{\varepsilon}}$,%
\begin{align*}
c &  >\left(  T_{m-1,n}+\frac{X_{m}}{n}+\frac{X_{m}+\kappa}{\sqrt{n}\sigma
}\right)  (\omega)\\
&  \geq\left(  T_{m-1,n}+\frac{X_{m}}{n}+\frac{X_{m}-E_{Q}[X_{m}%
|\mathcal{G}_{m-1}]}{\sqrt{n}\sigma}\right)  (\omega)\text{.}%
\end{align*}
By Lemma \ref{lemma-sym}, $H_{m,n}$ is decreasing on $(-\infty,c)$. Thus
\begin{align*}
&  I_{D_{m-1,n}^{\overline{\delta}_{\varepsilon},2}\cap N_{m,n}^{\overline
{\delta}_{\varepsilon}}}H_{m,n}\left(  T_{m-1,n}+\frac{X_{m}}{n}+\frac
{Y_{m}^{Q}}{\sqrt{n}}\right)  \\
\geq &  I_{D_{m-1,n}^{\overline{\delta}_{\varepsilon},2}\cap N_{m,n}%
^{\overline{\delta}_{\varepsilon}}}H_{m,n}\left(  T_{m-1,n}+\frac{X_{m}}%
{n}+\frac{\widetilde{Z}_{m}^{n}}{\sqrt{n}}\right)  \text{.}%
\end{align*}
Therefore,{
\begin{align*}
&  \sum_{m=1}^{n}\left\{  \sup\limits_{Q\in\mathcal{P}}E_{Q}\left[
H_{m,n}\left(  T_{m-1,n}+\frac{X_{m}}{n}+\frac{\widetilde{Z}_{m}^{n}}{\sqrt
{n}}\right)  \right]  \right.  \\
&  \left.  \qquad\qquad\qquad\qquad\qquad-\sup\limits_{Q\in\mathcal{P}}%
E_{Q}\left[  H_{m,n}\left(  T_{m-1,n}+\frac{X_{m}}{n}+\frac{Y_{m}^{Q}}%
{\sqrt{n}}\right)  \right]  \right\}  \\
\leq &  \sum_{m=1}^{n}\sup\limits_{Q\in\mathcal{P}}E_{Q}\left[  H_{m,n}\left(
T_{m-1,n}+\frac{X_{m}}{n}+\frac{\widetilde{Z}_{m}^{n}}{\sqrt{n}}\right)
\right.  \\
&  \qquad\qquad\qquad\qquad\qquad\qquad\qquad\ -H_{m,n}\left(  T_{m-1,n}%
+\frac{X_{m}}{n}+\frac{Y_{m}^{Q}}{\sqrt{n}}\right)  \Bigg]\\
\leq &  \sum_{m=1}^{n}\sup\limits_{Q\in\mathcal{P}}E_{Q}\left[  I_{D_{m-1,n}%
^{\overline{\delta}_{\varepsilon},3}}\left[  H_{m,n}\left(  T_{m-1,n}%
+\frac{X_{m}}{n}+\frac{\widetilde{Z}_{m}^{n}}{\sqrt{n}}\right)  \right.
\right.  \\
&  \qquad\qquad\qquad\qquad\qquad\qquad\quad\quad-H_{m,n}\left(
T_{m-1,n}+\frac{X_{m}}{n}+\frac{Y_{m}^{Q}}{\sqrt{n}}\right)  \bigg]\Bigg]\\
&  +\sum_{m=1}^{n}\sup\limits_{Q\in\mathcal{P}}E_{Q}\left[  I_{\{\left\vert
X_{m}\right\vert >\tfrac{\sigma n}{\sigma+\sqrt{n}}\overline{\delta
}_{\varepsilon}-\kappa\}}\left\vert H_{m,n}\left(  T_{m-1,n}+\frac{X_{m}}%
{n}+\frac{\widetilde{Z}_{m}^{n}}{\sqrt{n}}\right)  \right.  \right.  \\
&  \qquad\qquad\qquad\qquad\qquad\qquad\qquad\ -H_{m,n}\left(  T_{m-1,n}%
+\frac{X_{m}}{n}+\frac{Y_{m}^{Q}}{\sqrt{n}}\right)  \bigg\vert\Bigg]\\
\leq &  \widetilde{I}_{n}^{1}+\widetilde{I}_{n}^{2}+\widetilde{I}_{n}^{3},
\end{align*}
} where $\widetilde{I}_{n}^{1},\widetilde{I}_{n}^{2},\widetilde{I}_{n}^{3}$
are defined by{\small
\begin{align*}
\widetilde{I}_{n}^{1}\equiv &  \sum_{m=1}^{n}\sup\limits_{Q\in\mathcal{P}%
}E_{Q}\left[  {I_{D_{m-1,n}^{\overline{\delta}_{\varepsilon},3}}}\left(
F(\widetilde{\mu}_{m}^{n},m,n)-F\left(  E_{Q}[X_{m}|\mathcal{G}_{m-1}%
],m,n\right)  \right)  \right]  ,\\
\widetilde{I}_{n}^{2}\equiv &  \sum_{m=1}^{n}\sup\limits_{Q\in\mathcal{P}%
}E_{Q}\left[  2\Vert\varphi\Vert I_{\{\left\vert X_{m}\right\vert
>\tfrac{\sigma n}{\sigma+\sqrt{n}}\overline{\delta}_{\varepsilon}-\kappa
\}}\right]  \qquad(\text{ where }\Vert\varphi\Vert=\sup_{x\in\mathbb{R}%
}\varphi(x)),\\
\widetilde{I}_{n}^{3}\equiv &  \sum_{m=1}^{n}\sup\limits_{Q\in\mathcal{P}%
}E_{Q}\left[  \left\vert H_{m,n}\left(  T_{m-1,n}+\frac{X_{m}}{n}%
+\frac{\widetilde{Z}_{m}^{n}}{\sqrt{n}}\right)  -F(\widetilde{\mu}_{m}%
^{n},m,n)\right\vert \right]  \\
&  +\sum_{m=1}^{n}\sup\limits_{Q\in\mathcal{P}}E_{Q}\left[  \left\vert
H_{m,n}\left(  T_{m-1,n}+\frac{X_{m}}{n}+\frac{Y_{m}^{Q}}{\sqrt{n}}\right)
-F\left(  E_{Q}[X_{m}|\mathcal{G}_{m-1}],m,n\right)  \right\vert \right]  .
\end{align*}
} By the Lindeberg condition (\ref{linder}), $\lim\limits_{n\rightarrow\infty
}\widetilde{I}_{n}^{2}=0$, and by the remainder estimate in the proof of Lemma
\ref{lemma-taylor}, $\lim\limits_{n\rightarrow\infty}\widetilde{I}_{n}^{3}=0$.
Finally, we prove that $\lim\limits_{n\rightarrow\infty}\widetilde{I}_{n}%
^{3}=0$:
\begin{align*}
\widetilde{I}_{n}^{1}= &  \sum_{m=1}^{n}\sup\limits_{Q\in\mathcal{P}}%
E_{Q}\left[  {I_{D_{m-1,n}^{\overline{\delta}_{\varepsilon},3}}}\left(
F(\widetilde{\mu}_{m}^{n},m,n)-F\left(  E_{Q}[X_{m}|\mathcal{G}_{m-1}%
],m,n\right)  \right)  \right]  \\
\leq &  \frac{1}{2n\sigma^{2}}\sum_{m=1}^{n}\sup\limits_{Q\in\mathcal{P}}%
E_{Q}\left[  H_{m,n}^{\prime\prime}(T_{m-1,n})(E_{Q}[X_{m}|\mathcal{G}%
_{m-1}]-\widetilde{\mu}_{m}^{n})^{2}I_{\{|T_{m-1,n}-c|\leq\overline{\delta
}_{\varepsilon}\}}\right]  \\
\leq &  {\frac{C\kappa}{n\sigma^{2}}}\sum_{m=1}^{n}\sup\limits_{Q\in
\mathcal{P}}E_{Q}\left[  \left\vert E_{Q}[X_{m}|\mathcal{G}_{m-1}%
]-\widetilde{\mu}_{m}^{n}\right\vert I_{\{|T_{m-1,n}-c|\leq\overline{\delta
}_{\varepsilon}\}}\right]  ,
\end{align*}
where $C$ is the uniform bound given in Lemma \ref{lemma-ddp} (2). Thus
$\limsup_{n\rightarrow\infty}\widetilde{I}_{n}^{1}\leq \frac{C\kappa\epsilon}{\sigma^2}$,
where $\varepsilon$ is arbitrary. This proves (\ref{symproof-2-step0}) and
completes the proof of part (2).\hfill\hfill$\blacksquare$ }

{\normalsize \medskip\  }

{\normalsize The next lemma is used in extending the two theorems from the
special case (\ref{kappa}) to general \underline{$\mu$} and $\overline{\mu}$.
}

\begin{lemma}
{\normalsize \label{lemma-shift}For any $\kappa>0$ and $c\in\mathbb{R},$
\[
\mathbb{E}_{\left[  -\kappa,\kappa\right]  }\left[  \varphi\left(
c+B_{1}\right)  \right]  \newline=\mathbb{E}_{\left[  -\kappa+c,\;\kappa
+c\right]  }\left[  \varphi\left(  B_{1}\right)  \right]  ,
\]
where $\mathbb{E}_{\left[  -\kappa+c,\;\kappa+c\right]  }\left[
\varphi\left(  B_{1}\right)  \right]  $ is defined in (\ref{BSDE}).  }
\end{lemma}

{\normalsize \noindent\textbf{Proof:} $\mathbb{E}_{[-\kappa+c,\kappa
+c]}[\varphi(B_{1})]=Y_{0}^{c}$, where $(Y_{t}^{c},Z_{t}^{c})$ solves
\begin{align*}
Y_{t}^{c}=  &  \varphi(B_{1})+\int_{t}^{1}\max\limits_{-\kappa+c\leq\mu
\leq\kappa+c}(\mu Z_{s}^{c})ds-\int_{t}^{1}Z_{s}^{c}dB_{s}\\
=  &  \varphi(B_{1})+\int_{t}^{1}\left(  \max\limits_{-\kappa\leq\mu\leq
\kappa}(\mu Z_{s}^{c})+cZ_{s}^{c}\right)  ds-\int_{t}^{1}Z_{s}^{c}%
dB_{s},\ 0\leq t\leq1\text{.}%
\end{align*}
where the last equality is due to
\[
\max\limits_{-\kappa+c\leq\mu\leq\kappa+c}(\mu z)=(\kappa+c)z^{+}%
-(-\kappa+c)z^{-}=\max\limits_{-\kappa\leq\mu\leq\kappa}(\mu z)+cz.
\]
Let $Q^{c}$ be the probability measure satisfying{
\[
E_{P^{\ast}}\left[  \frac{dQ^{c}}{dP^{\ast}}|\mathcal{F}_{t}\right]
=\exp\left\{  -\frac{c^{2}t}{2}+cB_{t}\right\}  ,\quad t\geq0\text{.}%
\]
} Then $W_{t}=B_{t}-ct$ is a Brownian motion under $Q^{c}$ and $(Y_{t}%
^{c},Z_{t}^{c})$ solves%
\[
Y_{t}^{c}=\varphi(c+W_{1})+\int_{t}^{1}{\max\limits_{-\kappa\leq\mu\leq\kappa
}}(\mu Z_{s}^{c})ds-\int_{t}^{1}Z_{s}^{c}dW_{s},\ 0\leq t\leq1\text{.}%
\]
Hence $Y_{0}^{c}=\mathbb{E}_{\left[  -\kappa,\kappa\right]  }\left[
\varphi\left(  c+B_{1}\right)  \right]  $.\hfill\hfill$\blacksquare$  }

{\normalsize \medskip}

{\normalsize Chen et al \cite{chenBSDE} derive closed-form solutions for a
class of BSDEs by using the properties of BSDEs and related PDEs. The next
lemma provides a simpler derivation for the special case where the terminal
value of the BSDE is a suitably defined indicator function.  }

\begin{lemma}
{\normalsize \label{explicit-solution} For any $a<b\in$$\mathbb{R}$$\ $and
$\kappa>0$,
\[
\mathbb{E}_{[-\kappa,\kappa]}[I_{[a,b]}(B_{1})]=\left\{
\begin{array}
[c]{lcc}%
\Phi_{-\kappa}\left(  -a\right)  -e^{-\kappa(b-a)}\Phi_{-\kappa}\left(
-b\right)  &  & \text{if }a+b\geq0\\
\Phi_{-\kappa}\left(  b\right)  -e^{-\kappa(b-a)}\Phi_{-\kappa}\left(
a\right)  &  & \text{if }a+b<0.
\end{array}
\right.
\]
and
\[
\mathcal{E}_{[-\kappa,\kappa]}[I_{[a,b]}(B_{1})]=\left\{
\begin{array}
[c]{lcc}%
\Phi_{\kappa}\left(  -a\right)  -e^{\kappa(b-a)}\;\Phi_{\kappa}\left(
-b\right)  &  & \text{if }a+b\geq0\\
\Phi_{\kappa}\left(  b\right)  -e^{\kappa(b-a)}\;\Phi_{\kappa}\left(  a\right)
&  & \text{if }a+b<0.
\end{array}
\right.
\]
}
\end{lemma}

{\normalsize \noindent\noindent\textbf{Proof}: For $\kappa>0$, let
{\footnotesize
\begin{align*}
\mathcal{P} \equiv\left\{  Q^{v}:E_{P^{\ast}}[\frac{dQ^{v}}{dP^{\ast}
}|\mathcal{F}_{1}]=e^{-\frac{1}{2}\int_{0}^{1}v_{s} ^{2}ds+\int_{0}^{1}%
v_{s}dB_{s}}, (v_{t})\text{ is }\mathcal{F}_{t}\text{-adapted and } \sup
_{s\in[0,1]}|v_{s}|\leq\kappa\right\}  .
\end{align*}
} }

{\normalsize Let $\varphi=I_{[a,b]}$, then by \cite[Theorem 2.2]{CE} or
\cite[Lemma 3]{chenpeng},
\[
\mathbb{E}_{[-\kappa,\kappa]}[\varphi\left(  B_{1}\right)  ]=\sup
\limits_{Q\in\mathcal{P}}E_{Q}\left[  \varphi( B_{1}) \right]  =\sup
\limits_{|v|\leq\kappa}E_{Q^{v}}\left[  \varphi\left(  B_{1}^{v}+\int_{0}%
^{1}v_{s}ds\right)  \right]  \text{,}%
\]
where $B_{t}^{v}\equiv B_{t}-\int_{0}^{t}v_{s}ds$ is the Brownian motion under
$Q^{v}$.  }

{\normalsize Let $(v_{s})$ be any $\mathcal{F}_{t}$-adapted process valued in
$[-\kappa,\kappa]$, and consider the following BSDEs:%
\begin{align*}
Y_{t}= &  \varphi(B_{1})+\int_{t}^{1}\max\limits_{-\kappa\leq\mu\leq\kappa
}(\mu Z_{s})ds-\int_{t}^{1}Z_{s}dB_{s}\\
Y_{t}^{v}= &  \varphi(\overline{B}_{1}^{v})+\int_{t}^{1}v_{s}Z_{s}^{v}%
ds-\int_{t}^{1}Z_{s}^{v}d\overline{B}_{s}^{v}\\
= &  \varphi(\overline{B}_{1}^{v})-\int_{t}^{1}Z_{s}^{v}dB_{s}\text{
}\\
Y_{t}^{\prime}= &  \varphi(\overline{B}_{1}^{v})+\int_{t}^{1}\max
\limits_{-\kappa\leq\mu\leq\kappa}(\mu Z_{s}^{\prime})ds-\int_{t}^{1}%
Z_{s}^{\prime}d\overline{B}_{s}^{v}\text{,}%
\end{align*}
where $\overline{B}_{t}^{v}\equiv B_{t}+\int_{0}^{t}v_{s}ds$. Clearly,
$Y_{0}=Y_{0}^{\prime}\geq Y_{0}^{v}$, and thus%
\[
\sup\limits_{|v|\leq\kappa}\mathbb{E}_{Q^{v}}[\varphi(B_{1}^{v}+\int_{0}%
^{1}v_{s}ds)]=Y_{0}\geq\sup\limits_{|v|\leq\kappa}Y_{0}^{v}=\sup
\limits_{|v|\leq\kappa}E_{P^{\ast}}[\varphi(B_{1}+\int_{0}^{1}v_{s}ds)]
\]
Let $(X_{t}^{v,x})$ and $(X_{t}^{\ast,x})$ be the solutions respectively of
\begin{align*}
X_{t}^{v,x}  & =x+\int_{0}^{t}v_{s}ds+B_{t},\quad0\leq t\leq1\text{, and}\\
X_{t}^{\ast,x}  & =x-\kappa\int_{0}^{t}sgn\left(  X_{s}^{\ast,x}\right)
ds+B_{t},\quad0\leq t\leq1.
\end{align*}
By the comparison theorem for stochastic differential equations \cite[Thm.
2.1]{watanabe1977},
\begin{align*}
&  \sup\limits_{|v|\leq\kappa}P^{\ast}\left(  a\leq B_{1}+\int_{0}^{1}%
v_{s}ds\leq b\right)  \\
= &  \sup\limits_{|v|\leq\kappa}P^{\ast}\left(  |X_{1}^{v,c}|\leq\frac{b-a}%
{2}\right)  =P^{\ast}\left(  |X_{1}^{\ast,c}|\leq\frac{b-a}{2}\right)
\text{,}%
\end{align*}
where $c=-\frac{a+b}{2}$. }

{\normalsize On the other hand, let $(\alpha_{s})$ be any $\mathcal{F}_{t}%
$-adapted process valued in $[-\kappa,\kappa]$, and $(\overline{X}_{t}%
^{\alpha,x})$ be the solution of
\[
\overline{X}_{t}^{\alpha,x}=x+\int_{0}^{t}\alpha_{s}ds+B_{t}^{v},\quad0\leq
t\leq1
\]
and let $(\overline{X}_{t}^{\ast,x})$ be the solution of%
\[
\overline{X}_{t}^{\ast,x}=x-\kappa\int_{0}^{t}sgn\left(  \overline{X}%
_{s}^{\ast,x}\right)  ds+B_{t}^{v},\quad0\leq t\leq1.
\]
Then
\begin{align*}
&  Q^{v}\left(  a\leq B_{1}^{v}+\int_{0}^{1}\alpha_{s}ds\leq b\right)  \\
= &  Q^{v}\left(  |\overline{X}_{1}^{\alpha,c}|\leq\frac{b-a}{2}\right)  \leq
Q^{v}\left(  |\overline{X}_{1}^{\ast,c}|\leq\frac{b-a}{2}\right)  =P^{\ast
}\left(  |X_{1}^{\ast,c}|\leq\frac{b-a}{2}\right)  \text{,}%
\end{align*}
and
\begin{align*}
&  \sup\limits_{|v|\leq\kappa}Q^{v}\left(  a\leq B_{1}^{v}+\int_{0}^{1}%
v_{s}ds\leq b\right)  \\
\leq &  \sup\limits_{|\alpha|\leq\kappa}Q^{v}\left(  a\leq B_{1}^{v}+\int
_{0}^{1}\alpha_{s}ds\leq b\right)  \leq P^{\ast}\left(  |X_{1}^{\ast,c}%
|\leq\frac{b-a}{2}\right)  .
\end{align*}
That is,
\[
\sup\limits_{|v|\leq\kappa}Q^{v}\left(  a\leq B_{1}^{v}+\int_{0}^{1}%
v_{s}ds\leq b\right)  =P^{\ast}\left(  |X_{1}^{\ast,c}|\leq\frac{b-a}%
{2}\right)
\]
By \cite[Proposition 5.1]{KS1984}, the transition probability density of
$(X_{t}^{\ast,x})$ is given by (for all $t\in(0,1]$, $z\in$$\mathbb{R}$),
\[
q_{x}(t,z)=\frac{1}{\sqrt{2\pi t}}e^{-\frac{(x-z)^{2}+2\kappa
t(|z|-|x|)+\kappa^{2}t^{2}}{2t}}+\kappa e^{-2\kappa|z|}\int_{|x|+|z|-\kappa
t}^{\infty}\frac{1}{\sqrt{2\pi t}}e^{-\frac{u^{2}}{2t}}du\text{.}\quad
\]
Thus
\begin{align*}
P^{\ast}\left(  |X_{1}^{\ast,c}|\leq\frac{b-a}{2}\right)   &  =\int
_{\frac{a-b}{2}}^{\frac{b-a}{2}}q_{c}(1,z)dz\\
&  =\left\{
\begin{array}
[c]{lc}%
\Phi_{-\kappa}\left(  -a\right)  -e^{-\kappa(b-a)}\Phi_{-\kappa}\left(
-b\right)   & \text{if }a+b\geq0\\
\Phi_{-\kappa}\left(  b\right)  -e^{-\kappa(b-a)}\Phi_{-\kappa}\left(
a\right)   & \text{if }a+b<0.
\end{array}
\right.
\end{align*}
}

{\normalsize The rest can be proven in the same way. \hfill\hfill
$\blacksquare$  }

\subsection{{\protect\normalsize Proof of Theorem \ref{thm-CLT}}}

{\normalsize It is enough to prove (\ref{CLT}). We prove it for $\varphi\in
C_{b}^{\infty}(\mathbb{R})$. This suffices because any $\varphi\in
{C([-\infty,\infty])}$ can be approximated uniformly by a sequence of
functions in $C_{b}^{\infty}(\mathbb{R})$ (see Approximation Lemma in
\cite[Ch. VIII]{feller}).  }

{\normalsize Let%
\[
S_{0}\equiv0,\;S_{n}\equiv\sum_{i=1}^{n}X_{i},~S_{n}^{Q}\equiv\sum_{i=1}%
^{n}Y_{i}^{Q}\text{,~~}Y_{i}^{Q}\equiv\frac{1}{\sigma}\left(  X_{i}%
-E_{Q}[X_{i}|\mathcal{G}_{i-1}]\right)  \text{.}%
\]
First we prove that
\begin{equation}
\lim_{\epsilon\rightarrow0}\lim_{n\rightarrow\infty}\left\vert \sup
\limits_{Q\in\mathcal{P}}E_{Q}\left[  \varphi\left(  \frac{S_{n}}{n}%
+\frac{S_{n}^{Q}}{\sqrt{n}}\right)  \right]  -\mathbb{E}_{g_{\epsilon}}\left[
\varphi\left(  B_{1}\right)  \right]  \right\vert =0\text{.} \label{eq3-1}%
\end{equation}
By the definition of $\{H_{m,n}\}_{m=0}^{n}$,
{
\begin{align*}
&  \sup\limits_{Q\in\mathcal{P}}E_{Q}\left[  \varphi\left(  \frac{S_{n}}%
{n}+\frac{S_{n}^{Q}}{\sqrt{n}}\right)  \right]  -\mathbb{E}_{g_{\epsilon}%
}\left[  \varphi\left(  B_{1}\right)  \right] \\
=  &  \sup\limits_{Q\in\mathcal{P}}E_{Q}\left[  H_{n,n}\left(  \frac{S_{n}}%
{n}+\frac{S_{n}^{Q}}{\sqrt{n}}\right)  \right]  -H_{0,n}(0)\\
=  &  \sup\limits_{Q\in\mathcal{P}}E_{Q}\left[  H_{n,n}\left(  \frac{S_{n}}%
{n}+\frac{S_{n}^{Q}}{\sqrt{n}}\right)  \right]  -\sup\limits_{Q\in\mathcal{P}%
}E_{Q}\left[  H_{n-1,n}\left(  \frac{S_{n-1}}{n}+\frac{S_{n-1}^{Q}}{\sqrt{n}%
}\right)  \right] \\
&  +\sup\limits_{Q\in\mathcal{P}}E_{Q}\left[  H_{n-1,n}\left(  \frac{S_{n-1}%
}{n}+\frac{S_{n-1}^{Q}}{\sqrt{n}}\right)  \right]  -\sup\limits_{Q\in
\mathcal{P}}E_{Q}\left[  H_{n-2,n}\left(  \frac{S_{n-2}}{n}+\frac{S_{n-2}^{Q}%
}{\sqrt{n}}\right)  \right] \\
&  +\ldots\\
&  +\sup\limits_{Q\in\mathcal{P}}E_{Q}\left[  H_{m,n}\left(  \frac{S_{m}}%
{n}+\frac{S_{m}^{Q}}{\sqrt{n}}\right)  \right]  -\sup\limits_{Q\in\mathcal{P}%
}E_{Q}\left[  H_{m-1,n}\left(  \frac{S_{m-1}}{n}+\frac{S_{m-1}^{Q}}{\sqrt{n}%
}\right)  \right] \\
&  +\ldots+ \sup\limits_{Q\in\mathcal{P}}E_{Q}\left[  H_{1,n}\left(
\frac{S_{1}}{n}+\frac{S_{1}^{Q}}{\sqrt{n}}\right)  \right]  -H_{0,n}(0)\\
=  &  \sum\limits_{m=1}^{n}\left\{  \sup\limits_{Q\in\mathcal{P}}E_{Q}\left[
H_{m,n}\left(  \frac{S_{m}}{n}+\frac{S_{m}^{Q}}{\sqrt{n}}\right)  \right]
-\sup\limits_{Q\in\mathcal{P}}E_{Q}\bigg[ H_{m-1,n}\bigg( \frac{S_{m-1}}%
{n}+\frac{S_{m-1}^{Q}}{\sqrt{n}}\bigg) \bigg] \right\} \\
=  &  \sum\limits_{m=1}^{n}\left\{  \sup\limits_{Q\in\mathcal{P}}E_{Q}\left[
H_{m,n}\left(  \frac{S_{m}}{n}+\frac{S_{m}^{Q}}{\sqrt{n}}\right)  \right]
-\sup\limits_{Q\in\mathcal{P}}E_{Q}\bigg[ L_{m,n}\bigg( \frac{S_{m-1}}%
{n}+\frac{S_{m-1}^{Q}}{\sqrt{n}}\bigg) \bigg] \right\} \\
&  +\sum_{m=1}^{n}\left\{  \sup\limits_{Q\in\mathcal{P}}E_{Q}\left[
L_{m,n}\left(  \frac{S_{m-1}}{n}+\frac{S_{m-1}^{Q}}{\sqrt{n}}\right)  \right]
\right. \\
&  \left.  \qquad\qquad\qquad\qquad\qquad\qquad\qquad-\sup\limits_{Q\in
\mathcal{P}}E_{Q}\left[  H_{m-1,n}\left(  \frac{S_{m-1}}{n}+\frac{S_{m-1}^{Q}%
}{\sqrt{n}}\right)  \right]  \right\} \\
\equiv &  I_{1n}+I_{2n}\text{,}%
\end{align*}
} where $L_{m,n}(x)=H_{m,n}\left(  x\right)  +\frac{1}{n}g_{0}\left(
H_{m,n}^{\prime}(x)\right)  +\frac{1}{2n}H_{m,n}^{\prime\prime}(x).$  }

{\normalsize By Lemma \ref{lemma-taylor}, if $T_{m,n}=\frac{S_{m}}{n}%
+\frac{S_{m}^{Q}}{\sqrt{n}}$, then {\small
\begin{align*}
|I_{1n}|  &  \leq\sum\limits_{m=1}^{n}\left\vert \sup\limits_{Q\in\mathcal{P}%
}E_{Q}\left[  H_{m,n}\left(  \frac{S_{m}}{n}+\frac{S_{m}^{Q}}{\sqrt{n}%
}\right)  \right]  -\sup\limits_{Q\in\mathcal{P}}E_{Q}\left[  L_{m,n}\left(
\frac{S_{m-1}}{n}+\frac{S_{m-1}^{Q}}{\sqrt{n}}\right)  \right]  \right\vert \\
&  \rightarrow0\text{ ~as }n\rightarrow\infty\text{.}%
\end{align*}
Furthermore, by Lemmas \ref{lemma-ddp} and \ref{lemma-gexpectation}, as
$n\rightarrow\infty$,%
\begin{align*}
|I_{2n}|  &  \leq\sum_{m=1}^{n}\sup\limits_{Q\in\mathcal{P}}E_{Q}\left[
\left\vert L_{m,n}\left(  \frac{S_{m-1}}{n}+\frac{S_{m-1}^{Q}}{\sqrt{n}%
}\right)  -H_{m-1,n}\left(  \frac{S_{m-1}}{n}+\frac{S_{m-1}^{Q}}{\sqrt{n}%
}\right)  \right\vert \right] \\
\leq &  \sum_{m=1}^{n}\sup\limits_{x\in\mathbb{R}}\left\vert L_{m,n}%
(x)-H_{m-1,n}(x)\right\vert \\
=  &  \sum_{m=1}^{n}\sup\limits_{x\in\mathbb{R}}\left\vert H_{m,n}\left(
x\right)  +\frac{1}{n}g_{0}\left(  H_{m,n}^{\prime}(x)\right)  +\frac{1}%
{2n}H_{m,n}^{\prime\prime}(x) \right. \\
&  \qquad\qquad\qquad\qquad\qquad\qquad\qquad\qquad\quad- \mathbb{E}%
_{g_{\epsilon}}\left[  H_{m,n}\left(  x+B_{\frac{m}{n}}-B_{\frac{m-1}{n}%
}\right)  \right]  \bigg\vert\\
\leq &  \sum_{m=1}^{n}\sup\limits_{x\in\mathbb{R}}\left\vert H_{m,n}\left(
x\right)  +\frac{1}{n}g_{\epsilon}\left(  H_{m,n}^{\prime}(x)\right)
+\frac{1}{2n}H_{m,n}^{\prime\prime}(x) - \mathbb{E}_{g_{\epsilon}}\left[
H_{m,n}\left(  x+B_{\frac{1}{n}}\right)  \right]  \right\vert \\
&  +\frac{1}{n}\sum_{m=1}^{n}\sup\limits_{x\in\mathbb{R}}\left\vert
g_{\epsilon}\left(  H_{m,n}^{\prime}(x)\right)  -g_{0}\left(  H_{m,n}^{\prime
}(x)\right)  \right\vert \\
\leq &  \sum_{m=1}^{n}\sup\limits_{x\in\mathbb{R}}\left\vert \mathbb{E}%
_{g_{\epsilon}}\left[  H_{m,n}\left(  x+B_{\frac{1}{n}}\right)  \right]
-H_{m,n}\left(  x\right)  -\frac{1}{n}g_{\epsilon}\left(  H_{m,n}^{\prime
}(x)\right)  -\frac{1}{2n}H_{m,n}^{\prime\prime}(x)\right\vert \\
&  +2\kappa\epsilon\text{,}%
\end{align*}
which sum converges to $2\kappa\epsilon$. This proves (\ref{eq3-1}). }  }

{\normalsize From the standard estimates for BSDEs \cite[Proposition 2.1]%
{KPQ},
\begin{align*}
\left\vert \mathbb{E}_{g_{\epsilon}}\left[  \varphi\left(  B_{1}\right)
\right]  -\mathbb{E}_{g_{0}}\left[  \varphi\left(  B_{1}\right)  \right]
\right\vert ^{2}  &  \leq{\widehat{C}}E_{P^{\ast}}\left[  \left(  \int_{0}%
^{1}\left\vert g_{\epsilon}(Z_{s}^{\epsilon})-g_{0}(Z_{s}^{\epsilon
})\right\vert ds\right)  ^{2}\right] \\
&  <{\widehat{C}4\kappa^{2}}\epsilon^{2}\text{,}%
\end{align*}
where {$\widehat{C}>0$} is a constant. Combine with (\ref{eq3-1}) to obtain%
\begin{align*}
&  \lim_{n\rightarrow\infty}\left\vert \sup\limits_{Q\in\mathcal{P}}%
E_{Q}\left[  \varphi\left(  \frac{S_{n}}{n}+\frac{S_{n}^{Q}}{n}\right)
\right]  -\mathbb{E}_{g_{0}}\left[  \varphi\left(  B_{1}\right)  \right]
\right\vert \\
\leq &  \lim_{\epsilon\rightarrow0}\lim_{n\rightarrow\infty}\left\vert
\sup\limits_{Q\in\mathcal{P}}E_{Q}\left[  \varphi\left(  \frac{S_{n}}{n}%
+\frac{S_{n}^{Q}}{n}\right)  \right]  -\mathbb{E}_{g_{\epsilon}}\left[
\varphi\left(  B_{1}\right)  \right]  \right\vert \\
&  +\lim_{\epsilon\rightarrow0}\left\vert \mathbb{E}_{g_{\epsilon}}\left[
\varphi\left(  B_{1}\right)  \right]  -\mathbb{E}_{g_{0}}\left[
\varphi\left(  B_{1}\right)  \right]  \right\vert \text{.}%
\end{align*}
The latter sum equals $0$, thus completing the proof under condition
(\ref{kappa}).  }

{\normalsize Finally, we describe the proof for general \underline{$\mu$} and
$\overline{\mu}$. Let $Y_{i}=X_{i}-\tfrac{\overline{\mu}+\underline{\mu}}{2}$
and $\kappa=\frac{\overline{\mu}-\underline{\mu}}{2}$. Then
\[
{\mathbb{E}}[Y_{i}\mid\mathcal{G}_{i-1}]=\frac{\overline{\mu}-\underline{\mu}%
}{2}=\kappa,\quad{\mathcal{E}}[Y_{i}\mid\mathcal{G}_{i-1}]=-\frac
{\overline{\mu}-\underline{\mu}}{2}=-\kappa.
\]
Apply the above result to $(Y_{i})$ and $\widehat{\varphi}$, $\widehat
{\varphi}(x)=\varphi\left(  x+\tfrac{\overline{\mu}+\underline{\mu}}%
{2}\right)  $, to obtain%
\begin{align*}
&  \lim_{n\rightarrow\infty}\sup\limits_{Q\in\mathcal{P}}E_{Q}\left[
\varphi\left(  \frac{S_{n}}{n}+\frac{S_{n}^{Q}}{n}\right)  \right] \\
=  &  \lim_{n\rightarrow\infty}\sup_{Q\in\mathcal{P}}E_{Q}\left[
\varphi\left(  \frac{1}{n}\sum\limits_{i=1}^{n}X_{i}+\frac{1}{\sqrt{n}}%
\sum\limits_{i=1}^{n}\frac{1}{\sigma}\left(  X_{i}-E_{Q}[X_{i}|\mathcal{G}%
_{i-1}]\right)  \right)  \right] \\
=  &  \lim_{n\rightarrow\infty}\sup_{Q\in\mathcal{P}}E_{Q}\left[
\varphi\left(  \frac{1}{n}\sum\limits_{i=1}^{n}Y_{i}+\frac{\overline{\mu
}+\underline{\mu}}{2}+\frac{1}{\sqrt{n}}\sum\limits_{i=1}^{n}\frac{1}{\sigma
}\left(  Y_{i}-E_{Q}[Y_{i}|\mathcal{G}_{i-1}]\right)  \right)  \right] \\
=  &  \lim_{n\rightarrow\infty}\sup_{Q\in\mathcal{P}}E_{Q}\left[
\widehat{\varphi}\left(  \frac{1}{n}\sum\limits_{i=1}^{n}Y_{i}+\frac{1}%
{\sqrt{n}}\sum\limits_{i=1}^{n}\frac{1}{\sigma}\left(  Y_{i}-E_{Q}%
[Y_{i}|\mathcal{G}_{i-1}]\right)  \right)  \right] \\
=  &  \mathbb{E}_{\left[  \tfrac{\underline{\mu}-\overline{\mu}}{2}%
,\tfrac{\overline{\mu}-\underline{\mu}}{2}\right]  }\left[  \widehat{\varphi
}\left(  B_{1}\right)  \right] \\
=  &  \mathbb{E}_{\left[  \tfrac{\underline{\mu}-\overline{\mu}}{2}%
,\tfrac{\overline{\mu}-\underline{\mu}}{2}\right]  }\left[  \varphi\left(
\frac{\overline{\mu}+\underline{\mu}}{2}+B_{1}\right)  \right]  =\mathbb{E}%
_{\left[  \underline{\mu},\overline{\mu}\right]  }\left[  \varphi\left(
B_{1}\right)  \right]  ,
\end{align*}
where the last equality is due to Lemma \ref{lemma-shift}. This completes the
proof.\hfill\hfill\hfill\hfill\ $\blacksquare$  }

\begin{remark}
{\normalsize \label{remark-CLT2b}Straightforward modifications of the
preceding arguments deliver a proof of Theorem \ref{thm-CLT2}. The key is
modification of Lemma \ref{lemma-taylor} (see Remark \ref{remark-CLT2} and
Appendix \ref{app-CLT2}). The remaining arguments are similar to those given
above and are omitted.  }
\end{remark}

\subsection{{\protect\normalsize Proof of Theorem \ref{thm-special}}}

{\normalsize \noindent\textbf{Proof of (1):}\textbf{ } Let {$\varphi\in
C([-\infty,\infty])$} be symmetric with {center $c\in\mathbb{R}$} and
decreasing on $(c,\infty)$. {The result is clear if $\varphi$ is globally
constant. Thus we assume that $\varphi$ is not a globally constant function.
}Then $\varphi$ can be approximated uniformly by $\varphi_{h}$ defined by%
\begin{equation}
\varphi_{h}(x)=\int_{-\infty}^{\infty}\frac{1}{\sqrt{2\pi}}\varphi
(x+hy)e^{-\tfrac{y^{2}}{2}}dy\text{,} \label{phih}%
\end{equation}
and (see Appendix \ref{app-special}), $\varphi_{h}$ is symmetric with center
$c$, and satisfies
\[
sgn(\varphi_{h}^{\prime}(x))=-sgn(x-c)\text{, ~}\forall h>0\text{.}%
\]
}

{\normalsize
}

{\normalsize Consider the special case (\ref{kappa}). Let $\{H_{m,n}%
\}_{m=0}^{n}$ be defined via (\ref{BSDEep-1}) using $\varphi$, where, without
loss of generality we assume $\varphi\in C_{b}^{3}(\mathbb{R})$; otherwise, we
can use $\varphi_{h}$ defined in (\ref{phih}).  }

{\normalsize Let
\[
S_{0}\equiv0,\;S_{n}\equiv\sum_{i=1}^{n}X_{i},~\overline{S}_{n}\equiv
\sum_{i=1}^{n}Z_{i}^{n}\text{,~~}Z_{i}^{n}\equiv\frac{1}{\sigma}\left(
X_{i}-\mu_{i}^{n}\right)  \text{.}%
\]
First, prove that
\begin{equation}
\lim_{\epsilon\rightarrow0}\lim_{n\rightarrow\infty}\left\vert \sup
\limits_{Q\in\mathcal{P}}E_{Q}\left[  \varphi\left(  \frac{S_{n}}{n}%
+\frac{\overline{S}_{n}}{\sqrt{n}}\right)  \right]  -\mathbb{E}_{g_{\epsilon}%
}\left[  \varphi\left(  B_{1}\right)  \right]  \right\vert =0\text{.}
\label{limit}%
\end{equation}
Argue as follows:
\begin{align*}
&  \sup\limits_{Q\in\mathcal{P}}E_{Q}\left[  {\varphi}\left(  \frac{S_{n}}%
{n}+\frac{\overline{S}_{n}}{\sqrt{n}}\right)  \right]  -\mathbb{E}%
_{g_{\epsilon}}\left[  {\varphi}\left(  B_{1}\right)  \right] \\
=  &  \sup\limits_{Q\in\mathcal{P}}E_{Q}\left[  H_{n,n}\left(  \frac{S_{n}}%
{n}+\frac{\overline{S}_{n}}{\sqrt{n}}\right)  \right]  -H_{0,n}(0)\\
=  &  \sup\limits_{Q\in\mathcal{P}}E_{Q}\left[  H_{n,n}\left(  \frac{S_{n}}%
{n}+\frac{\overline{S}_{n}}{\sqrt{n}}\right)  \right]  -\sup\limits_{Q\in
\mathcal{P}}E_{Q}\left[  H_{n-1,n}\left(  \frac{S_{n-1}}{n}+\frac{\overline
{S}_{n-1}}{\sqrt{n}}\right)  \right] \\
&  +\;\sup\limits_{Q\in\mathcal{P}}E_{Q}\left[  H_{n-1,n}\left(  \frac
{S_{n-1}}{n}+\frac{\overline{S}_{n-1}}{\sqrt{n}}\right)  \right]
-\sup\limits_{Q\in\mathcal{P}}E_{Q}\left[  H_{n-2,n}\left(  \frac{S_{n-2}}%
{n}+\frac{\overline{S}_{n-2}}{\sqrt{n}}\right)  \right] \\
&  +\ldots\\
&  +\;\sup\limits_{Q\in\mathcal{P}}E_{Q}\left[  H_{m,n}\left(  \frac{S_{m}}%
{n}+\frac{\overline{S}_{m}}{\sqrt{n}}\right)  \right]  -\sup\limits_{Q\in
\mathcal{P}}E_{Q}\left[  H_{m-1,n}\left(  \frac{S_{m-1}}{n}+\frac{\overline
{S}_{m-1}}{\sqrt{n}}\right)  \right] \\
&  +\ldots\\
&  +\;\sup\limits_{Q\in\mathcal{P}}E_{Q}\left[  H_{1,n}\left(  \frac{S_{1}}%
{n}+\frac{\overline{S}_{1}}{\sqrt{n}}\right)  \right]  -H_{0,n}(0)\\
=  &  \sum\limits_{m=1}^{n}\left\{  \sup\limits_{Q\in\mathcal{P}}E_{Q}\left[
H_{m,n}\left(  \frac{S_{m}}{n}+\frac{\overline{S}_{m}}{\sqrt{n}}\right)
\right]  -\sup\limits_{Q\in\mathcal{P}}E_{Q}\left[  L_{m,n}\left(
\frac{S_{m-1}}{n}+\frac{\overline{S}_{m-1}}{\sqrt{n}}\right)  \right]
\right\} \\
&  +\sum_{m=1}^{n}\left\{  \sup\limits_{Q\in\mathcal{P}}E_{Q}\left[
L_{m,n}\left(  \frac{S_{m-1}}{n}+\frac{\overline{S}_{m-1}}{\sqrt{n}}\right)
\right]  \right. \\
&  \left.  \qquad\qquad\qquad\qquad\qquad\qquad\qquad-\sup\limits_{Q\in
\mathcal{P}}E_{Q}\left[  H_{m-1,n}\left(  \frac{S_{m-1}}{n}+\frac{\overline
{S}_{m-1}}{\sqrt{n}}\right)  \right]  \right\} \\
\equiv &  J_{1n}+J_{2n},
\end{align*}
where $L_{m,n}(x)=H_{m,n}\left(  x\right)  +\frac{1}{n}g_{0}\left(
H_{m,n}^{\prime}(x)\right)  +\frac{1}{2n}H_{m,n}^{\prime\prime}(x).$  }

{\normalsize By Lemma \ref{lemma-taylor2}(1), with $T_{m,n}=\frac{S_{m}}%
{n}+\frac{\overline{S}_{m}}{\sqrt{n}}$, we have{\small
\begin{align*}
|J_{1n}|  &  \leq\sum\limits_{m=1}^{n}\left\vert \sup\limits_{Q\in\mathcal{P}%
}E_{Q}\left[  H_{m,n}\left(  \frac{S_{m}}{n}+\frac{\overline{S}_{m}}{\sqrt{n}%
}\right)  \right]  -\sup\limits_{Q\in\mathcal{P}}E_{Q}\left[  L_{m,n}\left(
\frac{S_{m-1}}{n}+\frac{\overline{S}_{m-1}}{\sqrt{n}}\right)  \right]
\right\vert \\
&  \rightarrow0\text{ ~as }n\rightarrow\infty\text{.}%
\end{align*}
}{As in the proof of Theorem \ref{thm-CLT}, we have $|J_{2n}|\rightarrow0,~$as
$n\rightarrow\infty$. Hence, (\ref{limit}) holds. Combine it with standard
estimate for BSDEs to complete the proof under condition (\ref{kappa}).}  }

{\normalsize For the case of general $\underline{\mu}$ and $\overline{\mu}$,
let $Y_{i}=X_{i}-\tfrac{\overline{\mu}+\underline{\mu}}{2}$. Then
\[
{\mathbb{E}[Y_{i}\mid\mathcal{G}_{i-1}]=\frac{\overline{\mu}-\underline{\mu}%
}{2},\quad\mathcal{E}[Y_{i}\mid\mathcal{G}_{i-1}]=-\frac{\overline{\mu
}-\underline{\mu}}{2}.}%
\]
Apply the above result for $(Y_{i})$ to $\widehat{\varphi}$, $\widehat
{\varphi}(x)=\varphi\left(  x+\tfrac{\overline{\mu}+\underline{\mu}}%
{2}\right)  $, to obtain%
\begin{align*}
&  \lim_{n\rightarrow\infty}\sup_{Q\in\mathcal{P}}E_{Q}\left[  \varphi\left(
\frac{1}{n}\sum\limits_{i=1}^{n}X_{i}+\frac{1}{\sqrt{n}}\sum\limits_{i=1}%
^{n}\frac{1}{\sigma}\left(  X_{i}-\mu_{i}^{n}\right)  \right)  \right]  \\
= &  \lim_{n\rightarrow\infty}\sup_{Q\in\mathcal{P}}E_{Q}\left[
\varphi\left(  \frac{1}{n}\sum\limits_{i=1}^{n}Y_{i}+\frac{\overline{\mu
}+\underline{\mu}}{2}+\frac{1}{\sqrt{n}}\sum\limits_{i=1}^{n}\frac{1}{\sigma
}\left(  Y_{i}-\left(  \mu_{i}^{n}-\frac{\overline{\mu}+\underline{\mu}}%
{2}\right)  \right)  \right)  \right]  \\
= &  \lim_{n\rightarrow\infty}\sup_{Q\in\mathcal{P}}E_{Q}\left[
\widehat{\varphi}\left(  \frac{1}{n}\sum\limits_{i=1}^{n}Y_{i}+\frac{1}%
{\sqrt{n}}\sum\limits_{i=1}^{n}\frac{1}{\sigma}\left(  Y_{i}-\gamma_{i}%
^{n}\right)  \right)  \right]  \\
= &  \mathbb{E}_{\left[  \tfrac{\underline{\mu}-\overline{\mu}}{2}%
,\tfrac{\overline{\mu}-\underline{\mu}}{2}\right]  }\left[  \widehat{\varphi
}\left(  B_{1}\right)  \right]  \\
= &  \mathbb{E}_{\left[  \tfrac{\underline{\mu}-\overline{\mu}}{2}%
,\tfrac{\overline{\mu}-\underline{\mu}}{2}\right]  }\left[  \varphi\left(
\frac{\overline{\mu}+\underline{\mu}}{2}+B_{1}\right)  \right]  \\
= &  \mathbb{E}_{\left[  \underline{\mu},\overline{\mu}\right]  }\left[
\varphi\left(  B_{1}\right)  \right]  ,
\end{align*}
where the last equality is due to Lemma \ref{lemma-shift}. Also,
\[
\gamma_{m}^{n}=\tfrac{\overline{\mu}-\underline{\mu}}{2}I_{\widehat{A}%
_{m-1,n}}+\tfrac{\underline{\mu}-\overline{\mu}}{2}I_{\widehat{A}_{m-1,n}^{c}%
},\text{ ~and}%
\]
\begin{align*}
\widehat{A}_{m-1,n} &  =\left\{  \tfrac{1}{n}\sum\limits_{i=1}^{m-1}%
Y_{i}+\tfrac{1}{\sqrt{n}}\sum\limits_{i=1}^{m-1}\tfrac{1}{\sigma}\left(
Y_{i}-\gamma_{i}^{n}\right)  \leq-\tfrac{\overline{\mu}+\underline{\mu}}%
{2}+c\right\}  \\
&  =\bigg\{\tfrac{1}{n}\sum\limits_{i=1}^{m-1}X_{i}+\tfrac{1}{\sqrt{n}}%
\sum\limits_{i=1}^{m-1}\tfrac{1}{\sigma}\big(X_{i}-\gamma_{i}^{n}%
-\tfrac{\overline{\mu}+\underline{\mu}}{2}\big)\leq-\tfrac{\overline{\mu
}+\underline{\mu}}{2}\big(1-\tfrac{m-1}{n}\big)+c\bigg\}\text{.}%
\end{align*}
Thus $\widehat{A}_{0,n}=A_{0,n}$, and
\[
\gamma_{1}^{n}+\tfrac{\overline{\mu}+\underline{\mu}}{2}=\overline{\mu
}I_{\widehat{A}_{0,n}}+\underline{\mu}I_{\widehat{A}_{0,n}^{c}}=\overline{\mu
}I_{A_{0,n}}+\underline{\mu}I_{A_{0,n}^{c}}=\mu_{1}^{n}\text{.}%
\]
By induction, $A_{m-1,n}=\widehat{A}_{m-1,n}$, for $m\geq1$, and%
\[
\gamma_{m}^{n}+\tfrac{\overline{\mu}+\underline{\mu}}{2}=\overline{\mu
}I_{\widehat{A}_{m-1,n}}+\underline{\mu}I_{\widehat{A}_{m-1,n}^{c}}%
=\overline{\mu}I_{A_{m-1,n}}+\underline{\mu}I_{A_{m-1,n}^{c}}=\mu_{m}%
^{n}\text{.}%
\]
This completes the proof of (\ref{sym-clt1}). }

{\normalsize By standard limiting arguments, (\ref{sym-clt1}) can be extended
to indicator functions for intervals. Then (\ref{sym-indi-clt}) follows from
Lemmas \ref{lemma-shift} and \ref{explicit-solution}.  }

{\normalsize \medskip}

{\normalsize \noindent\textbf{Proof of (2):} In light of Lemma
\ref{lemma-taylor2}(2), the proof of part (2) is similar to that of (1) and is
omitted.\hfill\hfill\hfill\hfill\ $\blacksquare$  }

{\normalsize
}

{\normalsize \newpage}

{\normalsize \appendix
}

\section{{\protect\normalsize Supplementary Appendix}}

\subsection{{\protect\normalsize Rectangularity\label{app-rect}}}

{\normalsize Let $\mathcal{P}\subset\Delta\left(  \Omega,\mathcal{G}\right)  $
be rectangular. All measures in $\mathcal{P}$ are equivalent on each
$\mathcal{G}_{n}$ and relations between $\mathcal{G}_{n}$-measurable r.v.s
should be understood to hold $P_{0}$-a.s. for some fixed measure $P_{0}$ in
$\mathcal{P}$. $\mathcal{H}$ denotes the set of r.v.s $X$ on $\left(
\Omega,\mathcal{G}\right)  $ satisfying $\sup_{Q\in\mathcal{P}}E_{Q}%
[|X|]<\infty$.  }

\begin{lemma}
{\normalsize \label{Lem_increasing_seq} For any $X\in\mathcal{H}$ and any $n$,
there is a sequence $E_{P_{i}}[X|\mathcal{G}_{n}]$ in $\{E_{P}[X|\mathcal{G}%
_{n}]:P\in\mathcal{P}\}$ such that $ess\sup_{P\in\mathcal{P}}E_{P}%
[X|\mathcal{G}_{n}]$ is the increasing limit of $E_{P_{i}}[X|\mathcal{G}_{n}%
]$.  }
\end{lemma}

{\normalsize \noindent\textbf{Proof:} We prove that $\{E_{P}[X|\mathcal{G}%
_{n}]:P\in\mathcal{P}\}$ is an upward-directed set. Then the result follows
from \cite[Theorem A.32]{FS}.  }

{\normalsize \noindent Let $Q_{1},Q_{2}\in\mathcal{P}$, and $\forall$%
$B\in\mathcal{G}_{n}$, $\forall\omega=(\omega^{(n)},\omega_{(n+1)})\in\Omega$,
$\forall D\in\mathcal{G}_{(n+1)}$, define
\[
\lambda(\omega^{(n)},D)=\left\{
\begin{array}
[c]{ll}%
Q_{1}(\prod_{1}^{n}\Omega_{i}\times D|\mathcal{G}_{n})(\omega^{(n)}%
)~\ \text{if}~\omega\in B, & \\
Q_{2}(\prod_{1}^{n}\Omega_{i}\times D|\mathcal{G}_{n})(\omega^{(n)}%
)~\ \text{if}~\omega\in B^{c}, &
\end{array}
\right.  \text{ \ \ and }p_{n}=Q_{1}|_{\mathcal{G}_{n}}.
\]
Then, $p_{n}\in\mathcal{P}_{0,n}$ and $\lambda$ is a $\mathcal{P}$-kernel. By
rectangularity, $P\in\mathcal{P}$, where, for $A\in\mathcal{G}$,%
\begin{align*}
&  P(A)=\int_{\prod_{1}^{n}\Omega_{i}}\int_{\prod_{n+1}^{\infty}\Omega_{i}%
}I_{A}(\omega^{(n)},\omega_{(n+1)})\lambda(\omega^{(n)},d\omega_{(n+1)}%
)p_{n}(d\omega^{(n)})\\
&  =\int_{\prod_{1}^{n}\Omega_{i}}\int_{\prod_{n+1}^{\infty}\Omega_{i}%
}I_{A\cap B}(\omega^{(n)},\omega_{(n+1)})Q_{1}(\prod_{i=1}^{n}\Omega_{i}\times
d\omega_{(n+1)}|\mathcal{G}_{n})(\omega^{(n)})p_{n}(d\omega^{(n)})\\
&  +\int_{\prod_{1}^{n}\Omega_{i}}\int_{\prod_{n+1}^{\infty}\Omega_{i}%
}I_{A\cap B^{c}}(\omega^{(n)},\omega_{(n+1)})Q_{2}(\prod_{i=1}^{n}\Omega
_{i}\times d\omega_{(n+1)}|\mathcal{G}_{n})(\omega^{(n)})p_{n}(d\omega
^{(n)})\text{.}%
\end{align*}
}

{\normalsize Consider the probability measure $\widetilde{P}$ with
Radon-Nikodym derivative%
\[
\frac{d\widetilde{P}}{dP_{0}}=\frac{dQ_{1}}{dP_{0}}I_{B}+\frac{(\frac{dQ_{1}%
}{dP_{0}})_{n}}{(\frac{dQ_{2}}{dP_{0}})_{n}}\frac{dQ_{2}}{dP_{0}}I_{B^{c}},
\]
where $(\frac{dQ_{1}}{dP_{0}})_{n}$ means $E_{P_{0}}[\frac{dQ_{1}}{dP_{0}%
}|_{\mathcal{G}_{n}}]$. We claim that $P(A)=\widetilde{P}(A)$, $\forall
A\in\mathcal{G}$. Indeed, by the definitions, for all $n$, $P(A)=\widetilde
{P}(A)$, $\forall A\in\mathcal{C}_{n}$, where
\[
\mathcal{C}_{n}=\{A^{(n)}\times A_{(n+1)}:A^{(n)}\in\mathcal{G}_{n}%
,~A_{(n+1)}\in\mathcal{G}_{(n+1)}\}\text{.}%
\]
Since $\mathcal{C}$$_{n}$ is a $\pi$ class, and satisfies $\sigma
(\mathcal{C}_{n})=\mathcal{G}$, $\ P$ and $\widetilde{P}$ are identical on
$\mathcal{G}$.  }

{\normalsize Note that $(\tfrac{d\widetilde{P}}{dP_{0}})_{n}=(\tfrac{dQ_{1}%
}{dP_{0}})_{n},$ and, by Bayes rule,  }

{\normalsize
\begin{equation}%
\begin{split}
E_{P}[X|\mathcal{G}_{n}]  &  =E_{\widetilde{P}}[X|\mathcal{G}_{n}]\\
&  =E_{P_{0}}\left[  X\frac{d\widetilde{P}}{dP_{0}}|\mathcal{G}_{n}\right]
\left(  (\tfrac{dQ_{1}}{dP_{0}})_{n}\right)  ^{-1}\\
&  =E_{P_{0}}\left[  X\left(  \frac{dQ_{1}}{dP_{0}} I_{B}+\frac{(\frac{dQ_{1}%
}{dP_{0}})_{n}}{(\frac{dQ_{2}}{dP_{0}})_{n}}\frac{dQ_{2}}{dP_{0}} I_{B^{c}%
}\right)  |\mathcal{G}_{n}\right]  [(\tfrac{dQ_{1}}{dP_{0}})_{n}]^{-1}\\
&  = I_{B}E_{Q_{1}}\left[  X|\mathcal{G}_{n}\right]  + I_{B^{c}}E_{Q_{2}%
}\left[  X|\mathcal{G}_{n}\right]  \text{.}%
\end{split}
\label{stable by bifurcation}%
\end{equation}
}

{\normalsize \noindent If $B=\{\omega\in\Omega:E_{Q_{1}}[X|\mathcal{G}%
_{n}](\omega)>E_{Q_{2}}[X|\mathcal{G}_{n}](\omega)\}$, then
\[
E_{P}[X|\mathcal{G}_{n}]=ess\sup\{E_{Q_{1}}\left[  X|\mathcal{G}_{n}\right]
,E_{Q_{2}}\left[  X|\mathcal{G}_{n}\right]  \}\text{. \ \ \ }\blacksquare
\]
\hfill\hfill}

{\normalsize \medskip\bigskip}

{\normalsize \noindent\noindent}

{\normalsize \noindent\textbf{Proof of Lemma \ref{lemma-rect}}: (i) For any
$\omega^{(n)} \in\prod_{1}^{n}\Omega_{i}$ and $B\in\mathcal{G}_{(n+1)}$,
define
\[
\lambda(\omega^{(n)},B)=Q(\prod_{i=1}^{n}\Omega_{i}\times B|\mathcal{G}%
_{n})(\omega^{(n)})\text{ and }p_{n}=R|_{\mathcal{G}_{n}}\text{.}%
\]
Then $p_{n}\in\mathcal{P}_{0,n}$ and $\lambda$ is a $\mathcal{P}$-kernel. By
rectangularity, $P\in\mathcal{P}$, where, for $A\in\mathcal{G}$,
\[
P(A)=\int_{\prod_{1}^{n}\Omega_{i}}\int_{\prod_{n+1}^{\infty}\Omega_{i}}
I_{A}(\omega^{(n)},\omega_{(n+1)})Q(\prod_{i=1}^{n}\Omega_{i}\times
d\omega_{(n+1)}|\mathcal{G}_{n})(\omega^{(n)})p_{n}(d\omega^{(n)})\text{.}%
\]
Consider the probability measure $\widetilde{P}$ with Radon-Nikodym derivative%
\[
\frac{d\widetilde{P}}{dP_{0}}=\frac{(\frac{dR}{dP_{0}})_{n}}{(\frac{dQ}%
{dP_{0}})_{n}}\frac{dQ}{dP_{0}},
\]
where $(\frac{dR}{dP_{0}})_{n}$ means $E_{P_{0}}[\frac{dR}{dP_{0}%
}|_{\mathcal{G}_{n}}]$. Argue as in the proof of Lemma
\ref{Lem_increasing_seq}, to show that $P$ and $\widetilde{P}$ are identical
on $\mathcal{G}$. Note that $(\tfrac{d\widetilde{P}}{dP_{0}})_{n}=(\tfrac
{dR}{dP_{0}})_{n},$ and, by Bayes rule, for any $m<n$ and $X\in\mathcal{H}$,
\begin{align*}
E_{P}[X|\mathcal{G}_{m}]=E_{\widetilde{P}}[X|\mathcal{G}_{m}]  &
=E_{\widetilde{P}}[E_{\widetilde{P}}[X|\mathcal{G}_{n}]|\mathcal{G}_{m}]\\
&  =E_{\widetilde{P}}\left[  E_{P_{0}}\left[  X\frac{d\widetilde{P}}{dP_{0}%
}|\mathcal{G}_{n}\right]  \left(  (\tfrac{dR}{dP_{0}})_{n}\right)
^{-1}|\mathcal{G}_{m}\right] \\
&  =E_{\widetilde{P}}\left[  E_{P_{0}}\left[  X\frac{(\frac{dR}{dP_{0}})_{n}%
}{(\frac{dQ}{dP_{0}})_{n}}\frac{dQ}{dP_{0}}|\mathcal{G}_{n}\right]  \left(
(\tfrac{dR}{dP_{0}})_{n}\right)  ^{-1}|\mathcal{G}_{m}\right] \\
&  =E_{\widetilde{P}}\left[  E_{Q}\left[  X|\mathcal{G}_{n}\right]
|\mathcal{G}_{m}\right]  =E_{R}\left[  E_{Q}\left[  X|\mathcal{G}_{n}\right]
|\mathcal{G}_{m}\right]  \text{.}~
\end{align*}
}

{\normalsize \medskip}

{\normalsize \noindent(ii) can be proven using (\ref{stable by bifurcation}).
}

{\normalsize \smallskip}

{\normalsize \noindent(iii) By Lemma \ref{Lem_increasing_seq}, there exist
increasing sequences $\{E_{Q_{i}}[\phi(X)|\mathcal{G}_{n}]\}$ and $\{E_{P_{j}%
}[\mathbb{E}[X|\mathcal{G}_{n}]|\mathcal{G}_{m}]|\}$, with $Q_{i},P_{j}%
\in\mathcal{P}$ for all $i$ and $j$, such that
\begin{align*}
\mathbb{E}[X|\mathcal{G}_{n}]  &  =\lim_{i\rightarrow\infty}E_{Q_{i}%
}[X|\mathcal{G}_{n}],\text{ and}\\
\mathbb{E}[\mathbb{E}[X|\mathcal{G}_{n}]|\mathcal{G}_{m}]  &  =\lim
_{j\rightarrow\infty}E_{P_{j}}[\mathbb{E}[X|\mathcal{G}_{n}]|\mathcal{G}%
_{m}]\text{.}%
\end{align*}
By the monotone convergence theorem and (i),
\begin{align*}
\mathbb{E}[\mathbb{E}[X|\mathcal{G}_{n}]|\mathcal{G}_{m}]  &  =\lim
_{j\rightarrow\infty}\lim_{i\rightarrow\infty}E_{P_{j}}[E_{Q_{i}%
}[X|\mathcal{G}_{n}]|\mathcal{G}_{m}]\\
&  \leq ess\sup_{P\in\mathcal{P}}E_{P}[X|\mathcal{G}_{m}]\\
&  =\mathbb{E}[X|\mathcal{G}_{m}].
\end{align*}
For the reverse inequality, we have%
\begin{align*}
\mathbb{E}[X|\mathcal{G}_{m}]  &  =ess\sup_{P\in\mathcal{P}}E_{P}%
[E_{P}[X|\mathcal{G}_{n}]|\mathcal{G}_{m}]\\
&  \leq ess\sup_{P\in\mathcal{P}}E_{P}[ess\sup_{Q\in\mathcal{P}}%
E_{Q}[X|\mathcal{G}_{n}]|\mathcal{G}_{m}]\\
&  =\mathbb{E}[\mathbb{E}[X|\mathcal{G}_{n}]|\mathcal{G}_{m}]\text{. \ \ \ \ }%
\end{align*}
\hfill\hfill}

{\normalsize \noindent(iv) can be proven using (iii). \hfill\hfill
$\blacksquare$  }

\subsection{{\protect\normalsize IID model: Lemma \ref{lemma-IID}%
\label{app-IID}}}

{\normalsize Part (i) was proven in the text. (iii) follows from (i) and Lemma
\ref{lemma-rect}.  }

{\normalsize \noindent For (iv), use (i) and (iii) to argue that, for
example,
\begin{align*}
&  \sup_{Q\in\mathcal{P}^{IID}}E_{Q}\left[  (X_{n}-E_{Q}[X_{n}|\mathcal{G}%
_{n-1}])^{2}|\mathcal{G}_{n-1}\right] \\
&  =\sup_{Q\in\mathcal{P}^{IID}}\left\{  E_{Q}\left[  X_{n}^{2}|\mathcal{G}%
_{n-1}\right]  -(E_{Q}[X_{n}|\mathcal{G}_{n-1}])^{2}\right\} \\
&  =\sup_{q\in\mathcal{L}}\left\{  E_{q}[X_{n}^{2}]-(E_{q}[X_{n}%
])^{2}\right\}  =\sup_{q\in\mathcal{L}}E_{q}\left[  (\overline{X}%
-E_{q}\overline{X})^{2}\right]  \text{. }\hfill\
\end{align*}
}

{\normalsize The equivalence on each $\mathcal{G}_{n}$ stated in (ii) is
proven by induction. Let $P,Q\in$$\mathcal{P}^{IID}$. Equivalence on
$\mathcal{G}_{1}$ is due to the equivalence of measures in $\mathcal{L}$.
Suppose $P$ and $Q$ are equivalent on $\mathcal{G}_{n-1}$, and prove
equivalence on $\mathcal{G}_{n}$. \newline Let $A\in$$\mathcal{G}_{n}$,
$P(A)=0$. Then{%
\[
E_{P}\left[  E_{P}\left[  I_{A}|\mathcal{G}_{n-1}\right]  \right]
=0\Leftrightarrow\ P\left(  \{E_{P}[I_{A}|\mathcal{G}_{n-1}]>0\}\right)  =0.
\]
} By the equivalence of measures in $\mathcal{L}$, $\forall$$\omega^{(n-1)}%
\in\prod_{1}^{n-1}\Omega_{i}$,
\[
\left\{  \omega^{(n-1)}:\ E_{P}[I_{A}|\mathcal{G}_{n-1}](\omega^{(n-1)}%
)>0\right\}  =\left\{  \omega^{(n-1)}:\ E_{Q}[I_{A}|\mathcal{G}_{n-1}%
](\omega^{(n-1)})>0\right\}  \text{.}%
\]
Given also equivalence of $P$ and $Q$ on $\mathcal{G}_{n-1}$, conclude that
\begin{align*}
Q\left(  \{E_{Q}[I_{A}|\mathcal{G}_{n-1}]>0\}\right)   &  =P\left(
\{E_{Q}[I_{A}|\mathcal{G}_{n-1}]>0\}\right) \\
&  =P\left(  \{E_{P}[I_{A}|\mathcal{G}_{n-1}]>0\}\right)  =0\text{,}%
\end{align*}
and hence $Q(A)=0$. \hfill\hfill$\blacksquare$  }

\subsection{{\protect\normalsize Some details for proof of Theorem
\ref{thm-special}\label{app-special}}}

{\normalsize Let {$\varphi\in C([-\infty,\infty])$} be symmetric with center
$c\in\mathbb{R}$ and decreasing on $(c,\infty)$. Define $\varphi_{h}$, for
$h>0$, by (\ref{phih}). Here we prove that: \newline(i) $\varphi_{h}$ is
symmetric with center $c$; and (ii) $sgn(\varphi_{h}^{\prime}(x))=-sgn(x-c)$.
}

{\normalsize \noindent\textbf{Proof:}\textbf{ } (i) By the definition of
$\varphi_{h}$,
\begin{align*}
\varphi_{h}(x+c)=  &  \int_{-\infty}^{\infty}\frac{1}{\sqrt{2\pi}}%
\varphi(x+c+hy)e^{-\tfrac{y^{2}}{2}}dy\\
=  &  \int_{-\infty}^{\infty}\frac{1}{\sqrt{2\pi}}\varphi(-x+c-hy)e^{-\tfrac
{y^{2}}{2}}dy\\
=  &  \int_{-\infty}^{\infty}\frac{1}{\sqrt{2\pi}}\varphi(-x+c+hy)e^{-\tfrac
{y^{2}}{2}}dy\\
=  &  \varphi_{h}(-x+c)
\end{align*}
(ii) Compute that%
\[
\varphi_{h}^{\prime}(x)=\int_{-\infty}^{\infty}\frac{1}{\sqrt{2\pi}h^{3}%
}\varphi(x+y)ye^{-\tfrac{y^{2}}{2h^{2}}}dy\text{.}%
\]
Since $\varphi_{h}$ is symmetric with $c$, we have for any $x>c$,{
\begin{align*}
\varphi_{h}^{\prime}(x)=  &  \int_{-\infty}^{\infty}\frac{1}{\sqrt{2\pi}h^{3}%
}\varphi(x+y)ye^{-\tfrac{y^{2}}{2h^{2}}}dy\\
=  &  \int_{0}^{\infty}\frac{1}{\sqrt{2\pi}h^{3}}\varphi(c+y+x-c)ye^{-\tfrac
{y^{2}}{2h^{2}}}dy\\
&  +\int_{-\infty}^{0}\frac{1}{\sqrt{2\pi}h^{3}}\varphi(c+y+x-c)ye^{-\tfrac
{y^{2}}{2h^{2}}}dy\\
=  &  \int_{0}^{\infty}\frac{1}{\sqrt{2\pi}h^{3}}\varphi(c+y+x-c)ye^{-\tfrac
{y^{2}}{2h^{2}}}dy\\
&  -\int_{0}^{\infty}\frac{1}{\sqrt{2\pi}h^{3}}\varphi(c+y+c-x)ye^{-\tfrac
{y^{2}}{2h^{2}}}dy\\
=  &  \int_{0}^{\infty}\frac{1}{\sqrt{2\pi}h^{3}}\left(  \varphi
(c+y+x-c)-\varphi(c+y+c-x)\right)  ye^{-\tfrac{y^{2}}{2h^{2}}}dy\\
<  &  0
\end{align*}
} Thus $sgn(\varphi_{h}^{\prime}(x))=-sgn(x-c)$.\hfill\hfill$\blacksquare$  }

\subsection{{\protect\normalsize Some details for hypothesis
testing\label{app-hyp}}}

{\normalsize Both $\left(  X_{i}\right)  $ and $\left(  Y_{i}\right)  $
described in section \ref{section-hyp} conform to the IID model, with the
common variance $\sigma^{2}$ and mean intervals $\left[  \underline{\mu
},\overline{\mu}\right]  $ and $\left[  \underline{\mu}-\theta,\overline{\mu
}-\theta\right]  $ respectively. (The text considers the special case $\left[
\underline{\mu}-\theta,\overline{\mu}-\theta\right]  =\left[  -\kappa
,\kappa\right]  $.) Let $\mu_{m}^{n}$ be defined by the form of {(\ref{mumn})}
appropriate for $\left(  X_{i}\right)  $ and denote by $\gamma_{m}^{n}$ the
corresponding variables appropriate for $\left(  Y_{i}\right)  $. Here we
prove (\ref{Mtheta}), for which it suffices to show that
\begin{equation}
\mu_{m}^{n}=\theta+\gamma_{m}^{n}\text{.} \label{mu-gamma}%
\end{equation}
}

{\normalsize By Theorem \ref{thm-special}(1), if $\varphi$ is decreasing on
$(c,+\infty)$, then
\[
\lim\limits_{n\rightarrow\infty}\sup_{Q\in\mathcal{P}}E_{Q}\left[
\varphi\left(  \frac{1}{n}\sum\limits_{i=1}^{n}X_{i}+\frac{1}{\sqrt{n}}%
\sum\limits_{i=1}^{n}\frac{1}{\sigma}\left(  X_{i}-\mu_{i}^{n}\right)
\right)  \right]  =\mathbb{E}_{\left[  \underline{\mu},\overline{\mu}\right]
}[\varphi\left(  B_{1}\right)  ],
\]
where, by (\ref{mumn}), $\mu_{m}^{n}=\overline{\mu}I_{A_{m-1,n}}%
+\underline{\mu}I_{A_{m-1,n}^{c}}$ \ and%
\[
A_{m-1,n}=\left\{  \frac{1}{n}\sum\limits_{i=1}^{m-1}X_{i}+\frac{1}{\sqrt{n}%
}\sum\limits_{i=1}^{m-1}\frac{1}{\sigma}\left(  X_{i}-\mu_{i}^{n}\right)
\leq-\frac{\overline{\mu}+\underline{\mu}}{2}\left(  {1-\frac{m-1}{n}}\right)
+c\right\}  .
\]
}

{\normalsize
Let $\phi(x)=\varphi(x+\theta)$. Then $\phi$ is symmetric with center
$\widehat{c}=c-\theta$.
Theorem \ref{thm-special}(1) applied to $(Y_{i})$ yields%
\[
\lim\limits_{n\rightarrow\infty}\sup_{Q\in\mathcal{P}}E_{Q}\left[  \phi\left(
\frac{1}{n}\sum\limits_{i=1}^{n}Y_{i}+\frac{1}{\sqrt{n}}\sum\limits_{i=1}%
^{n}\frac{1}{\sigma}\left(  Y_{i}-\gamma_{i}^{n}\right)  \right)  \right]
=\mathbb{E}_{\left[  \underline{\mu}-\theta,\overline{\mu}-\theta\right]
}[\phi\left(  B_{1}\right)  ]
\]
where $\gamma_{m}^{n}=(\overline{\mu}-\theta)I_{\widehat{A}_{m-1,n}%
}+(\underline{\mu}-\theta)I_{\widehat{A}_{m-1,n}^{c}}$ for $m=1,\cdots,n$,
and, for $m\geq1$,
\[
\widehat{A}_{m-1,n}=\left\{  \tfrac{1}{n}\sum\limits_{i=1}^{m-1}Y_{i}%
+\tfrac{1}{\sqrt{n}}\sum\limits_{i=1}^{m-1}\tfrac{1}{\sigma}\left(
Y_{i}-\gamma_{i}^{n}\right)  \leq-\left(  \tfrac{\overline{\mu}+\underline
{\mu}}{2}-\theta\right)  \left(  1-\tfrac{m-1}{n}\right)  +\widehat
{c}\right\}  .
\]
Replace $Y_{i}$ by $X_{i}-\theta$ to obtain%
\begin{align*}
\widehat{A}_{m-1,n}  &  = \left\{  \tfrac{1}{n}\sum\limits_{i=1}^{m-1}%
Y_{i}+\tfrac{1}{\sqrt{n}}\sum\limits_{i=1}^{m-1}\tfrac{1}{\sigma}\left(
Y_{i}-\gamma_{i}^{n}\right)  \leq-\left(  \tfrac{\overline{\mu}+\underline
{\mu}}{2}-\theta\right)  \left(  1-\tfrac{m-1}{n}\right)  +\widehat{c}\right\}
\\
&  = \left\{  \tfrac{1}{n}\sum\limits_{i=1}^{m-1}X_{i}+\tfrac{1}{\sqrt{n}}%
\sum\limits_{i=1}^{m-1}\tfrac{1}{\sigma}\left(  X_{i}-\gamma_{i}^{n}%
-\theta\right)  \leq- \tfrac{\overline{\mu}+\underline{\mu}}{2} \left(
1-\tfrac{m-1}{n}\right)  +c\right\}
\end{align*}
Thus $A_{0,n}=\widehat{A}_{0,n}$, and
\[
\gamma_{1}^{n}+\theta=\overline{\mu}I_{\widehat{A}_{0,n}}+\underline{\mu
}I_{\widehat{A}_{0,n}^{c}}=\overline{\mu}I_{A_{0,n}}+\underline{\mu}%
I_{A_{0,n}^{c}}=\mu_{1}^{n}%
\]
By induction, $A_{m-1,n}=\widehat{A}_{m-1,n}$, for $m\geq1$, and%
\[
\gamma_{m}^{n}+\theta=\overline{\mu}I_{\widehat{A}_{m-1,n}}+\underline{\mu
}I_{\widehat{A}_{m-1,n}^{c}}=\overline{\mu}I_{A_{m-1,n}}+\underline{\mu
}I_{A_{m-1,n}^{c}}=\mu_{m}^{n}\text{.}%
\]
}

\subsection{{\protect\normalsize Proof of Theorem \ref{thm-CLT2}%
\label{app-CLT2}}}

{\normalsize As noted previously (Remarks \ref{remark-CLT2} and
\ref{remark-CLT2b}), a suitably modified version of Lemma \ref{lemma-taylor}
is the key to proof of Theorem \ref{thm-CLT2}. Here we outline a proof of the
modified lemma. \ We prove it in two steps.  }

{\normalsize \medskip}

{\normalsize \noindent\noindent}

{\normalsize \noindent\textbf{Step 1}: For every $m\geq1$, {let}{
$\ \theta_{m}$ be a $\mathcal{G}_{m-1}$-measurable r.v. satisfying
\[
|\theta_{m}|\leq\kappa.
\]
} We prove that
\begin{equation}
\lim_{n\rightarrow\infty}\sum_{m=1}^{n}\left\vert \sup\limits_{Q\in
\mathcal{P}}E_{Q}\left[  H_{m,n}\left(  T_{m-1,n}+\frac{X_{m}-\theta_{m}%
}{\sigma\sqrt{n}}\right)  \right]  -\sup\limits_{Q\in\mathcal{P}}E_{Q}\left[
f(\theta_{m},m,n)\right]  \right\vert =0, \label{remainder-esti-without1}%
\end{equation}
where $f(\theta_{m},m,n)$ is given by: $f(\theta_{m},m,n)=${\footnotesize
\begin{equation}
H_{m,n}(T_{m-1,n})+H_{m,n}^{\prime}(T_{m-1,n})\left(  \frac{X_{m}-\theta_{m}%
}{\sigma\sqrt{n}}\right)  +\frac{1}{2}H_{m,n}^{\prime\prime}(T_{m-1,n})\left(
\frac{X_{m}-\theta_{m}}{\sigma\sqrt{n}}\right)  ^{2}. \label{f}%
\end{equation}
}
Let $x=T_{m-1,n}$ and $y=\frac{X_{m}-\theta_{m}}{\sigma\sqrt{n}}$ in
inequality (\ref{le0}), and obtain%
\begin{align*}
&  \sum_{m=1}^{n}\left\vert \sup\limits_{Q\in\mathcal{P}}E_{Q}\left[
H_{m,n}\left(  T_{m-1,n}+\frac{X_{m}-\theta_{m}}{\sigma\sqrt{n}}\right)
\right]  -\sup\limits_{Q\in\mathcal{P}}E_{Q}\left[  f(\theta_{m},m,n)\right]
\right\vert \\
\leq &  r_{1}(\overline{\epsilon},n)+r_{2}(C,n)\text{,}%
\end{align*}
where%
\begin{align*}
r_{1}(\overline{\epsilon},n)  &  :=\overline{\epsilon}\sum\limits_{m=1}%
^{n}\sup\limits_{Q\in\mathcal{P}}E_{Q}\left[  \left\vert \frac{X_{m}%
-\theta_{m}}{\sigma\sqrt{n}}\right\vert ^{2}I_{\left\{  |\frac{X_{m}%
-\theta_{m}}{\sigma\sqrt{n}}|<\delta\right\}  }\right] \\
r_{2}(C,n)  &  :=C\sum\limits_{m=1}^{n}\sup\limits_{Q\in\mathcal{P}}%
E_{Q}\left[  \left\vert \frac{X_{m}-\theta_{m}}{\sigma\sqrt{n}}\right\vert
^{2}I_{\left\{  |\frac{X_{m}-\theta_{m}}{\sigma\sqrt{n}}|\geq\delta\right\}
}\right]  .
\end{align*}
It is readily proven that, for sufficiently large $n$,%
\begin{align*}
r_{1}(\overline{\epsilon},n)\leq &  \frac{\overline{\epsilon}}{\sigma^{2}%
}\left(  \sigma^{2}+4\kappa^{2}\right) \\
r_{2}(C,n)\leq &  \frac{2C}{n\sigma^{2}}\sum\limits_{m=1}^{n}\sup
\limits_{Q\in\mathcal{P}}E_{Q}\left[  \left\vert X_{m}\right\vert
^{2}I_{\left\{  |\frac{X_{m}-\theta_{m}}{\sigma\sqrt{n}}|\geq\delta\right\}
}\right] \\
&  +\frac{2C}{\sigma^{2}n}\sum\limits_{m=1}^{n}\sup\limits_{Q\in\mathcal{P}%
}E_{Q}\left[  \left\vert \theta_{m}\right\vert ^{2}I_{\left\{  |\frac
{X_{m}-\theta_{m}}{\sigma\sqrt{n}}|\geq\delta\right\}  }\right] \\
\leq &  \frac{2C}{n\sigma^{2}}\sum\limits_{m=1}^{n}\sup\limits_{Q\in
\mathcal{P}}E_{Q}\left[  \left\vert X_{m}\right\vert ^{2}I_{\left\{
|X_{m}|>\sigma\sqrt{n}\delta-\kappa\right\}  }\right] \\
&  +\frac{2C}{\sigma^{2}n}\frac{\kappa^{2}}{\delta^{2}}\sum\limits_{m=1}%
^{n}\sup\limits_{Q\in\mathcal{P}}E_{Q}\left[  \left\vert \frac{X_{m}%
-\theta_{m}}{\sigma\sqrt{n}}\right\vert ^{2}\right]  .
\end{align*}
By the finiteness of $\kappa,\sigma$ and the Lindeberg condition
(\ref{linder}),%
\[
\lim\limits_{\overline{\epsilon}\rightarrow0}\lim\limits_{n\rightarrow\infty
}\left(  {r_{1}(\overline{\epsilon},n)}+r_{2}(C,n)\right)  =0\text{,}%
\]
which proves (\ref{remainder-esti-without1}).  }

{\normalsize \medskip}

{\normalsize \noindent\noindent}

{\normalsize \noindent\textbf{Step 2}: We take $\theta_{m}$$=E_{Q}%
[X_{m}|\mathcal{G}_{m-1}]$ in (\ref{f}). Then for all $n\geq m\geq1,$
\[
\sup\limits_{Q\in\mathcal{P}}E_{Q}\left[  f(\theta_{m},m,n)\right]
=\sup\limits_{Q\in\mathcal{P}}E_{Q}\left[  H_{m,n}(T_{m-1,n})+\tfrac{1}%
{2n}H_{m,n}^{\prime\prime}(T_{m-1,n})\right]  .
\]
In fact,
\begin{align*}
&  \sup\limits_{Q\in\mathcal{P}}E_{Q}\left[  f(\theta_{m},m,n)\right] \\
=  &  \sup\limits_{Q\in\mathcal{P}}E_{Q}\left[  H_{m,n}(T_{m-1,n}%
)+H_{m,n}^{\prime}(T_{m-1,n})\left(  \frac{X_{m}-E_{Q}[X_{m}|\mathcal{G}%
_{m-1}]}{\sigma\sqrt{n}}\right)  \right. \\
&  \qquad\qquad\qquad\qquad\quad\;\left.  +\frac{1}{2}H_{m,n}^{\prime\prime
}(T_{m-1,n})\left(  \frac{X_{m}-E_{Q}[X_{m}|\mathcal{G}_{m-1}]}{\sigma\sqrt
{n}}\right)  ^{2}\right] \\
=  &  \sup\limits_{Q\in\mathcal{P}}E_{Q}\left[  H_{m,n}(T_{m-1,n}%
)+H_{m,n}^{\prime}(T_{m-1,n})E_{Q}\left[  \bigg(\frac{X_{m}-E_{Q}%
[X_{m}|\mathcal{G}_{m-1}]}{\sigma\sqrt{n}}\bigg)|\mathcal{G}_{m-1}\right]
\right. \\
&  \qquad\qquad\qquad\quad\ \left.  +\frac{1}{2n}H_{m,n}^{\prime\prime
}(T_{m-1,n})E_{Q}\left[  \left(  \frac{X_{m}-E_{Q}[X_{m}|\mathcal{G}_{m-1}%
]}{\sigma\sqrt{n}}\right)  ^{2}|\mathcal{G}_{m-1}\right]  \right] \\
=  &  \sup\limits_{Q\in\mathcal{P}}E_{Q}\left[  H_{m,n}(T_{m-1,n})+\frac
{1}{2n}H_{m,n}^{\prime\prime}(T_{m-1,n})\right] \\
=  &  \sup\limits_{Q\in\mathcal{P}}E_{Q}\left[  L_{m,n}(T_{m-1,n})\right]
\text{.}%
\end{align*}
The last equality follows when $g_{0}$ in (\ref{eq-1}) equals 0.  }

\subsection{{\protect\normalsize Proof of LLN: Corollary \ref{cor-LLN}%
\label{app-LLN}}}

{\normalsize Here we prove Corollary \ref{cor-LLN}, showing how it can be
derived from our main result Theorem \ref{thm-CLT}, or more precisely, from
the following slight generalization.  }

\begin{thm}
{\normalsize \label{thm-CLTgen} Adopt the assumptions of Theorem
\ref{thm-CLT}. Then, for any $\varphi\in C\left(  \left[  -\infty
,\infty\right]  \right)  $, $\beta\geq0$ and $\alpha>0$,%
\begin{equation}
\lim\limits_{n\rightarrow\infty}\sup_{Q\in\mathcal{P}}E_{Q}\left[
\varphi\left(  \frac{\beta}{n}{\sum_{i=1}^{n}X_{i}}+\frac{\alpha}{\sqrt{n}%
}\sum\limits_{i=1}^{n}\frac{1}{\sigma}{(X_{i}-E_{Q}[X_{i}|\mathcal{G}_{i-1}%
])}\right)  \right]  =\mathbb{E}_{g}[\varphi\left(  \alpha B_{1}\right)  ],
\label{CLTgen}%
\end{equation}
where the right side of this equation is defined to be $Y_{0}$, given that
$(Y_{t},Z_{t})$ is the solution of the BSDE
\begin{equation}
Y_{t}=\varphi\left(  \alpha B_{1}\right)  +\int_{t}^{1}g(Z_{s})ds-\int_{t}%
^{1}Z_{s}dB_{s},\;0\leq t\leq1, \label{BSDEgen}%
\end{equation}
Here {$g(z):=\frac{\beta}{\alpha}\max\limits_{\underline{\mu}\leq\mu
\leq\overline{\mu}}(\mu z),$} and $(B_{t})$ is a standard Brownian motion.  }
\end{thm}

{\normalsize \noindent\textbf{Proof:} Change variables to $\tilde{X_{i}}%
=\frac{\beta}{\alpha}X_{i}$ and let $\tilde{\varphi}(x)=\varphi(\alpha x)$.
Then \newline$\mathbb{E}[\tilde{X_{i}}]=\frac{\overline{\mu}\beta}{\alpha}$,
$\mathcal{E}$$[\tilde{X_{i}}]=\frac{\underline{\mu}\beta}{\alpha}$ and their
variance is $\left(  \frac{\beta\sigma}{\alpha}\right)  ^{2}$. Apply Theorem
\ref{thm-CLT} to obtain  }

{\normalsize
\begin{align*}
&  \sup_{Q\in\mathcal{P}}E_{Q}\left[  \varphi\left(  \frac{\beta}{n}%
{\sum_{i=1}^{n}X_{i}}+\frac{\alpha}{\sqrt{n}\sigma}{\sum_{i=1}^{n}(X_{i}%
-E_{Q}[X_{i}|\mathcal{G}_{i-1}])}\right)  \right] \\
&  =\lim\limits_{n\rightarrow\infty}\sup_{Q\in\mathcal{P}}E_{Q}\left[
\tilde{\varphi}\left(  \frac{1}{n}{\sum_{i=1}^{n}\tilde{X}_{i}}+\frac{\alpha
}{\sqrt{n}\beta\sigma}{\sum_{i=1}^{n}(\tilde{X}_{i}-E_{Q}[\tilde{X}%
_{i}|\mathcal{G}_{i-1}])}\right)  \right] \\
&  =\mathbb{E}_{g}[\tilde{\varphi}\left(  B_{1}\right)  ]=\mathbb{E}%
_{g}[\varphi\left(  \alpha B_{1}\right)  ]. \
\end{align*}
\hfill\hfill$\blacksquare$  }

{\normalsize \noindent\noindent}

{\normalsize \noindent\textbf{Proof of Corollary \ref{cor-LLN}}: It suffices
to take $\varphi\in C_{b}^{\infty}(\mathbb{R})$. Let
\[
S_{0}=0,\;S_{n}=\sum_{i=1}^{n}X_{i},~S_{n}^{Q}=\sum_{i=1}^{n}Y_{i}%
^{Q}\text{,~~}Y_{i}^{Q}=\frac{1}{\sigma}\left(  X_{i}-E_{Q}[X_{i}%
|\mathcal{G}_{i-1}]\right)  \text{.}%
\]
Then, {for any $Q\in\mathcal{P}$},
\[
E_{Q}[Y_{i}^{Q}|\mathcal{G}_{i-1}]=0\text{ and }E_{Q}[(Y_{i}^{Q}%
)^{2}|\mathcal{G}_{i-1}]=1\text{ for all }i\text{.}%
\]
}

{\normalsize \medskip}

{\normalsize \noindent\noindent}

{\normalsize \noindent\textbf{Step 1: }Prove that
\[
\lim\limits_{\alpha\rightarrow0}\lim\limits_{n\rightarrow\infty}\sup
_{Q\in\mathcal{P}}E_{Q}\left[  \varphi\left(  \frac{S_{n}}{n}+\alpha
\frac{S_{n}^{Q}}{\sqrt{n}}\right)  \right]  =\lim\limits_{n\rightarrow\infty
}\sup_{Q\in\mathcal{P}}E_{Q}\left[  \varphi\left(  \frac{S_{n}}{n}\right)
\right]  .
\]
Since $\varphi\in C_{b}^{\infty}(\mathbb{R}),$ $\varphi$ is uniformly
Lipschitz continuous, i.e. $\exists C>0$ such that \newline$|\varphi
(x+y)-\varphi(x)|\leq C|x-y|$ for $x,y\in$$\mathbb{R}$$.$ Thus
\begin{align*}
&  \left\vert \sup_{Q\in\mathcal{P}}E_{Q}\left[  \varphi\left(  \frac{S_{n}%
}{n}+\alpha\frac{S_{n}^{Q}}{\sqrt{n}}\right)  \right]  -\sup_{Q\in\mathcal{P}%
}E_{Q}\left[  \varphi\left(  \frac{S_{n}}{n}\right)  \right]  \right\vert \\
\leq &  \frac{\alpha C}{\sqrt{n}}\sup_{Q\in\mathcal{P}}E_{Q}\left[  \left\vert
S_{n}^{Q}\right\vert \right] \\
\leq &  \frac{\alpha C}{\sqrt{n}}\left(  \sup_{Q\in\mathcal{P}}E_{Q}\left[
\left(  S_{n-1}^{Q}+Y_{n}^{Q}\right)  ^{2}\right]  \right)  ^{\frac{1}{2}}\\
=  &  \frac{\alpha C}{\sqrt{n}}\left(  \sup_{Q\in\mathcal{P}}E_{Q}\left[
(S_{n-1}^{Q})^{2}+2S_{n-1}^{Q}Y_{n}^{Q}+(Y_{n}^{Q})^{2}\right]  \right)
^{\frac{1}{2}}\\
=  &  \frac{\alpha C}{\sqrt{n}}\left(  \sup_{Q\in\mathcal{P}}E_{Q}\left[
(S_{n-1}^{Q})^{2}+2S_{n-1}^{Q}E_{Q}[Y_{n}^{Q}|\mathcal{G}_{n-1}]+E_{Q}%
[(Y_{n}^{Q})^{2}|\mathcal{G}_{n-1}]\right]  \right)  ^{\frac{1}{2}}\\
=  &  \frac{\alpha C}{\sqrt{n}}\left(  \sup_{Q\in\mathcal{P}}E_{Q}\left[
(S_{n-1}^{Q})^{2}\right]  +1\right)  ^{\frac{1}{2}}=\cdots=\frac{\alpha
C}{\sqrt{n}}\left(  n\right)  ^{\frac{1}{2}}\\
=  &  \alpha C\rightarrow0\text{\ as }\alpha\rightarrow0.
\end{align*}
}

{\normalsize \noindent\textbf{Step 2: }Prove that if $\beta=1$ and
{$g(z)=\frac{1}{\alpha}\max\limits_{\underline{\mu}\leq x\leq\overline{\mu}%
}(xz),$} then
\begin{equation}
\lim\limits_{\alpha\rightarrow0}\mathbb{E}_{g}[\varphi\left(  \alpha
B_{1}\right)  ]=\sup\limits_{\underline{\mu}\leq x\leq\overline{\mu}}%
\varphi\left(  x\right)  . \label{Step2}%
\end{equation}
Let
\begin{align*}
&  \mathcal{P}\equiv\\
&  \left\{  Q^{v}:E_{P^{\ast}}[\frac{dQ^{v}}{dP^{\ast} }|\mathcal{F}%
_{1}]=e^{-\frac{1}{2}\int_{0}^{1}\frac{v_{s} ^{2}}{\alpha^{2}}ds+\int_{0}%
^{1}\frac{v_{s}}{\alpha}dB_{s}}, (v_{t})\text{ is }\mathcal{F}_{t}%
\text{-adapted and } v\in[\underline{\mu},\overline{\mu}]\right\}
\end{align*}
where $v\in[\underline{\mu},\overline{\mu}]$ is in the sense of $\inf_{0\leq
s\leq1}v_{s}\ge\underline{\mu}$ and $\sup_{0\leq s\leq1}v_{s}\leq\overline
{\mu}$ a.s..  }

{\normalsize By \cite[Theorem 2.2]{CE} or \cite[Lemma 3]{chenpeng},
\begin{align}
\mathbb{E}_{g}[\varphi\left(  \alpha B_{1}\right)  ]  &  =\sup\limits_{Q\in
\mathcal{P}}E_{Q}\left[  \varphi\left(  \alpha B_{1}\right)  \right]
\nonumber\\
&  =\sup\limits_{\underline{\mu}\leq v\leq\overline{\mu}}E_{Q^{v}}\left[
\varphi\left(  \alpha B_{1}\right)  \right] \nonumber\\
&  =\sup\limits_{\underline{\mu}\leq v\leq\overline{\mu}}E_{Q^{v}}\left[
\varphi\left(  \alpha\left(  B_{1}-\int_{0}^{1}\frac{v_{s}}{\alpha}ds+\int
_{0}^{1}\frac{v_{s}}{\alpha}ds\right)  \right)  \right] \nonumber\\
&  =\sup\limits_{\underline{\mu}\leq v\leq\overline{\mu}}E_{Q^{v}}\left[
\varphi\left(  \alpha B_{1}^{v}+\int_{0}^{1}v_{s}ds\right)  \right]  \text{,}
\label{Step2.05}%
\end{align}
where $B_{t}^{v}\equiv B_{t}-\int_{0}^{t}\frac{v_{s}}{\alpha}ds$ is the
Brownian motion under $Q^{v}$.  }

{\normalsize We now prove that
\begin{equation}
\lim\limits_{\alpha\rightarrow0}\left\vert \sup\limits_{\underline{\mu}\leq
v\leq\overline{\mu}}E_{Q^{v}}\left[  \varphi\left(  \alpha B_{1}^{v}+\int
_{0}^{1}v_{s}ds\right)  \right]  -\sup\limits_{\underline{\mu}\leq
v\leq\overline{\mu}}E_{Q^{v}}\left[  \varphi\left(  \int_{0}^{1}%
v_{s}ds\right)  \right]  \right\vert =0. \label{Step2.1}%
\end{equation}
Because $\varphi$ has Lipschitz constant C$>0$,
\begin{align*}
&  \left\vert \sup\limits_{\underline{\mu}\leq v\leq\overline{\mu}}E_{Q^{v}%
}\left[  \varphi\left(  \alpha B_{1}^{v}+\int_{0}^{1}v_{s}ds\right)  \right]
-\sup\limits_{\underline{\mu}\leq v\leq\overline{\mu}}E_{Q^{v}}\left[
\varphi\left(  \int_{0}^{1}v_{s}ds\right)  \right]  \right\vert \\
&  \leq\sup\limits_{\underline{\mu}\leq v\leq\overline{\mu}}E_{Q^{v}}\left[
\left\vert \varphi\left(  \alpha B_{1}^{v}+\int_{0}^{1}v_{s}ds\right)
-\varphi\left(  \int_{0}^{1}v_{s}ds\right)  \right\vert \right] \\
&  \leq C\sup\limits_{\underline{\mu}\leq v\leq\overline{\mu}}E_{Q^{v}}\left[
\left\vert \alpha B_{1}^{v}\right\vert \right]  \leq\alpha C\rightarrow
0\quad\text{as}\ \alpha\rightarrow0,
\end{align*}
because $(B_{t}^{v})$ is $Q^{v}$-Brownian motion and $E_{Q^{v}}\left[
\left\vert B_{1}^{v}\right\vert \right]  \leq1.$ This proves (\ref{Step2.1}).
}

{\normalsize We now prove (\ref{Step2}): For any $x\in\lbrack\underline{\mu
},\overline{\mu}]$, let $v_{s}=x,$ $s\in\lbrack0,1]$. Then
\begin{equation}
\sup\limits_{\underline{\mu}\leq v\leq\overline{\mu}}E_{Q^{v}}\left[
\varphi\left(  \int_{0}^{1}v_{s}ds\right)  \right]  \geq\sup
\limits_{\underline{\mu}\leq x\leq\overline{\mu}}E_{Q^{v}}\left[
\varphi\left(  x\right)  \right]  =\sup\limits_{\underline{\mu}\leq
x\leq\overline{\mu}}\varphi\left(  x\right)  . \label{ineq-1}%
\end{equation}
In addition, since $\underline{\mu}\leq\inf_{s\in[0,1]} v_{s} \leq\sup
_{s\in[0,1]} v_{s} \leq\overline{\mu}$ a.s., we have  }

{\normalsize
\[
\sup\limits_{\underline{\mu}\leq v\leq\overline{\mu}}\varphi\left(  \int
_{0}^{1}v_{s}(\omega)ds\right)  \leq\sup\limits_{\underline{\mu}\leq
x\leq\overline{\mu}}\varphi\left(  x\right)  , \text{ a.s. } .
\]
Therefore,
\begin{align*}
\sup\limits_{\underline{\mu}\leq v\leq\overline{\mu}}E_{Q^{v}}\left[
\varphi\left(  \int_{0}^{1}v_{s}ds\right)  \right]   &  \leq\sup
\limits_{\underline{\mu}\leq v\leq\overline{\mu}}E_{Q^{v}}\left[
\sup\limits_{\underline{\mu}\leq v\leq\overline{\mu}}\varphi\left(  \int
_{0}^{1}v_{s}ds\right)  \right] \\
&  \leq\sup\limits_{\underline{\mu}\leq x\leq\overline{\mu}}\varphi\left(
x\right)  \text{,}%
\end{align*}
which implies, given (\ref{ineq-1}), that
\[
\sup\limits_{\underline{\mu}\leq v\leq\overline{\mu}}E_{Q^{v}}\left[
\varphi\left(  \int_{0}^{1}v_{s}ds\right)  \right]  =\sup\limits_{\underline
{\mu}\leq x\leq\overline{\mu}}\varphi\left(  x\right)  .
\]
From (\ref{Step2.05}) and (\ref{Step2.1}), we have
\begin{align*}
&  \lim\limits_{\alpha\rightarrow0}\mathbb{E}_{g}[\varphi\left(  \alpha
B_{1}\right)  ]=\lim\limits_{\alpha\rightarrow0}\sup\limits_{\underline{\mu
}\leq v\leq\overline{\mu}}E_{Q^{v}}\left[  \varphi\left(  \alpha B_{1}%
^{v}+\int_{0}^{1}v_{s}ds\right)  \right] \\
=  &  \lim\limits_{\alpha\rightarrow0}\left\{  \sup\limits_{\underline{\mu
}\leq v\leq\overline{\mu}}E_{Q^{v}}\left[  \varphi\left(  \alpha B_{1}%
^{v}+\int_{0}^{1}v_{s}ds\right)  \right]  -\sup\limits_{\underline{\mu}\leq
v\leq\overline{\mu}}E_{Q^{v}}\left[  \varphi\left(  \int_{0}^{1}%
v_{s}ds\right)  \right]  \right\} \\
&  +\sup\limits_{\underline{\mu}\leq x\leq\overline{\mu}}\varphi\left(
x\right) \\
=  &  \sup\limits_{\underline{\mu}\leq x\leq\overline{\mu}}\varphi\left(
x\right)  .
\end{align*}
The proof of Step 2 is complete.  }

{\normalsize Finally, let $\alpha\rightarrow0$ on both sides of (\ref{CLTgen})
and apply Steps 1 and 2.\hfill\hfill$\blacksquare$  }

\end{document}